\documentclass[a4paper,11pt]{amsart}
\usepackage{amssymb,amsfonts,amsxtra,amscd}
\usepackage[all]{xy}
\usepackage{fullpage}
\theoremstyle{plain}
\newtheorem{theorem}{Theorem}[section]
\newtheorem{lemma}[theorem]{Lemma}
\newtheorem{cor}[theorem]{Corollary}
\newtheorem{prop}[theorem]{Proposition}
\theoremstyle{definition}
\newtheorem{defi}[theorem]{Definition}
\newtheorem{example}[theorem]{Example}
\theoremstyle{remark}
\newtheorem{rem}[theorem]{Remark}
\numberwithin{equation}{section}
\newcommand{\ci}{\ensuremath{C_\infty}}
\newcommand{\ai}{\ensuremath{A_\infty}}
\newcommand{\li}{\ensuremath{L_\infty}}
\newcommand{\csalg}[1]{\ensuremath{\widehat{S}\Sigma #1^*}}
\newcommand{\ctalg}[1]{\ensuremath{\widehat{T}\Sigma #1^*}}
\newcommand{\clalg}[1]{\ensuremath{\widehat{L}\Sigma #1^*}}
\newcommand{\col}[2][i]{\ensuremath{\prod_{#1=1}^\infty (\Sigma {#2}^*)^{\cotimes #1}}}
\newcommand{\deof}[2][]{\ensuremath{\Omega^1_\mathrm{#1}(#2)}}
\newcommand{\drof}[2][]{\ensuremath{DR^1_\mathrm{#1}(#2)}}
\newcommand{\drnof}[2][]{\ensuremath{\overline{DR}^1_\mathrm{#1}(#2)}}
\newcommand{\drzf}[2][]{\ensuremath{DR^0_\mathrm{#1}(#2)}}
\newcommand{\drtf}[2][]{\ensuremath{DR^2_\mathrm{#1}(#2)}}
\newcommand{\de}[2][]{\ensuremath{\Omega^\bullet_\mathrm{#1}(#2)}}
\newcommand{\dr}[2][]{\ensuremath{DR^\bullet_\mathrm{#1}(#2)}}
\newcommand{\clac}[3][\bullet]{\ensuremath{C_\mathrm{CE}^{#1}(#2,#3)}}
\newcommand{\hlac}[3][\bullet]{\ensuremath{H_\mathrm{CE}^{#1}(#2,#3)}}
\newcommand{\choch}[3][\bullet]{\ensuremath{C^{#1}_\mathrm{Hoch}(#2,#3)}}
\newcommand{\hhoch}[3][\bullet]{\ensuremath{H^{#1}_\mathrm{Hoch}(#2,#3)}}
\newcommand{\caq}[3][\bullet]{\ensuremath{C^{#1}_\mathrm{Harr}(#2,#3)}}
\newcommand{\haq}[3][\bullet]{\ensuremath{H^{#1}_\mathrm{Harr}(#2,#3)}}
\newcommand{\cbr}[2][\bullet]{\ensuremath{C^{#1}_\mathrm{Bar}(#2)}}
\newcommand{\hbr}[2][\bullet]{\ensuremath{H^{#1}_\mathrm{Bar}(#2)}}
\newcommand{\cclac}[2][\bullet]{\ensuremath{CC^{#1}_\mathrm{CE}(#2)}}
\newcommand{\hclac}[2][\bullet]{\ensuremath{HC^{#1}_\mathrm{CE}(#2)}}
\newcommand{\cchoch}[2][\bullet]{\ensuremath{CC^{#1}_\mathrm{Hoch}(#2)}}
\newcommand{\hchoch}[2][\bullet]{\ensuremath{HC^{#1}_\mathrm{Hoch}(#2)}}
\newcommand{\ccnhoch}[2][\bullet]{\ensuremath{\overline{CC}^{#1}_\mathrm{Hoch}(#2)}}
\newcommand{\hcnhoch}[2][\bullet]{\ensuremath{\overline{HC}^{#1}_\mathrm{Hoch}(#2)}}
\newcommand{\ccaq}[2][\bullet]{\ensuremath{CC^{#1}_\mathrm{Harr}(#2)}}
\newcommand{\hcaq}[2][\bullet]{\ensuremath{HC^{#1}_\mathrm{Harr}(#2)}}
\newcommand{\ccnaq}[2][\bullet]{\ensuremath{\overline{CC}^{#1}_\mathrm{Harr}(#2)}}
\newcommand{\hcnaq}[2][\bullet]{\ensuremath{\overline{HC}^{#1}_\mathrm{Harr}(#2)}}
\newcommand{\ctsygan}[2][\bullet\bullet]{\ensuremath{CC^{#1}_\mathrm{Tsygan}(#2)}}

\newcommand{\connes}[2][\bullet\bullet]{\ensuremath{CC^{#1}_\mathrm{Connes}(#2)}}
\newcommand{\nconnes}[2][\bullet\bullet]{\ensuremath{\overline{CC}^{#1}_\mathrm{Connes}(#2)}}

\newcommand{\echoch}[4][\bullet]{\ensuremath{C\langle #4 \rangle^{#1}_\mathrm{Hoch}(#2,#3)}}
\newcommand{\ehhoch}[4][\bullet]{\ensuremath{H\langle #4 \rangle^{#1}_\mathrm{Hoch}(#2,#3)}}
\newcommand{\ecchoch}[3][\bullet]{\ensuremath{CC\langle #3 \rangle^{#1}_\mathrm{Hoch}(#2)}}
\newcommand{\ehchoch}[3][\bullet]{\ensuremath{HC\langle #3 \rangle^{#1}_\mathrm{Hoch}(#2)}}
\newcommand{\eccnhoch}[3][\bullet]{\ensuremath{\overline{CC}\langle #3 \rangle^{#1}_\mathrm{Hoch}(#2)}}
\newcommand{\ehcnhoch}[3][\bullet]{\ensuremath{\overline{HC}\langle #3 \rangle^{#1}_\mathrm{Hoch}(#2)}}
\newcommand{\thiso}[1][]{\ensuremath{\Theta_\mathrm{#1}}}

\newcommand{\zetiso}[1][]{\ensuremath{\zeta_\mathrm{#1}}}
\newcommand{\cotimes}{\ensuremath{\hat{\otimes}}}
\newcommand{\gf}{\ensuremath{\mathbb{K}}}
\newcommand{\fpsa}{\ensuremath{\gf\langle\langle\tau,\boldsymbol{t}\rangle\rangle}}
\newcommand{\pvect}{\ensuremath{\mathcal{P}Mod_\gf}}
\newcommand{\inlim}[2]{\ensuremath{\varprojlim_{#1} #2}}
\newcommand{\dilim}[2]{\ensuremath{\varinjlim_{#1} #2}}
\newcommand{\innprod}[1]{\ensuremath{\langle -,- \rangle:#1 \otimes #1 \to \gf}}
\newcommand{\noproof}{\begin{flushright} \ensuremath{\square} \end{flushright}}

\DeclareMathOperator{\Hom}{Hom}
\DeclareMathOperator{\Der}{Der}
\DeclareMathOperator{\Coder}{Coder}
\DeclareMathOperator{\Aut}{Aut}

\DeclareMathOperator{\ord}{ord}
\DeclareMathOperator{\ad}{ad}
\DeclareMathOperator{\id}{id}

\begin{document}
\begin{abstract}
This paper builds a general framework in which to study cohomology theories of strongly homotopy algebras, namely $A_\infty, C_\infty$ and $L_\infty$-algebras. This framework is based on noncommutative geometry as expounded by Connes and Kontsevich. The developed machinery is then used to establish a general form of Hodge decomposition of Hochschild and cyclic cohomology of $C_\infty$-algebras. This generalizes and puts in a conceptual framework previous work by Loday and Gerstenhaber-Schack.
\end{abstract}
\title{Cohomology theories for homotopy algebras and noncommutative geometry}
\author{Alastair Hamilton \and Andrey Lazarev}
\maketitle
\tableofcontents

\section{Introduction}

An $A_\infty$-algebra is a generalisation of an associative algebra introduced in \cite{staha2} for the purposes of studying $H$-spaces. It was originally defined via a system of higher multiplication maps satisfying a series of complicated relations. One can similarly define $C_\infty$ and $L_\infty$-algebras as $\infty$-generalisations of commutative and Lie algebras respectively.

More recently, $L_\infty$ and $A_\infty$-algebras have found applications in mathematical physics, particularly in string field theory and the theory of topological $\Sigma$-models, cf. \cite{zwiebach} and \cite{azksch}. In addition, $\infty$-algebras with an invariant inner product were introduced in \cite{kontfd} and \cite{kontsg} and were shown to have a close relation with graph homology and therefore to the intersection theory on the moduli spaces of complex curves and invariants of differentiable manifolds. A short, informal introduction to graph homology is contained in \cite{voronov}, more substantial accounts are in \cite{vogtmann}, \cite{hamgraph} and \cite{ccs}.

In this paper we give a detailed analysis of the cohomology theories associated to $A_\infty, L_\infty$ and $C_\infty$-algebras. It is likely that much of our work could be extended to algebras over a Koszul quadratic operad, however we see no advantage in working in this more general context since our most interesting applications are concerned with $C_\infty$-algebras and do not generalise to other operads.

One of the important ideas that we wish to advocate is working with a geometric definition of an $\infty$-algebra, see Definition \ref{def_infstr}, as a formal supermanifold together with a homological vector field. This idea is not new, cf. for example \cite{azksch} or \cite{lazmod}, however we feel that it has not been used to its full potential. Apart from the obvious advantage of being concise, this definition allows one to engage the powerful apparatus of noncommutative differential geometry which could be quite beneficial as we hope to demonstrate.

Very roughly, noncommutative geometry considers noncommutative algebras as if they were algebras of functions on `noncommutative spaces'. Noncommutative geometry, as a branch of mathematics, appeared in the mid-eighties, mostly thanks to the fundamental works of Connes, see \cite{connes} and references therein. Nowadays the term `noncommutative geometry' has many different meanings; it is studied in the context of $K$-theory of $C^*$-algebras; in measure theory, theory of foliations, characteristic classes of manifolds.
Quantum groups, cf. \cite{Kassel} also form part of noncommutative geometry analogous to the theory of algebraic groups in the conventional algebraic geometry.

There are also attempts to build a general theory of noncommutative algebraic schemes, cf. \cite{konroz1, konroz2}, \cite{kapranov}. Furthermore, one can also consider `operadic' geometry, i.e. geometry where `algebras of functions' are algebras over operads (i.e. associative, Lie, Poisson, etc.) This point of view is contained in \cite{ginzsg}; it stems from Kontsevich's work \cite{kontsg}. This is also the point of view taken in the present paper, except we restrict ourselves with considering only three operads; namely those governing commutative, associative and Lie algebra. Because of that the term `operad' is not mentioned explicitly in the main text.
In addition, the noncommutative geometry considered in this paper is formal: our model for the noncommutative algebra of functions is an inverse limit of finite-dimensional nilpotent algebras. The most important special case (which is the only one needed for the purposes of study of homotopy algebras) is the case of a pro-free algebra, i.e. the formal completion of a free (commutative, associative or Lie) algebra.

Our main application is concerned with the Hodge decomposition of the Hochschild and cyclic Hochschild cohomology of a \ci-algebra and generalises the work of previous authors; cf. \cite{loday}, \cite{gersch} and \cite{natsch}. Our geometric approach allowed us to considerably simplify the combinatorics present in the cited sources while working in the more general context of $C_\infty$-algebras. We actually get something new even for the usual (strictly commutative) algebras, namely the Hodge decomposition of the Hochschild, cyclic Hochschild and bar cohomology in the \emph{nonunital} case.

The paper is organised as follows. Sections \ref{sec_ncgeom} and \ref{relations} introduce the basics of formal noncommutative geometry in the commutative, associative and Lie worlds. This is largely a detailed exposition of a part of Kontsevich's paper \cite{kontfd}. The language of topological modules and topological algebras is used here and throughout the paper and we collected the necessary facts in Appendix \ref{app_todual}. Section \ref{prerequisites} deals with the definitions and basic properties of $\infty$-algebras, emphasising the geometrical viewpoint. Here we consider $\mathbb{Z}$-graded $\infty$-algebras; however all the results carry over with obvious modifications to the $\mathbb{Z}/2$-graded framework used in the cited work of Kontsevich. In section \ref{minimal} we prove the analogue of Kadeishvili's minimality theorem for $C_\infty$-algebras.

In section \ref{sec_iachom} Chevalley-Eilenberg, Hochschild and Harrison cohomology theories for $\infty$-algebras are defined along with their cyclic counterparts. The Hodge decomposition for $C_\infty$-algebras is established in sections \ref{sec_hdgdec} and \ref{cyclichodge}.

{\it Acknowledgement.} We are grateful to Jim Stasheff for many useful discussions and comments and to
Martin Markl who drew our attention to the works \cite{markl1} and \cite{markl2}.

\subsection{Notation and conventions} \label{sec_notcon}

Throughout the paper our ground ring $\gf$ will be an evenly graded commutative ring containing the
field $\mathbb{Q}$. For most applications it will be enough to assume that $\gf$ is a field or the rational numbers, however, our approach is designed to accommodate \emph{deformation theory} (which is not considered in the present paper) and to do this we need to allow ground rings which are not necessarily fields. \gf-algebras and \gf-modules will simply be called algebras and modules. We will assume that all of our \gf-modules are obtained from $k$-vector spaces by extension of scalars where $k=\mathbb{Q}$ or any other subfield of $\gf$. All of our tensors will be taken over the ground ring {\gf} unless stated otherwise.

Given a graded module V we define the tensor algebra $TV$ by
\[ TV:=\gf \oplus V \oplus V^{\otimes 2} \oplus \ldots \oplus V^{\otimes n} \oplus \ldots. \]
We define the symmetric algebra $SV$ as the quotient of $TV$ by the relation $x \otimes y = (-1)^{|x||y|}y \otimes x$. Finally we define the free Lie algebra $LV$ as the Lie subalgebra of the commutator algebra $TV$ which consists of linear combinations of Lie monomials.

We will use $~\widehat{ \ }~$ to denote completion. Given a profinite graded module $V$ we can define the completed versions of $SV$, $TV$ and $LV$. For instance the completed tensor algebra would be
\[ \widehat{T}V :=\prod_{i=0}^{\infty} V^{\cotimes i} \]
where $\cotimes$ denotes the completed tensor product. The completed symmetric algebra $\widehat{S}V$ is then the quotient of $\widehat{T}V$ by the usual relations and $\widehat{L}V$ is the Lie subalgebra of $\widehat{T}V$ consisting of all convergent (possibly uncountably) infinite linear combinations of Lie monomials. Further details on formal objects and constructions can be found in Appendix \ref{app_todual}. We will assume that all the formal commutative and associative $\gf$-algebras considered in the main text are augmented.

Given a \emph{formal} graded commutative, associative or Lie algebra $X$, the Lie algebra consisting of all \emph{continuous} derivations
\[ \xi : X \to X \]
is denoted by $\Der(X)$. In order to emphasise our use of geometrical ideas in this paper we will use the term `vector field' synonymously with `continuous derivation'. The group consisting of all invertible \emph{continuous} (commutative, associative or Lie) algebra homomorphisms
\[ \phi: X \to X \]
will be denoted by $\Aut(X)$. Again, in order to emphasise the geometrical approach we will call an `invertible continuous homomorphism' of formal graded commutative, associative or Lie algebras a `diffeomorphism'.

On many occasions we will want to deal with commutative, associative and Lie algebras simultaneously. When we do so we will often abuse the terminology, for instance calling a bimodule a `module' or calling a Lie algebra an `algebra'.

Given a formal graded associative algebra $X$, the module of commutators $[X,X]$ is defined as the module consisting of all convergent (possibly uncountably) infinite linear combinations of elements of the form,
\[ [x,y]:=x \cdot y - (-1)^{|x||y|}y \cdot x; \quad x,y \in X \]
We will often denote the module of commutators by $[-,-]$ when the context makes it clear what $X$ is.

Given a profinite graded module $V$ we can place a grading on $\widehat{T}V$ which is different from the grading which is naturally inherited from $V$. We say an element $x \in \widehat{T}V$ has homogeneous \emph{order} $n$ if $x \in V^{\cotimes n}$. This also defines a grading by order on $\widehat{S}V$ and $\widehat{L}V$ since they are quotients and submodules respectively of $\widehat{T}V$.

We say a continuous endomorphism (linear map) $f:\widehat{T}V \to \widehat{T}V$ has homogeneous \emph{order} $n$ if it takes any element of order $i$ to an element of order $i+n-1$. Any vector field $\xi \in \Der(\widehat{T}V)$ could be written in the form
\begin{equation} \label{eqn_vecfrm}
\xi = \xi_0 + \xi_1 + \xi_2 + \ldots + \xi_n + \ldots,
\end{equation}
where $\xi_i$ is a vector field of order $i$. We say that $\xi$ \emph{vanishes at zero} if $\xi_0 = 0$.

Likewise any diffeomorphism $\phi \in \Aut(\widehat{T}V)$ could be written in the form
\[ \phi = \phi_1 + \phi_2 + \phi_3 + \ldots + \phi_n + \ldots, \]
where $\phi_i$ is an endomorphism of order $i$. We call $\phi$ a \emph{pointed diffeomorphism} if $\phi_1 = \id$. Similarly, we could make the same definitions and observations if we were to replace $\widehat{T}V$ with $\widehat{S}V$ or $\widehat{L}V$ in the above.

We will denote the symmetric group on $n$ letters by $S_n$ and the cyclic group of order $n$ by $Z_n$. Given a module $M$ over a group $G$ the module of coinvariants will be denoted by $M_G$ and the module of invariants by $M^G$.

Given a graded profinite module $V$, its completed tensor power $V^{\cotimes n}$ has a continuous action of the cyclic group $Z_n$ on it which is the restriction of the canonical action of $S_n$ to the subgroup $Z_n \cong \langle (n\,n-1\ldots 2\,1) \rangle \subset S_n$. If we define the algebra $\Lambda$ as
\[ \Lambda:= \prod_{n=1}^\infty\mathbb{Z}[Z_n] \]
then the action of each $Z_n$ on $V^{\cotimes n}$ for $n\geq 1$ gives $\prod_{n=1}^\infty V^{\cotimes n}$ the structure of a left $\Lambda$-module.

Let $z_n$ denote the generator of $Z_n$ corresponding to the cycle $(n\,n-1\ldots 2\,1)$. We define $z \in \Lambda$ by the formula,
\[ z:=\sum_{n=1}^\infty z_n. \]
Let the norm operator $N_n \in \mathbb{Z}[Z_n]$ be the element given by the formula,
\[ N_n:=1+z_n+z_n^2 +\ldots + z_n^{n-1}. \]
We define $N \in \Lambda$ by the formula,
\[ N:= \sum_{n=1}^\infty N_n. \]
The operators $z$ and $N$ will be used regularly throughout the paper.

Our convention will be to always work with cohomologically graded objects and we consequently define the suspension of a graded module $V$ as $\Sigma V$ where $\Sigma V^i:=V^{i+1}$. We define the desuspension of $V$ as $\Sigma^{-1}V$ where $\Sigma^{-1}V^{i}:= V^{i-1}$. The term `differential graded algebra' will be abbreviated as `DGA'.

Given a graded module $V$, we denote the graded \gf-linear dual
\[ \Hom_\gf(V,\gf) \]
by $V^*$. If $V$ is free with basis $\{v_\alpha, \alpha \in I\}$, then we denote the dual basis of $V^*$ by
\[ \{v_\alpha^*, \alpha \in I\}. \]
For the sake of clarity, when we write $\Sigma V^*$ we mean the graded module $\Hom_\gf(\Sigma V,\gf)$.
\section{Formal Noncommutative Geometry} \label{sec_ncgeom}

In this section we will collect the definitions and fundamental facts about formal commutative and noncommutative geometry. Although much of this work can be done in the generality of working with an algebra over an operad, our applications will be concerned with the three particular theories of (non)commutative geometry corresponding to commutative, associative and Lie algebras and we will describe each theory in its own detail.

The starting point is to define the (non)commutative 1-forms in each of the three settings. From this we can construct a differential envelope of our commutative, associative or Lie algebra. We can then construct the de Rham complex from the differential envelope and introduce the contraction and Lie operators. We conclude the section with some basic facts from noncommutative geometry. The theory we consider here will be the formal version of the theory considered by Kontsevich in \cite{kontsg}. The reader can refer to Appendix \ref{app_todual} for a discussion of the formal objects and constructions that we use.

\begin{defi} \label{def_nconfm}
\
\begin{enumerate}
\item[(a)]
Let $A$ be a formal graded commutative algebra. Consider the module $A \cotimes (A/\gf)$ and write $x \cotimes y$ as $x\cotimes dy$. The module of commutative 1-forms \deof[Com]{A} is defined as the quotient of $A \cotimes (A/\gf)$ by the relations
\[ x \cotimes d(yz)=(-1)^{|z|(|y|+|x|)}zx\cotimes dy + xy\cotimes dz. \]
This is a formal left $A$-module via the action
\[ a \cdot x\cotimes dy:=ax\cotimes dy. \]
\item[(b)]
Let $A$ be a formal graded associative algebra. The module of noncommutative 1-forms \deof[Ass]{A} is defined as
\[ \deof[Ass]{A}:= A \cotimes (A/\gf). \]
Let us write $x\cotimes y$ as $x\cotimes dy$. \deof[Ass]{A} has the structure of a formal $A$-bimodule via the actions
\begin{align*}
& a \cdot x\cotimes dy:= ax \cotimes dy, \\
& x \cotimes dy \cdot a:= x \cotimes d(ya) - xy \cotimes da.
\end{align*}
\item[(c)]
Let $g$ be a formal graded Lie algebra. Consider the module $\widehat{\mathcal{U}}(g) \cotimes g$ and write $x\cotimes y=x\cotimes dy$. Then the module of Lie 1-forms \deof[Lie]{g} is defined as the quotient of $\widehat{\mathcal{U}}(g) \cotimes g$ by the relations
\[ x \cotimes d[y,z]=xy \cotimes dz -(-1)^{|z||y|}xz \cotimes dy. \]
This is a formal left $g$-module via the action
\[ x \cdot (y \cotimes dz):= xy \cotimes dz. \]
\end{enumerate}
\end{defi}

\begin{rem}
When formulating theorems and proofs it will often be our policy to discuss these three theories simultaneously where this is practical. When we do this we will omit the subscripts $\mathrm{Com}$, $\mathrm{Ass}$ and $\mathrm{Lie}$. It should be understood that the reader should choose the appropriate subscript/construction depending on whether they wish to work with commutative, associative or Lie algebras.
\end{rem}

Let X be either a formal graded commutative, associative or Lie algebra. The following proposition shows that the module of 1-forms \deof{X} could be introduced as the unique $X$-module representing the functor $\Der^0(X,-)$ which sends a formal graded $X$-module $M$ to the module consisting of all continuous derivations $\xi:X \to M$ of degree zero. This proposition is the formal analogue of Lemma 5.5 of \cite{ginzsg}:

\begin{prop}
Let $X$ be either a formal graded commutative, associative or Lie algebra, then the map $d:X \to \Omega^1(X)$ is a derivation of degree zero. Given any formal graded $X$-module $M$ there is an isomorphism which is natural in both variables:
\begin{displaymath}
\begin{array}{ccc}
\Der^0(X,M) & \cong & \Hom^0_{X\mathrm{-module}}(\Omega^1(X),M), \\
\partial & \mapsto & [dx \mapsto \partial x]. \\
\end{array}
\end{displaymath}
\end{prop}
\noproof

This proposition could be summarised by the following diagram:
\[ \xymatrix{X \ar^{d}[rr] \ar_{\partial}[rd] && \Omega^1(X) \ar@{-->}^{\exists ! \phi}[ld] \\ & M} \]

\begin{rem}
Note that if $X$ is an associative algebra then the word `$X$-module' should be replaced with the word `$X$-bimodule'.
\end{rem}

Now we want to extend the 1-forms $\Omega^1(X)$ to the module of forms $\Omega^*(X)$ by forming a differential envelope of $X$.

\begin{defi} \label{def_ncfrms}
\
\begin{enumerate}
\item[(a)]
Let $A$ be a formal graded commutative algebra. The module of commutative forms \de[Com]{A} is defined as
\[ \de[Com]{A}:= \widehat{S}_A(\Sigma^{-1}\deof[Com]{A})=A \times \prod_{i=1}^\infty (\underbrace{\Sigma^{-1}\deof[Com]{A}\underset{A}{\cotimes}\ldots\underset{A}{\cotimes}\Sigma^{-1}\deof[Com]{A}}_{i \text{ factors}})_{S_i}. \]
Since \deof[Com]{A} is a module over the commutative algebra $A$, \de[Com]{A} has the structure of a formal graded commutative algebra whose multiplication is the standard multiplication on the completed symmetric algebra $\widehat{S}_A(\Sigma^{-1}\deof[Com]{A})$. The map $d:A \to \deof[Com]{A}$ lifts uniquely to give \de[Com]{A} the structure of a formal DGA.
\item[(b)]
Let $A$ be a formal graded associative algebra. The module of noncommutative forms \de[Ass]{A} is defined as
\[ \de[Ass]{A}:=\widehat{T}_A(\Sigma^{-1}\deof[Ass]{A})=A \times \prod_{i=1}^\infty \underbrace{\Sigma^{-1}\deof[Ass]{A} \underset{A}{\cotimes} \ldots \underset{A}{\cotimes} \Sigma^{-1}\deof[Ass]{A}}_{i \text{ factors}}. \]
Since \deof[Ass]{A} is an $A$-bimodule, \de[Ass]{A} has the structure of a formal associative algebra whose multiplication is the standard associative multiplication on the tensor algebra $\widehat{T}_A(\Sigma^{-1}\deof[Ass]{A})$. The map $d:A \to \deof[Ass]{A}$ lifts uniquely to give \de[Ass]{A} the structure of a formal DGA.
\item[(c)]
Let $g$ be a formal graded Lie algebra. The module of Lie forms \de[Lie]{g} is defined as
\[ \de[Lie]{g}:=g\ltimes \widehat{L}(\Sigma^{-1}\deof[Lie]{g}) \]
where the action of $g$ on $\widehat{L}(\Sigma^{-1}\deof[Lie]{g})$ is the restriction of the standard action of $g$ on $\widehat{T}(\Sigma^{-1}\deof[Lie]{g})$ to the Lie subalgebra of Lie monomials. The map $d:g \to \deof[Lie]{g}$ lifts uniquely to give \de[Lie]{g} the structure of a formal DGLA.
\end{enumerate}
\end{defi}

\begin{rem}
Since \deof{X} is a formal $X$-module, $\Sigma^{-1}\deof{X}$ is a formal $X$-module via the action
\[ x \cdot \Sigma^{-1}y:=(-1)^{|x|}\Sigma^{-1}(x \cdot y); \quad x \in X, y \in \deof{X}. \]
The additional signs appear because of the Koszul sign rule. This results in the map $d:\de{X} \to \de{X}$ being a graded derivation of degree one.
\end{rem}

\begin{rem}
The module of forms \de{X} inherits a grading from the graded modules $X$ and \deof{X}. In fact it has a bigrading. An element $x \in X$ has bidegree $(0,|x|)$ and an element $x \cotimes dy \in \deof{X}$ has bidegree $(1,|x|+|y|)$. This implicitly defines a bigrading on the whole of \de{X}. The map $d:\de{X} \to \de{X}$ has bidegree $(1,0)$ in this bigrading. The natural grading on \de{X} inherited from the graded modules $X$ and \deof{X} coincides with the total grading of this bigrading.
\end{rem}

The following proposition says that \de{X} could be uniquely defined as the differential envelope of $X$. This proposition is proved in the operadic and nonformal context in Proposition 5.6 of \cite{ginzsg}.

\begin{prop} \label{prop_difenv}
Let $X$ be either a formal graded commutative, associative or Lie algebra and let $M$ be either a formal differential graded commutative, associative or Lie algebra respectively. There is a natural adjunction isomorphism:
\[ \Hom_\mathrm{Alg}(X,M) \cong \Hom_\mathrm{DGA}(\de{X},M). \]
The adjunction isomorphism is given by the following diagram:
\begin{displaymath}
\xymatrix{X \ar@<-0.25ex>@{^{(}->}[rr] \ar^{\phi}[rd] && \de{X} \ar@{-->}^{\exists !\psi}[ld] \\ & M}
\end{displaymath}
\end{prop}
\noproof

\begin{rem} \label{rem_defunc}
It follows from Proposition \ref{prop_difenv} that $\Omega^\bullet$ is functorial. Suppose that both $X$ and $Y$ are either formal graded commutative, associative or Lie algebras. By Proposition \ref{prop_difenv}, any continuous algebra homomorphism $\phi:X \to Y$ lifts uniquely to a continuous differential graded algebra homomorphism $\phi^*:\de{X} \to \de{Y}$.
\end{rem}

Next we want to introduce contraction and Lie operators onto the module of forms \de{X}.

\begin{defi} \label{def_ncoper}
Let $X$ be either a formal graded commutative, associative or Lie algebra and let $\xi:X \to X$ be a vector field:
\begin{enumerate}
\item[(i)]
We can define a vector field $L_\xi: \de{X} \to \de{X}$ of bidegree $(0,|\xi|)$, called the Lie derivative, by the formula;
\begin{displaymath}
\begin{array}{lc}
L_\xi(x):=\xi(x), & x \in X; \\
L_\xi(dx):=(-1)^{|\xi|}d(\xi(x)), & x \in X. \\
\end{array}
\end{displaymath}
\item[(ii)]
We can define a vector field $i_\xi:\de{X} \to \de{X}$ of bidegree $(-1,|\xi|)$, called the contraction operator, by the formula;
\begin{displaymath}
\begin{array}{lc}
i_\xi(x):=0, & x \in X; \\
i_\xi(dx):=\xi(x), & x \in X. \\
\end{array}
\end{displaymath}
\end{enumerate}
\end{defi}

These operators satisfy certain important identities which are summarised by the following lemma.

\begin{lemma} \label{lem_schf}
Let both $X$ and $Y$ be either formal graded commutative, associative or Lie algebras, let $\xi:X \to X$ and $\gamma:X \to X$ be vector fields and let $\phi:X \to Y$ be a diffeomorphism, then we have the following identities:
\begin{displaymath}
\begin{array}{rl}
\textnormal{(i)} & L_\xi=[i_\xi,d]. \\
\textnormal{(ii)} & [L_\xi,i_\gamma]=i_{[\xi,\gamma]}. \\
\textnormal{(iii)} & L_{[\xi,\gamma]}=[L_\xi,L_\gamma]. \\
\textnormal{(iv)} & [i_\xi,i_\gamma]=0. \\
\textnormal{(v)} & [L_\xi,d]=0. \\
\textnormal{(vi)} & L_{\phi\xi\phi^{-1}}=\phi^* L_\xi \phi^{*-1}. \\
\textnormal{(vii)} & i_{\phi\xi\phi^{-1}}=\phi^* i_\xi \phi^{*-1}. \\
\end{array}
\end{displaymath}
\end{lemma}

\begin{proof}
Since \de{X} is generated by elements in $X$ and $d(X)$, we only need to establish the identities on these generators. Let us prove (i), the  Cartan homotopy formula. It should then be clear to the reader how to obtain the other identities.
\begin{align*}
& i_\xi d(x)=\xi(x), \quad di_\xi(x)=0; \\
& L_\xi(x)=\xi(x)=[i_\xi,d](x).
\end{align*}
Similarly,
\begin{align*}
& i_\xi d(dx)=0, \quad di_\xi(dx)=d(\xi(x)); \\
& L_\xi(dx)=(-1)^{|\xi|}d(\xi(x))=[i_\xi,d](dx).
\end{align*}
\end{proof}

We are now in a position to define the de Rham complex which will be the fundamental object in our applications of formal noncommutative geometry. For the purposes of working with commutative, associative or Lie algebras, one could use the definitions introduced by Kontsevich in \cite{kontsg}. These definitions were extended to the framework of operads by Getzler and Kapranov in \cite{getkap}. The functor $F: \mathcal{F}Alg \to \pvect$ (see Appendix \ref{app_todual} for definitions of the categories) is defined by the formula;
\[ F(X):=\widehat{S}^{\,2}(X)/(x\cotimes\mu(y,z)=\mu(x,y)\cotimes z), \]
where $X$ is a formal commutative, associative or Lie algebra and $\mu$ is the multiplication or Lie bracket respectively. This functor was originally introduced by Kontsevich in \cite{kontsg}. The de Rham complex is then defined as the result of applying the functor $F$ to the differential envelope \de{X}. The differential $d:\de{X} \to \de{X}$ induces the differential on the de Rham complex by applying the Leibniz rule, although in general the algebra structure will be lost and we will just have a complex.

\begin{defi} \label{def_derham}
\
\begin{enumerate}
\item[(a)]
Let $A$ be a formal graded commutative algebra. The de Rham complex \dr[Com]{A} is defined as
\[ \dr[Com]{A}:=\de[Com]{A}. \]
The differential is the same as the differential defined in Definition \ref{def_ncfrms}.
\item[(b)]
Let $A$ be a formal graded associative algebra. The de Rham complex \dr[Ass]{A} is defined as
\[ \dr[Ass]{A}:=\de[Ass]{A}/[\de[Ass]{A},\de[Ass]{A}]. \]
The differential is induced by the differential defined in Definition \ref{def_ncfrms}.
\item[(c)]
Let $g$ be a formal graded Lie algebra. The de Rham complex \dr[Lie]{g} is defined as the quotient of $\de[Lie]{g} \cotimes \de[Lie]{g}$ by the relations
\begin{displaymath}
\begin{array}{ll}
x \cotimes y = (-1)^{|x||y|}y \cotimes x; & x,y \in \de[Lie]{g}; \\
x \cotimes [y,z] = [x,y] \cotimes z; & x,y,z \in \de[Lie]{g}. \\
\end{array}
\end{displaymath}
The differential $d:\dr[Lie]{g} \to \dr[Lie]{g}$ is induced by the differential $d:\de[Lie]{g} \to \de[Lie]{g}$ by specifying
\[ d(x \cotimes y):=dx \cotimes y + (-1)^{|x|}x \cotimes dy; \quad x,y \in \de[Lie]{g}. \]
\end{enumerate}
\end{defi}

\begin{rem}
In the commutative and associative cases the definition in terms of the functor $F$ can be simplified, resulting in the definitions given above. The identification of $F(\de{X})$ with \dr{X} is given by
\begin{displaymath}
\begin{array}{ccc}
F(\de{X}) & \cong & \dr{X}, \\
x \cotimes y & \mapsto & x \cdot y. \\
\end{array}
\end{displaymath}
In the commutative case we see that \de[Com]{X} and \dr[Com]{X} actually coincide. In the Lie case we must resort to the slightly awkward definition of the de Rham complex given by the functor $F$.
\end{rem}

\begin{rem} \label{rem_ordgra}
Suppose we are given a profinite graded module $W$ and consider the commutative, associative or Lie algebras $\widehat{S}W$, $\widehat{T}W$ and $\widehat{L}W$. Recall from section \ref{sec_notcon} that these modules have a grading on them which we called the grading by order. An element $x \in \widehat{T}W$ has homogeneous order $n$ if $x \in W^{\cotimes n}$. We extend the grading by order on $\widehat{T}W$ to the whole of \dr[Ass]{\widehat{T}W} by stipulating that a 1-form $xdy$ has order $\ord(x)+\ord(y)-1$. This suffices to completely determine the grading by order on \dr[Ass]{\widehat{T}W}. For instance an $n$-form
\[ x_0 dx_1 \ldots dx_n \]
has order $\ord(x_0) + \ldots + \ord(x_n) - n$. We obtain a grading by order on \dr[Com]{\widehat{S}W} and \dr[Lie]{\widehat{L}W} in a similar and obvious manner.
\end{rem}

\begin{rem} \label{rem_drfunc}
It follows from Remark \ref{rem_defunc} that the de Rham complex construction is functorial. Suppose that both $X$ and $Y$ are either formal graded commutative, associative or Lie algebras and that $\phi:X \to Y$ is a continuous algebra homomorphism. The continuous differential graded algebra homomorphism $\phi^*:\de{X} \to \de{Y}$ induces a map between the de Rham complexes which, by an abuse of notation, we also denote by $\phi^*$;
\begin{equation} \label{eqn_drfunc}
\phi^*: \dr{X} \to \dr{Y}.
\end{equation}
For instance, in the Lie case the map $\phi^*:\dr[Lie]{g} \to \dr[Lie]{h}$ is given by;
\[ \phi^*(x \cotimes y):= \phi^*(x) \cotimes \phi^*(y); \quad x,y \in \de[Lie]{g}. \]
\end{rem}

The de Rham complex inherits a bigrading derived from the bigrading on \de{X}. Given a vector field $\xi:X \to X$ we can define Lie and contraction operators $L_\xi,i_\xi:\dr{X} \to \dr{X}$ of bidegrees $(0,|\xi|)$ and $(-1,|\xi|)$ respectively. These are induced by the Lie and contraction operators defined in Definition \ref{def_ncoper} and by an abuse of notation, we denote them by the same letters. It should be clear how $L_\xi$ and $i_\xi$ are defined except perhaps in the Lie case. In this case $L_\xi:\dr[Lie]{g} \to \dr[Lie]{g}$ is defined by specifying
\[ L_\xi(x \cotimes y):= L_\xi(x) \cotimes y + (-1)^{|x|}x \cotimes L_\xi(y); \quad x,y \in \de[Lie]{g}. \]
The operator $i_\xi:\dr[Lie]{g} \to \dr[Lie]{g}$ is defined in the same way. It is obvious that the operators $L_\xi,i_\xi:\dr{X} \to \dr{X}$ satisfy all the identities of Lemma \ref{lem_schf}.

We denote the $i$th fold suspension of the component of \dr{X} of bidegree $(i,\bullet)$ by $DR^i(X)$ and call these the $i$-forms. This definition gives us the identity;
\[ \dr{X} = DR^0(X) \times \Sigma^{-1}DR^1(X) \times \Sigma^{-2}DR^2(X) \times \ldots \times \Sigma^{-n}DR^n(X) \times \ldots. \]
The following lemma describing $DR^1(X)$ when $X$ is pro-free could be found in \cite{ginzsg}. There is also a description of $DR^0(X)$ in Proposition 4.9 of \cite{getkap} and in \cite{ginzsg}.

\begin{lemma} \label{lem_ofisom}
Let $W$ be a profinite graded module. We have the following isomorphisms of graded modules:
\begin{displaymath}
\begin{array}{rc}
\textnormal{(a)} & \thiso[Com] : W \cotimes \widehat{S}W \to \drof[Com]{\widehat{S}W}, \\
& \thiso[Com](x \cotimes y) := dx \cdot y. \\
\textnormal{(b)} & \thiso[Ass] : W \cotimes \widehat{T}W \to \drof[Ass]{\widehat{T}W}, \\
& \thiso[Ass](x \cotimes y) := dx \cdot y. \\
\textnormal{(c)} & \thiso[Lie] : W \cotimes \widehat{L}W \to \drof[Lie]{\widehat{L}W}, \\
& \thiso[Lie](x \cotimes y) := dx \cotimes y. \\
\end{array}
\end{displaymath}
\end{lemma}

\begin{proof}
Let us prove (c) since it is the most abstract. It should then be clear how to prove the other maps are isomorphisms. The 1-forms \drof[Lie]{\widehat{L}W} are generated by representatives of the form $x \cotimes y$ where $x \in \deof[Lie]{\widehat{L}W}$ and $y \in \widehat{L}W$. Using the Leibniz rule
\[ d[a,b]= (-1)^{|a|}[a,db] -(-1)^{(|a|+1)|b|}[b,da]; \quad a,b \in \widehat{L}W \]
we conclude that \deof[Lie]{\widehat{L}W} is generated by elements of the form $u \cdot dw\cotimes y$ where $u \in \widehat{\mathcal{U}}(\widehat{L}W)$, $w \in W$ and $y \in \widehat{L}W$. Now using the relation we imposed on \dr[Lie]{\widehat{L}W} in Definition \ref{def_derham},
\[ [a,b] \cotimes c = a \cotimes [b,c]; \quad a,b,c \in \de[Lie]{\widehat{L}W}; \]
we assert that \drof[Lie]{\widehat{L}W} is generated by representatives of the form $dw \cotimes y$ where $w \in W$ and $y \in \widehat{L}W$. It is now clear that \thiso[Lie] is an isomorphism.
\end{proof}

\begin{rem}
Note that since $W \cotimes \widehat{T}W = \prod_{i=1}^\infty W^{\cotimes i}$ we may also regard \thiso[Ass] as a map
\[ \thiso[Ass] : \prod_{i=1}^\infty W^{\cotimes i} \to \drof[Ass]{\widehat{T}W}. \]
\end{rem}

Given either a formal commutative, associative or Lie algebra $X$, we define the de Rham cohomology of $X$ to be the cohomology of the complex \dr{X}. We define $H^i(\dr{X})$ as
\[ H^i(\dr{X}):=\frac{\{ x \in DR^i(X):dx=0 \}}{d(DR^{i-1}(X))}. \]
The following result is the (non)commutative analogue of the Poincar\'e lemma and can be found in \cite{kontsg}.

\begin{lemma} \label{lem_poinca}
Let $W$ be a profinite graded module:
\begin{enumerate}
\item[(a)]
\begin{displaymath}
\begin{array}{lcll}
H^i(\dr[Com]{\widehat{S}W}) & = & 0, & i \geq 1; \\
H^0(\dr[Com]{\widehat{S}W}) & = & \gf. \\
\end{array}
\end{displaymath}
\item[(b)]
\begin{displaymath}
\begin{array}{lcll}
H^i(\dr[Ass]{\widehat{T}W}) & = & 0, & i \geq 1; \\
H^0(\dr[Ass]{\widehat{T}W}) & = & \gf. \\
\end{array}
\end{displaymath}
\item[(c)]
\[ H^\bullet(\dr[Lie]{\widehat{L}W}) = 0. \]
\end{enumerate}
\end{lemma}

\begin{proof}
Let us prove (b). It will then be clear that the same argument will work to establish the other cases. Choose a topological basis $\{x_i\}_{i \in I}$ for the profinite module $W$ and consider Euler's vector field $\xi:\widehat{T}W \to \widehat{T}W$ of degree zero and order one:
\[ \xi:= \sum_{i \in I}x_i\partial_{x_i}. \]
Since we are working in characteristic zero, $L_\xi:\dr[Ass]{\widehat{T}W} \to \dr[Ass]{\widehat{T}W}$ establishes a one-to-one correspondence on every component except on the component
\[ \gf \subset \widehat{T}W/[\widehat{T}W,\widehat{T}W]=\drzf[Ass]{\widehat{T}W} \]
of bidegree $(0,0)$ and order $0$, which it maps to zero. By the  Cartan homotopy formula of Lemma \ref{lem_schf} part (i) we know that $L_\xi$ is nullhomotopic, whence the result.
\end{proof}

\section{Relations between Commutative, Associative and Lie Geometries}\label{relations}

In this section we will discuss some aspects of noncommutative geometry which describe the relationship between the commutative, associative and Lie versions of the theory described in section \ref{sec_ncgeom}. We will only concern ourselves with describing this relationship for the pro-free algebras $\widehat{S}W$, $\widehat{T}W$ and $\widehat{L}W$. One of the main purposes of this section will be to lay the groundwork for sections \ref{sec_iachom}, \ref{sec_hdgdec}, and \ref{cyclichodge} and also to expand upon a theorem of Kontsevich \cite{kontsg}. This theorem will allow us to consider the closed commutative, associative or Lie 2-forms of the pro-free algebras $\widehat{S}U^*$, $\widehat{T}U^*$ or $\widehat{L}U^*$ respectively, as linear maps $\alpha:TU \to \gf$.

We first of all establish some fundamental identities for the operators $N$ and $z$ defined in section \ref{sec_notcon}. We then introduce a series of maps which we use to establish some relationships between the complexes \dr[Com]{\widehat{S}W}, \dr[Ass]{\widehat{T}W} and \dr[Lie]{\widehat{L}W}. We use these relationships to prove Theorem \ref{thm_tfmaps} and we use this theorem to interpret the closed 2-forms as linear maps. In particular we identify the 2-forms giving rise to symmetric bilinear forms.

\begin{lemma} \label{lem_nrmdif}
Let $W$ be a graded profinite module:
\begin{enumerate}
\item[(i)]
The following identity holds:
\begin{equation} \label{eqn_nrmdifdummy}
\widehat{T}W/[\widehat{T}W,\widehat{T}W] = \prod_{i=0}^\infty (W^{\cotimes i})_{Z_i}.
\end{equation}
\item[(ii)]
The following diagram commutes:
\begin{displaymath}
\xymatrix{ W \cotimes \widehat{T}W \ar^{\thiso[Ass]}_{\cong}[r] \ar@{=}[dd] & \drof[Ass]{\widehat{T}W} & \drzf[Ass]{\widehat{T}W} \ar_{d}[l] \ar@{=}[d] \\ & \widehat{T}W \ar@{->>}[r] & \widehat{T}W/[\widehat{T}W,\widehat{T}W] \\ \prod_{i=1}^\infty W^{\cotimes i} & \prod_{i=1}^\infty W^{\cotimes i} \ar_{N}[l] \ar@{^{(}->}[u] \ar@{->>}[r] & \prod_{i=1}^\infty (W^{\cotimes i})_{Z_i} \ar@{^{(}->}[u] \\ }
\end{displaymath}
\end{enumerate}
\end{lemma}

\begin{proof}
Let $x:=x_1\cotimes\ldots\cotimes x_n \in W^{\cotimes n} \subset \widehat{T}W$:
\begin{enumerate}
\item[(i)]
For $1\leq i \leq n-1$ we calculate;
\begin{equation} \label{eqn_cyccom}
\begin{split}
(1-z_n^i)\cdot x & = x_1\cotimes\ldots\cotimes x_n - (-1)^{(|x_1|+\ldots+|x_i|)(|x_{i+1}|+\ldots+|x_n|)}x_{i+1}\cotimes\ldots\cotimes x_n\cotimes x_1\cotimes \ldots x_i, \\
& = [x_1\cotimes\ldots\cotimes x_i,x_{i+1}\cotimes\ldots\cotimes x_n]. \\
\end{split}
\end{equation}
This calculation suffices to establish equation \eqref{eqn_nrmdifdummy}.

\item[(ii)]
\begin{equation} \label{eqn_nrmdif}
\begin{split}
d(x) &= \sum_{i=1}^n (-1)^{|x_1|+\ldots+|x_{i-1}|} x_1\cotimes\ldots\cotimes x_{i-1}\cdot dx_i \cdot x_{i+1}\cotimes\ldots\cotimes x_n, \\
& = \sum_{i=1}^n (-1)^{(|x_1|+\ldots+|x_{i-1}|)(|x_i|+\ldots+|x_n|)} dx_i \cdot x_{i+1}\cotimes\ldots\cotimes x_n \cotimes x_1 \cotimes \ldots \cotimes x_{i-1} \mod [-,-]. \\
N_n\cdot x & = \sum_{i=0}^{n-1} (-1)^{(|x_1|+\ldots+|x_i|)(|x_{i+1}|+\ldots+|x_n|)}x_{i+1}\cotimes\ldots\cotimes x_n\cotimes x_1 \cotimes \ldots \cotimes x_i. \\
\end{split}
\end{equation}
From this calculation we see that $\thiso[Ass](N_n \cdot x) = d(x)$ and therefore the diagram is commutative.
\end{enumerate}
\end{proof}

Our next task will be to identify the closed 2-forms of \dr[Com]{\widehat{S}U^*}, \dr[Ass]{\widehat{T}U^*} and \dr[Lie]{\widehat{L}U^*} with  submodules of the module consisting of all linear maps $\alpha:TU \to \gf$. In order to do this we will need to introduce maps between the de Rham complexes and establish some basic properties of them.

Let $W$ be a graded profinite module. The completed tensor algebra $\widehat{T}W$ is a Lie algebra under the commutator and the differential envelope \de[Ass]{\widehat{T}W} is a differential graded Lie algebra under the commutator. There is a series of canonical Lie algebra inclusions;
\[ \widehat{L}W \hookrightarrow \widehat{T}W \hookrightarrow \de[Ass]{\widehat{T}W}. \]
By Proposition \ref{prop_difenv} the composition of these inclusions extends to a map of differential graded Lie algebras
\[ l':\de[Lie]{\widehat{L}W} \to \de[Ass]{\widehat{T}W}. \]
One can easily check that the following formula holds for all $a,b,c \in \de[Ass]{\widehat{T}W}$:
\begin{equation} \label{eqn_modcom}
a[b,c]=[a,b]c \mod [-,-].
\end{equation}
This allows us to define a map $l: \dr[Lie]{\widehat{L}W} \to \dr[Ass]{\widehat{T}W}$ which is induced by the map $l'$ and given by
\begin{equation} \label{eqn_liasmp}
l(x\cotimes y):=l'(x)l'(y); \quad x,y \in \de[Lie]{\widehat{L}W}.
\end{equation}
It follows from equation \eqref{eqn_modcom} that this map is well defined modulo the relations of Definition \ref{def_derham} part (c).

We can define another map $p:\dr[Ass]{\widehat{T}W} \to \dr[Com]{\widehat{S}W}$ as follows: the commutative algebra $\widehat{S}W$ is obviously an associative algebra and the canonical projection $\widehat{T}W \twoheadrightarrow \widehat{S}W$ is a map of associative algebras. By Proposition \ref{prop_difenv} this map extends to a map of DGAs
\[ p':\de[Ass]{\widehat{T}W} \to \de[Com]{\widehat{S}W}. \]
Clearly this map is zero on commutators and hence induces a map
\begin{equation} \label{eqn_ascomp}
p: \dr[Ass]{\widehat{T}W} \to \dr[Com]{\widehat{S}W}.
\end{equation}

Next we need to introduce a map $i:\widehat{S}W \to \prod_{i=0}^\infty (W^{\cotimes i})^{S_i}$. This is defined as the canonical map identifying coinvariants with invariants and is given by
\begin{equation} \label{eqn_symiso}
i(x_1\cotimes\ldots\cotimes x_n):= \sum_{\sigma \in S_n} \sigma\cdot x_1\cotimes\ldots\cotimes x_n.
\end{equation}

Finally, we can describe a map $j:\drof[Com]{\widehat{S}W} \to \drof[Ass]{\widehat{T}W}$. This is defined by the following commutative diagram which uses the isomorphisms defined by Lemma \ref{lem_ofisom} and the map $i$ defined above by equation \eqref{eqn_symiso}:
\begin{equation} \label{fig_caofmp}
\xymatrix{ W \cotimes \widehat{S}W \ar^{\thiso[Com]}_{\cong}[r] \ar^{1 \cotimes i}[d] & \drof[Com]{\widehat{S}W} \ar^{j}[d] \\ W \cotimes \widehat{T}W \ar^{\thiso[Ass]}_{\cong}[r] & \drof[Ass]{\widehat{T}W} \\ }
\end{equation}

The map $l:\dr[Lie]{\widehat{L}W} \to \dr[Ass]{\widehat{T}W}$ defined by equation \eqref{eqn_liasmp} can be lifted to a unique map
\[ \frac{\drof[Lie]{\widehat{L}W}}{d(\drzf[Lie]{\widehat{L}W})} \to \frac{\drof[Ass]{\widehat{T}W}}{d(\drzf[Ass]{\widehat{T}W})}. \]
Similarly, the map $p:\drof[Ass]{\widehat{T}W}\to\drof[Com]{\widehat{S}W}$  of \eqref{eqn_ascomp} can be lifted to a unique map
\[ \frac{\drof[Ass]{\widehat{T}W}}{d(\drzf[Ass]{\widehat{T}W})} \to \frac{\drof[Com]{\widehat{S}W}}{d(\drzf[Com]{\widehat{S}W})}. \]
Abusing notation, we will denote the induced maps by the same symbols $l$ and $p$.

\begin{lemma} \label{lem_clsdtf}
Let $W$ be a graded profinite module:
\begin{enumerate}
\item[(i)]
The map $l:\drof[Lie]{\widehat{L}W} \to \drof[Ass]{\widehat{T}W}$ is injective.
\item[(ii)]
The map
\[ l:\frac{\drof[Lie]{\widehat{L}W}}{d(\drzf[Lie]{\widehat{L}W})} \to \frac{\drof[Ass]{\widehat{T}W}}{d(\drzf[Ass]{\widehat{T}W})} \]
is also injective.
\item[(iii)]
The map $j:\drof[Com]{\widehat{S}W}\to\drof[Ass]{\widehat{T}W}$ defined by diagram \eqref{fig_caofmp} is injective and the map $p:\drof[Ass]{\widehat{T}W}\to\drof[Com]{\widehat{S}W}$  is surjective.
\item[(iv)]
The map $j:\drof[Com]{\widehat{S}W}\to\drof[Ass]{\widehat{T}W}$ can be lifted to a unique map
\[ \frac{\drof[Com]{\widehat{S}W}}{d(\drzf[Com]{\widehat{S}W})} \to \frac{\drof[Ass]{\widehat{T}W}}{d(\drzf[Ass]{\widehat{T}W})} \]
which will be denoted by the same letter $j$ by the customary abuse of notation.
\item[(v)]
The map
\[ j:\frac{\drof[Com]{\widehat{S}W}}{d(\drzf[Com]{\widehat{S}W})} \to \frac{\drof[Ass]{\widehat{T}W}}{d(\drzf[Ass]{\widehat{T}W})} \]
is injective and the map
\[ p:\frac{\drof[Ass]{\widehat{T}W}}{d(\drzf[Ass]{\widehat{T}W})} \to \frac{\drof[Com]{\widehat{S}W}}{d(\drzf[Com]{\widehat{S}W})} \]
is surjective.
\end{enumerate}
All told, we have the following diagram of injections and surjections:
\begin{displaymath}
\xymatrix{ \drof[Com]{\widehat{S}W} \ar@{->>}[r] \ar@<-1ex>@{^{(}->}_{j}[d] & \frac{\drof[Com]{\widehat{S}W}}{d(\drzf[Com]{\widehat{S}W})} \ar@<-1ex>@{^{(}->}_{j}[d] \\  \drof[Ass]{\widehat{T}W} \ar@{->>}[r] \ar@<-1ex>@{->>}_{p}[u] & \frac{\drof[Ass]{\widehat{T}W}}{d(\drzf[Ass]{\widehat{T}W})} \ar@<-1ex>@{->>}_{p}[u] \\ \drof[Lie]{\widehat{L}W} \ar@{->>}[r] \ar@{^{(}->}^{l}[u] & \frac{\drof[Lie]{\widehat{L}W}}{d(\drzf[Lie]{\widehat{L}W})} \ar@{^{(}->}^{l}[u] \\ }
\end{displaymath}
\end{lemma}

\begin{proof}
From the identities;
\begin{enumerate}
\item[(a)] $\prod_{i=1}^\infty W^{\cotimes i} = W \cotimes \widehat{T}W$,
\item[(b)] $\drzf[Ass]{\widehat{T}W} = \widehat{T}W/[\widehat{T}W,\widehat{T}W]$,
\item[(c)] $\drzf[Lie]{\widehat{L}W} = \widehat{L}W \cotimes \widehat{L}W / ([x,y]\cotimes z = x \cotimes [y,z])$;
\end{enumerate}
it follows that there are maps;
\[ \pi_\mathrm{Ass} : W \cotimes \widehat{T}W \to \drzf[Ass]{\widehat{T}W}, \]
\[ \pi_\mathrm{Lie} : W \cotimes \widehat{L}W \to \drzf[Lie]{\widehat{L}W}; \]
which are just the canonically defined projection maps. Let $u:\widehat{L}W \to \widehat{T}W$ be the canonical inclusion of the Lie subalgebra $\widehat{L}W$ into $\widehat{T}W$. These maps fit together with the maps \thiso[Ass] and \thiso[Lie] defined by Lemma \ref{lem_ofisom} and the map $l$ defined by equation \eqref{eqn_liasmp} to form the following commutative diagram:
\begin{equation} \label{fig_clsdtf}
\xymatrix{ \drof[Ass]{\widehat{T}W} & W \cotimes \widehat{T}W \ar^{\pi_\mathrm{Ass}}[r] \ar_{\thiso[Ass]}^{\cong}[l] & \drzf[Ass]{\widehat{T}W} \\ \drof[Lie]{\widehat{L}W} \ar^{l}[u] & W \cotimes \widehat{L}W \ar^{1 \cotimes u}[u] \ar^{\pi_\mathrm{Lie}}[r] \ar_{\thiso[Lie]}^{\cong}[l] & \drzf[Lie]{\widehat{L}W} \ar^{l}[u] \\ }
\end{equation}
\begin{enumerate}
\item[(i)]
Since $u$ is injective we conclude from diagram \eqref{fig_clsdtf} that the map
\[ l:\drof[Lie]{\widehat{L}W} \to \drof[Ass]{\widehat{T}W} \]
is injective.
\item[(ii)]
Let $x \in \drof[Lie]{\widehat{L}W}$ be a 1-form of order $n$ and suppose that $l(x) \in d(\drzf[Ass]{\widehat{T}W})$. By Lemma \ref{lem_nrmdif} part (ii), \thiso[Ass] satisfies the identity,
\begin{equation} \label{eqn_clsdtfdummyb}
\thiso[Ass] \circ N = d \circ \pi_\mathrm{Ass}
\end{equation}
and therefore identifies $d(\drzf[Ass]{\widehat{T}W})$ with the module of invariants $\prod_{i=1}^\infty (W^{\cotimes i})^{Z_i}$. Using the commutativity of diagram \eqref{fig_clsdtf} we perform the following calculation:
\begin{displaymath}
\begin{split}
\thiso[Ass]^{-1}l(x) & = \frac{1}{n+1} N_{n+1} \cdot \thiso[Ass]^{-1}l(x), \\
l(x) & = \frac{1}{n+1} \thiso[Ass] N_{n+1} \thiso[Ass]^{-1}l(x),  \\
& = \frac{1}{n+1} ld\pi_\mathrm{Lie}\thiso[Lie]^{-1}(x). \\
\end{split}
\end{displaymath}
From part (i) of this lemma we conclude that $x=\frac{1}{n+1}d\pi_\mathrm{Lie}\thiso[Lie]^{-1}(x)$ and hence the map
\[ l: \frac{\drof[Lie]{\widehat{L}W}}{d(\drzf[Lie]{\widehat{L}W})} \to \frac{\drof[Ass]{\widehat{T}W}}{d(\drzf[Ass]{\widehat{T}W})} \]
is injective.
\item[(iii)]
One can easily verify the following equation from which the assertion follows immediately:
\begin{equation} \label{eqn_clsdtfdummya}
pj(dx_0 \cdot x_1 \cotimes \ldots \cotimes x_n) = n!dx_0 \cdot x_1 \cotimes \ldots \cotimes x_n.
\end{equation}
\item[(iv)]
Let $x:=x_1 \cotimes \ldots \cotimes x_{n+1} \in \widehat{T}W$ and let $\bar{x}$ denote its image in the quotient $\widehat{S}W =\drzf[Com]{\widehat{S}W}$, then
\begin{displaymath}
\begin{split}
\thiso[Ass]^{-1}jd(\bar{x}) & = \sum_{\sigma \in S_n}\sum_{\tau \in Z_{n+1}} \sigma\tau \cdot x, \\
& = \sum_{\sigma \in S_{n+1}} \sigma \cdot x, \\
& = \sum_{\tau \in Z_{n+1}}\sum_{\sigma \in S_n} \tau\sigma \cdot x. \\
\end{split}
\end{displaymath}
Using equation \eqref{eqn_clsdtfdummyb} we obtain;
\[ jd(\bar{x}) = d\pi_\mathrm{Ass}\left(\sum_{\sigma \in S_n} \sigma \cdot x\right), \]
hence the map defined by diagram \eqref{fig_caofmp} can be lifted to a map
\[ j:\frac{\drof[Com]{\widehat{S}W}}{d(\drzf[Com]{\widehat{S}W})} \to \frac{\drof[Ass]{\widehat{T}W}}{d(\drzf[Ass]{\widehat{T}W})}. \]
\item[(v)]
This is an immediate consequence of equation \eqref{eqn_clsdtfdummya}.
\end{enumerate}
\end{proof}

\begin{rem} \label{rem_vecseq}
Let $W$ be a graded profinite module. We will now describe a sequence of Lie algebra homomorphisms;
\[ \Der(\widehat{L}W) \hookrightarrow \Der(\widehat{T}W) \twoheadrightarrow \Der(\widehat{S}W). \]

Any vector field $\xi \in \Der(\widehat{T}W)$ is specified uniquely by its restriction $\xi: W \to \widehat{T}W$ (cf. Proposition \ref{prop_dergen}) and likewise for $\widehat{L}W$ and $\widehat{S}W$. It follows that any vector field $\xi \in \Der(\widehat{L}W)$ can be extended uniquely to a vector field
\[ \xi: \widehat{T}W \to \widehat{T}W. \]
In fact, this vector field is a Hopf algebra derivation (i.e. a derivation and a coderivation) when $\widehat{T}W$ is equipped with the cocommutative comultiplication (\emph{shuffle coproduct}) defined in Remark \ref{rem_ccmult}. Furthermore all continuous Hopf algebra derivations are obtained from $\Der(\widehat{L}W)$ in this way. This is because the Lie subalgebra of primitive elements of this Hopf algebra coincides with the Lie subalgebra $\widehat{L}W$. This establishes a one-to-one correspondence between $\Der(\widehat{L}W)$ and continuous Hopf algebra derivations on $\widehat{T}W$.

Any map $\xi:W \to \widehat{T}W$ gives rise to a map $\tilde{\xi}:W \to \widehat{S}W$ simply by composing it with the canonical projection $\widehat{T}W \to \widehat{S}W$. It follows that any vector field $\xi \in \Der(\widehat{T}W)$ can be lifted to a unique vector field
\[ \tilde{\xi}: \widehat{S}W \to \widehat{S}W \]
and that all the vector fields in $\Der(\widehat{S}W)$ are obtained from $\Der(\widehat{T}W)$ in this way.
\end{rem}

We will now formulate a lemma, the proof of which is a simple check, which describes how this sequence of Lie algebra homomorphisms interact with our maps $l$ and $p$:

\begin{lemma} \label{lem_plcomm}
Let $W$ be a graded profinite module:
\begin{enumerate}
\item[(i)]
Given any vector field $\xi:\widehat{L}W \to \widehat{L}W$, the following diagrams commute:
\begin{displaymath}
\xymatrix{ \dr[Ass]{\widehat{T}W} \ar^{i_\xi}[r] & \dr[Ass]{\widehat{T}W} \\ \dr[Lie]{\widehat{L}W} \ar^{i_\xi}[r] \ar@{^{(}->}^{l}[u] & \dr[Lie]{\widehat{L}W} \ar@{^{(}->}^{l}[u] \\ }
\qquad
\xymatrix{ \dr[Ass]{\widehat{T}W} \ar^{L_\xi}[r] & \dr[Ass]{\widehat{T}W} \\ \dr[Lie]{\widehat{L}W} \ar^{L_\xi}[r] \ar@{^{(}->}^{l}[u] & \dr[Lie]{\widehat{L}W} \ar@{^{(}->}^{l}[u] \\ }
\end{displaymath}
\item[(ii)]
Given any vector field $\xi:\widehat{T}W \to \widehat{T}W$, the following diagrams commute:
\begin{displaymath}
\xymatrix{ \dr[Com]{\widehat{S}W} \ar^{i_{\tilde{\xi}}}[r] & \dr[Com]{\widehat{S}W} \\ \dr[Ass]{\widehat{T}W} \ar^{i_\xi}[r] \ar@{->>}^{p}[u] & \dr[Ass]{\widehat{T}W} \ar@{->>}^{p}[u] \\ }
\qquad
\xymatrix{ \dr[Com]{\widehat{S}W} \ar^{L_{\tilde{\xi}}}[r] & \dr[Com]{\widehat{S}W} \\ \dr[Ass]{\widehat{T}W} \ar^{L_\xi}[r] \ar@{->>}^{p}[u] & \dr[Ass]{\widehat{T}W} \ar@{->>}^{p}[u] \\ }
\end{displaymath}
\end{enumerate}
\end{lemma}
\noproof

Using Lemma \ref{lem_clsdtf} we can prove the following interesting theorem. This result was proved in \cite{kontsg} for closed associative 2-forms and stated without a proof for closed Lie 2-forms. We will describe the analogous result in the commutative case as well. Given a graded profinite module $W$ we can consider the module of commutators
\[ [W,\textstyle{\prod_{i=0}^\infty (W^{\cotimes i})^{S_i}}] \subset [\widehat{T}W,\widehat{T}W] \subset \widehat{T}W. \]

\begin{theorem} \label{thm_tfmaps}
Let $W$ be a graded profinite module. There are isomorphisms:
\begin{enumerate}
\item[(i)]
\[ \zetiso[Ass]:\{ \omega \in \drtf[Ass]{\widehat{T}W} : d\omega = 0 \} \to [\widehat{T}W,\widehat{T}W]. \]
\item[(ii)]
\[ \zetiso[Com]:\{ \omega \in \drtf[Com]{\widehat{S}W} : d\omega = 0 \} \to [W,\textstyle{\prod_{i=0}^\infty (W^{\cotimes i})^{S_i}}]. \]
\item[(iii)]
\[ \zetiso[Lie]:\{ \omega \in \drtf[Lie]{\widehat{L}W} : d\omega = 0 \} \to [\widehat{L}W,\widehat{L}W]. \]
\end{enumerate}
\end{theorem}

\begin{proof}
Firstly if $X$ is one of the three algebras $\widehat{S}W$, $\widehat{T}W$ or $\widehat{L}W$ then by Lemma \ref{lem_poinca} there is an isomorphism of modules,
\[ f: \drof{X}/d(\drzf{X}) \to \{ \omega \in \drtf{X} : d\omega = 0 \}; \]
where $f$ is defined by the formula $f(x):=dx$. We will denote this isomorphism by $f_\mathrm{Com}$, $f_\mathrm{Ass}$ or $f_\mathrm{Lie}$ respectively.
\begin{enumerate}
\item[(i)]
First of all note that by Lemma \ref{lem_nrmdif}, the isomorphism \thiso[Ass] defined in Lemma \ref{lem_ofisom} induces a map,
\[ \bar{\Theta}_\mathrm{Ass}: \frac{\prod_{i=1}^\infty W^{\cotimes i}}{N\cdot\prod_{i=1}^\infty W^{\cotimes i}} \to \frac{\drof[Ass]{\widehat{T}W}}{d(\drzf[Ass]{\widehat{T}W})} \]
which is also an isomorphism.

Equation \eqref{eqn_cyccom} gives us the identity
\[ (1-z)\cdot\prod_{i=1}^\infty W^{\cotimes i} = [\widehat{T}W,\widehat{T}W]. \]
This allows us to define a module isomorphism
\[ g: \frac{\drof[Ass]{\widehat{T}W}}{d(\drzf[Ass]{\widehat{T}W})} \to [\widehat{T}W,\widehat{T}W] \]
by the formula $g(x):=(1-z)\cdot\bar{\Theta}_\mathrm{Ass}^{-1}(x)$. We then define the isomorphism \zetiso[Ass] by the formula
\[ \zetiso[Ass]:= g \circ f_\mathrm{Ass}^{-1}. \]
\item[(ii)]
By Lemma \ref{lem_clsdtf} there is an injection
\[ j: \frac{\drof[Com]{\widehat{S}W}}{d(\drzf[Com]{\widehat{S}W})} \hookrightarrow \frac{\drof[Ass]{\widehat{T}W}}{d(\drzf[Ass]{\widehat{T}W})}. \]
Composing $j$ with $g$ yields an injection $g\circ j$. It requires a simple check using the definition of $j$ and equation \eqref{eqn_cyccom} to see that the image of $g\circ j$ is $[W,\prod_{i=1}^\infty (W^{\cotimes i})^{S_i}]$. It follows that the map
\[ \zetiso[Com]:\{ \omega \in \drtf[Com]{\widehat{S}W} : d\omega = 0 \} \to [W,\textstyle{\prod_{i=0}^\infty (W^{\cotimes i})^{S_i}}] \]
defined by the formula $\zetiso[Com](x):=gjf_\mathrm{Com}^{-1}(x)$ is an isomorphism.
\item[(iii)]
By Lemma \ref{lem_clsdtf} there is an injection
\[l: \frac{\drof[Lie]{\widehat{L}W}}{d(\drzf[Lie]{\widehat{L}W})} \hookrightarrow \frac{\drof[Ass]{\widehat{T}W}}{d(\drzf[Ass]{\widehat{T}W})}. \]
Composing $l$ with $g$ yields an injection $g\circ l$. Again it requires a simple check using the definitions and equation \eqref{eqn_cyccom} to see that the image of $g \circ l$ is $[W,\widehat{L}W]$, however, using the Jacobi identity we obtain the equality
\[ [W,\widehat{L}W] = [\widehat{L}W,\widehat{L}W]. \]
It follows that the map
\[ \zetiso[Lie]:\{ \omega \in \drtf[Lie]{\widehat{L}W} : d\omega = 0 \} \to [\widehat{L}W,\widehat{L}W] \]
defined by the formula $\zetiso[Lie](x):=glf_\mathrm{Lie}^{-1}(x)$ is an isomorphism.
\end{enumerate}
\end{proof}

\begin{rem}
It follows from the definitions that the maps \zetiso[Com], \zetiso[Ass] and \zetiso[Lie] satisfy the following relations:
\begin{equation} \label{eqn_zetrel}
\begin{split}
\zetiso[Ass]\circ d & = (1-z) \circ \thiso[Ass]^{-1}, \\
\zetiso[Com] \circ d & = \zetiso[Ass] \circ d \circ j, \\
\zetiso[Lie] & = \zetiso[Ass] \circ l. \\
\end{split}
\end{equation}
\end{rem}

If we set $W:=U^*$ in Theorem \ref{thm_tfmaps}, where $U$ is a free graded module, then we can interpret the closed 2-forms as linear maps $(TU)^* \to \gf$. We now want to concentrate on 2-forms of order zero and determine what linear maps in $(TU)^*$ they give rise to under the maps \zetiso[Com], \zetiso[Ass] and \zetiso[Lie]. 2-forms of order zero are naturally closed since it is easy to see that any 2-form $\omega$ \emph{of order zero} has the following form:
\begin{displaymath}
\begin{array}{lcll}
\omega \in \drtf[Com]{\widehat{S}U^*} & \Rightarrow & \omega = \sum_i dx_i\cdot dy_i; & x_i,y_i \in U^*; \\
\omega \in \drtf[Ass]{\widehat{T}U^*} & \Rightarrow & \omega = \sum_i dx_i\cdot dy_i; & x_i,y_i \in U^*; \\
\omega \in \drtf[Lie]{\widehat{L}U^*} & \Rightarrow & \omega = \sum_i dx_i\cotimes dy_i; & x_i,y_i \in U^*. \\
\end{array}
\end{displaymath}

We have the following lemma which describes the maps that two forms of order zero give rise to:

\begin{lemma} \label{lem_tfmbmp}
Let $U$ be a free graded module. Every map in the following commutative diagram is an isomorphism:
\begin{displaymath}
\xymatrix{ \{ \omega \in \drtf[Com]{\widehat{S}U^*} : \ord(\omega) = 0 \} \ar_{j}[d] \ar^{\zetiso[Com]}[rd] \\ \{ \omega \in \drtf[Ass]{\widehat{T}U^*} : \ord(\omega) = 0 \} \ar^-{\zetiso[Ass]}[r] & (\Lambda^2 U)^* \\ \{ \omega \in \drtf[Lie]{\widehat{L}U^*} : \ord(\omega) = 0 \} \ar^{l}[u] \ar^{\zetiso[Lie]}[ru] \\ }
\end{displaymath}
\end{lemma}

\begin{proof}
This follows from Theorem \ref{thm_tfmaps} and the identity
\[ [U^*,U^*] = (1-z_2)\cdot \left(U^*\cotimes U^*\right) = (\Lambda^2 U)^*. \]
\end{proof}

This means there is a one-to-one correspondence between the module of 2-forms (commutative, associative or Lie) of order zero and the module consisting of all skew-symmetric bilinear forms
\[ \innprod{U}. \]
We will now recall the definition of a nondegenerate 2-form (cf. \cite{ginzsg}) and a nondegenerate bilinear form:

\begin{defi} \label{def_nondeg}
\
\begin{enumerate}
\item[(i)]
Let $X$ be either a formal graded commutative, associative or Lie algebra and let $\omega \in \drtf{X}$ be a 2-form. We say $\omega$ is nondegenerate if the following map is a bijection;
\begin{displaymath}
\begin{array}{ccc}
\Der(X) & \to & \drof{X}, \\
\xi & \mapsto & i_\xi(\omega). \\
\end{array}
\end{displaymath}
\item[(ii)]
Let $U$ be a free graded module of finite rank and let \innprod{U} be a bilinear form. We say that $\langle -,- \rangle$ is nondegenerate if the following map is a bijection;
\begin{displaymath}
\begin{array}{ccc}
U & \to & U^*, \\
u & \mapsto & [ x \mapsto \langle u,x \rangle ]. \\
\end{array}
\end{displaymath}
If in addition $\langle -,- \rangle$ is symmetric then we will call it an \emph{inner product} on $U$.
\end{enumerate}
\end{defi}

We have the following proposition relating the two notions:

\begin{prop} \label{prop_nondeg}
Let $U$ be a free graded module of finite rank and let $X$ be one of the three algebras $\widehat{S}U^*$, $\widehat{T}U^*$ or $\widehat{L}U^*$. Let $\omega \in \drtf{X}$ be a 2-form of order zero. Let $\langle -,- \rangle:=\zeta(\omega)$ be the skew-symmetric bilinear form corresponding to the 2-form $\omega$, then $\omega$ is nondegenerate if and only if $\langle -,- \rangle$ is nondegenerate.
\end{prop}

\begin{proof}
We will treat the three cases $X=\widehat{S}U^*$, $X=\widehat{T}U^*$ and $X=\widehat{L}U^*$ simultaneously. Let $x_1,\ldots,x_n$ be a basis of the free module $U$. There are coefficients $a_{ij} \in \gf$ such that
\[ \omega = \sum_{1\leq i,j \leq n} a_{ij}dx_i^* dx_j^*. \]
It follows from the definition of $\zeta$ that
\[ \langle a,b \rangle = \sum_{1 \leq i,j \leq n} a_{ij}\left[(-1)^{|x_i|(|x_j|+1)}x_j^*(a)x_i^*(b) - (-1)^{|x_i|}x_i^*(a)x_j^*(b)\right]. \]

Let us define the map $\Phi: \Der(X) \to \drof{X}$ by the formula
\[ \Phi(\xi):= i_\xi(\omega). \]
Let us also define a map $D:U \to U^*$ by the formula
\[ D(u):= [x \mapsto \langle u,x \rangle]. \]

We calculate $\Phi$ as follows: Let $\xi \in \Der(X)$, then
\begin{displaymath}
\begin{split}
\Phi(\xi) = i_\xi(\omega) & = \sum_{1 \leq i,j \leq n} a_{ij}[\xi(x_i^*)\cdot dx_j^* + (-1)^{(|x_i|+1)(|\xi|+1)}dx_i^*\cdot \xi(x_j^*)], \\
& =  \sum_{1 \leq i,j \leq n} a_{ij}[(-1)^{(|x_j|+1)(|x_i|+|\xi|)}dx_j^*\cdot \xi(x_i^*) + (-1)^{(|x_i|+1)(|\xi|+1)}dx_i^*\cdot \xi(x_j^*)]. \\
\end{split}
\end{displaymath}
We calculate $D$ as follows: For all $a \in U$, and $1\leq k \leq n$;
\begin{displaymath}
\begin{split}
D(x_k)[a] = \langle x_k,a \rangle & = \sum_{1 \leq i,j \leq n} a_{ij}\left[(-1)^{|x_i|(|x_j|+1)}x_j^*(x_k)x_i^*(a) - (-1)^{|x_i|}x_i^*(x_k)x_j^*(a)\right], \\
& = \sum_{1\leq i \leq n}\left[(-1)^{|x_i|(|x_k|+1)}a_{ik} - (-1)^{|x_k|}a_{ki}\right]x_i^*(a). \\
\end{split}
\end{displaymath}

Let $h:U^* \to X$ be the canonical inclusion of the submodule $U^*$ into $X$. Since any vector field $\xi \in \Der(X)$ is completely determined by its restriction $\xi\circ h$ (cf. Proposition \ref{prop_dergen}) we have the following isomorphism of graded modules:
\begin{displaymath}
\begin{array}{ccc}
\Der(X) & \cong & \Hom_\gf(U^*,X), \\
\xi & \mapsto & \xi\circ h. \\
\end{array}
\end{displaymath}
It follows from the preceding calculations that the following diagram is commutative:
\[ \xymatrix{ \drof{X} \ar^{\Theta^{-1}}_{\cong}[r] & U^* \otimes X \ar^{u \otimes x \mapsto (-1)^{|\omega|+(|u|+1)(|x|+1)}u \otimes x}[rr] && U^* \otimes X \\ \Der(X) \ar^{\Phi}[u] \ar_-{\cong}^-{\xi \mapsto \xi \circ h}[r] & \Hom_\gf(U^*,X) \ar@{=}[r] & U^{**} \otimes X \ar@{=}[r] & U \otimes X \ar_{D \otimes 1}[u] \\ } \]
We conclude that $\Phi$ is a bijection if and only if $D$ is a bijection.
\end{proof}

\begin{rem} \label{rem_canfrm}
For simplicity's sake, let us assume we are working over a field of characteristic zero which is closed under taking square roots and that $U$ is finite dimensional; then every homogeneous 2-form $\omega \in \drtf{X}$ (where $X=\widehat{S}U^*$, $\widehat{T}U^*$ or $\widehat{L}U^*$) of order zero has a canonical form.

Suppose that $\omega$ has \emph{even} degree, then there exist linearly independent vectors;
\begin{equation} \label{eqn_canfrmdummy}
p_1,\ldots,p_n;q_1,\ldots,q_n;x_1,\ldots,x_m \in U
\end{equation}
(where the $p_i$'s and the $q_i$'s are \emph{even} and the $x_i$'s are \emph{odd}) such that $\omega$ has the following form:
\[ \omega = \sum_{i=1}^n d p_i^* d q_i^* + \sum_{i=1}^m d x_i^* d x_i^*. \]
If $\omega$ is nondegenerate then \eqref{eqn_canfrmdummy} is a basis for $U$. The bilinear form $\langle -,- \rangle:=\zeta(\omega)$ is given by the formula:
\begin{displaymath}
\begin{array}{rcr}
\langle q_i,p_j \rangle = \frac{1}{2}\langle x_i,x_j \rangle & = & \delta_{ij}, \\
\langle x_i,p_j \rangle = \langle x_i,q_j \rangle = \langle p_i,p_j \rangle = \langle q_i,q_j \rangle & = & 0. \\
\end{array}
\end{displaymath}

Now suppose that $\omega$ has \emph{odd} degree, then there exist linearly independent vectors;
\begin{equation} \label{eqn_canfrmdummya}
x_1,\ldots,x_n;y_1,\ldots,y_n \in U
\end{equation}
(where the $x_i$'s have \emph{odd} degree and the $y_i$'s have \emph{even} degree) such that $\omega$ has the following form:
\[ \omega = \sum_{i=1}^n d x_i^* d y_i^*. \]
Again, if $\omega$ is nondegenerate then \eqref{eqn_canfrmdummya} is a basis for $U$. The bilinear form $\langle -,- \rangle:=\zeta(\omega)$ is given by the formula:
\begin{displaymath}
\begin{array}{rcr}
\langle x_i,y_j \rangle & = & \delta_{ij}, \\
\langle x_i,x_j \rangle = \langle y_i,y_j \rangle & = & 0. \\
\end{array}
\end{displaymath}
\end{rem}

\section{Infinity-algebra Prerequisites} \label{prerequisites}

In this section we will review the definitions of three types of infinity-algebra, namely {\ci}, {\ai} and {\li}. These are the (strong) homotopy generalisations of commutative, associative and Lie algebras respectively. We will also define the appropriate notion of a unital \ai and \ci-algebra. We shall utilise the duality between derivations and coderivations to define an $\infty$-structure as a homological vector field on a certain formal supermanifold. This approach will provide us with the natural framework in which to apply constructions from noncommutative geometry.

Let $V$ be a free graded module and choose a topological basis $\boldsymbol{t}:=\{t_i\}_{i \in I}$ of $\Sigma V^*$, then $\ctalg{V}=\gf\langle\langle\boldsymbol{t}\rangle\rangle$. Recall from section \ref{sec_notcon} that we say a vector field (continuous derivation) $\xi:\ctalg{V}\to\ctalg{V}$ vanishes at zero if it has the form
\[ \xi=\sum_{i \in I}A_i(\boldsymbol{t})\partial_{t_i}, \]
where the power series $A_i(\boldsymbol{t}), i \in I$ have vanishing constant terms.

We will now recall from \cite{azksch}, \cite{getjon} and \cite{lazmod} the definition of an \li, {\ai} and \ci-structure on a free graded module $V$. An \li, {\ai} or \ci-algebra is a free graded module together with an \li, {\ai} or \ci-structure. Given an element $x$ in a graded algebra $A$ we define the derivation $\ad x:A \to A$ by the formula $\ad x(y):=[x,y]$.

\begin{defi} \label{def_infstr}
Let $V$ be a free graded module:
\begin{enumerate}
\item[(a)]
An \li-structure on $V$ is a vector field
\[ m:\csalg{V} \to \csalg{V} \]
of degree one and vanishing at zero, such that $m^2=0$.
\item[(b)]
An \ai-structure on $V$ is a vector field
\[ m:\ctalg{V} \to \ctalg{V} \]
of degree one and vanishing at zero, such that $m^2=0$.
\item[(c)]
A \ci-structure on $V$ is a vector field
\[ m:\clalg{V} \to \clalg{V} \]
of degree one, such that $m^2=0$.
\end{enumerate}
\end{defi}

\begin{rem} \label{rem_infcel}
For most applications the definition just given will suffice. However for an arbitrary graded ring $\gf$ (as opposed to a field) it might lead to homotopy noninvariant constructions. E.g. the Hochschild cohomology of two weakly equivalent $A_\infty$-algebras may not be isomorphic. In order to avoid troubles like this the definition needs to be modified as follows. In all three cases the $\infty$-structure $m$ can be represented as
\[ m=m_1+m_2+\ldots+m_n+\ldots, \]
where $m_i$ is a vector field of order $i$. The condition
$m^2=0$ implies that $m_1^2=0$. In other words $V$ together with the self-map $m_1$ forms a complex of $\gf$-modules. We then require that this complex be \emph{cellular}, in the sense of  \cite{krizmay}. This requirement is extraneous when for example, $V$ is finitely generated and free over $\gf$ or when $\gf$ is a field.
\end{rem}

\begin{defi} \label{def_infmor}
Let $V$ and $U$ be free graded modules:
\begin{enumerate}
\item[(a)]
Let $m$ and $m'$ be \li-structures on $V$ and $U$ respectively. An \li-morphism from $V$ to $U$ is a continuous algebra homomorphism
\[ \phi:\csalg{U} \to \csalg{V} \]
of degree zero such that $\phi \circ m'=m \circ \phi$.
\item[(b)]
Let $m$ and $m'$ be \ai-structures on $V$ and $U$ respectively. An \ai-morphism from $V$ to $U$ is a continuous algebra homomorphism
\[ \phi:\ctalg{U} \to \ctalg{V} \]
of degree zero such that $\phi \circ m'=m \circ \phi$.
\item[(c)]
Let $m$ and $m'$ be \ci-structures on $V$ and $U$ respectively. A \ci-morphism from $V$ to $U$ is a continuous algebra homomorphism
\[ \phi:\clalg{U} \to \clalg{V} \]
of degree zero such that $\phi \circ m'=m \circ \phi$.
\end{enumerate}
\end{defi}

\begin{rem} \label{rem_calseq}
We have a diagram of functors
\begin{equation} \label{fig_calseq}
\ci\mathrm{-}algebras \longrightarrow \ai\mathrm{-}algebras \longrightarrow \li\mathrm{-}algebras
\end{equation}
which depends upon the sequence of Lie algebra homomorphisms defined in Remark \ref{rem_vecseq}.

Recall that any vector field $m:\clalg{V} \to \clalg{V}$ can be uniquely extended to a continuous Hopf algebra derivation (where \ctalg{V} is equipped with the cocommutative comultiplication (\emph{shuffle coproduct}) defined in Remark \ref{rem_ccmult})
\[ m : \ctalg{V} \to \ctalg{V} \]
and that all continuous Hopf algebra derivations on \ctalg{V} are obtained from $\Der(\clalg{V})$ in this manner. It follows that any \ci-structure gives rise to an \ai-structure in this way. Similarly, any continuous Lie algebra homomorphism $\phi:\clalg{U} \to \clalg{V}$ can be uniquely extended to a continuous Hopf algebra homomorphism
\[ \phi:\ctalg{U} \to \ctalg{V} \]
and all continuous Hopf algebra homomorphisms are obtained from continuous Lie algebra homomorphisms in this way. This is because the Lie subalgebra of primitive elements of our Hopf algebra \ctalg{V} coincides with the Lie subalgebra \clalg{V}.

The category $\ai\mathrm{-}algebras$ has a subcategory whose objects consist of the \ai-algebras whose \ai-structure is also a Hopf algebra derivation and whose morphisms are the \ai-morphisms which are also Hopf algebra homomorphisms. It follows from the discussion above that this category is isomorphic to $\ci\mathrm{-}algebras$.

Recall as well that any vector field $m:\ctalg{V} \to \ctalg{V}$ can be lifted to a unique vector field
\[ \tilde{m}: \csalg{V} \to \csalg{V}. \]
It follows that any \ai-structure gives rise to a \li-structure in this manner. Similarly any continuous algebra homomorphism $\phi:\ctalg{U} \to \ctalg{V}$ can be lifted to a unique continuous homomorphism
\[ \tilde{\phi}:\csalg{U} \to \csalg{V}. \]

Observe that given a minimal \ci-algebra, the corresponding \li-algebra under the maps of diagram \eqref{fig_calseq} is a trivial algebra (i.e. its multiplication sends everything to zero).
\end{rem}

We will now define an important type of {\ai} and \ci-algebra; unital {\ai} and \ci-algebras:

\begin{defi}
Let $V$ be a free graded module:
\begin{enumerate}
\item[(i)]
\begin{enumerate}
\item
We say that an \ai-structure $m:\ctalg{V} \to \ctalg{V}$ is unital if there is a distinguished element $1 \in V$ (the unit) of degree zero which can be extended to a basis $1,\{x_i\}_{i\in I}$ of $V$ such that $m$ has the form;
\begin{equation} \label{eqn_unitform}
m=A(\boldsymbol{t})\partial_\tau+\sum_{i\in I}B_i(\boldsymbol{t})\partial_{t_i} + \ad\tau - \tau^2\partial_\tau,
\end{equation}
where $\tau,\boldsymbol{t}:=\tau,\{t_i\}_{i\in I}$ is the topological basis of $\Sigma V^*$ which is dual to the basis $\Sigma 1,\{\Sigma x_i\}_{i\in I}$ of $\Sigma V$. In this case we say that the \ai-algebra $V$ is unital or that it has a unit.
\item
Suppose that $V$ and $U$ are two unital \ai-algebras. We say that an \ai-morphism $\phi:\ctalg{U} \to \ctalg{V}$ is unital if $\phi$ has the form;
\begin{displaymath}
\begin{array}{ccc}
\phi(\tau') & = & \tau + A(\boldsymbol{t}), \\
\phi(t'_i) & = & B_i(\boldsymbol{t}); \\
\end{array}
\end{displaymath}
where $\tau,\boldsymbol{t}$ and $\tau',\boldsymbol{t}'$ are the topological bases of $\Sigma V^*$ and $\Sigma U^*$ which are dual to the bases $\Sigma 1_V,\{\Sigma x_i\}_{i\in I}$ and $\Sigma 1_U,\{\Sigma x'_j\}_{j\in J}$ of $\Sigma V$ and $\Sigma U$ respectively.
\end{enumerate}
\item[(ii)]
\begin{enumerate}
\item
We say that a \ci-structure $m:\clalg{V} \to \clalg{V}$ is unital if the corresponding \ai-structure (see Remark \ref{rem_calseq}) is unital.
\item
Suppose that $V$ and $U$ are two unital \ci-algebras. We say that a \ci-morphism $\phi:\clalg{U} \to \clalg{V}$ is unital if the corresponding \ai-morphism (see Remark \ref{rem_calseq}) is unital.
\end{enumerate}
\end{enumerate}
\end{defi}

\begin{rem} \label{rem_infdef}
There is an alternative definition of an $\infty$-structure on a free graded module $V$. This is the definition that was originally introduced by Stasheff in \cite{staha1} and \cite{staha2}. According to this definition an $\infty$-structure on $V$ is a system of maps
\[ \check{m}_i:V^{\otimes i} \to V, i \geq 1 \]
(where $|\check{m}_i|=2-i$) satisfying the higher homotopy axioms (and possibly some graded symmetry axioms as well). For \ai-algebras for instance these axioms imply that $\check{m}_2$ is associative up to homotopy (the homotopy being provided by $\check{m}_3$) whilst for \li-algebras they imply that $\check{m}_2$ satisfies the Jacobi identity up to homotopy. In particular, in all three cases they specify that the map $\check{m}_1$ is a differential and a graded derivation with respect to the multiplication $\check{m}_2$.

There is also an alternative definition of $\infty$-morphisms in this context. These are a system of maps
\[ \check{\phi}_i:V^{\otimes i} \to U, i \geq 1 \]
(where $|\check{\phi}_i|=1-i$) satisfying certain compatibility conditions with the $\check{m}_i$ and $\check{m}'_i$'s which specify the $\infty$-structures on $V$ and $U$ respectively. In particular the map $\check{\phi}_1$ must be a map of complexes; $\check{\phi}_1\circ\check{m}_1=\check{m}'_1\circ\check{\phi}_1$. We say that the $\infty$-morphism given by this system of maps is a weak equivalence if the map
\[ \check{\phi}_1: (V,\check{m}_1) \to (U,\check{m}'_1) \]
is a quasi-isomorphism. If $\check{m}_1=0$ then we say that the \ai-algebra is minimal. It is easy to see that a weak equivalence of two minimal $\infty$-algebras is in fact an isomorphism.
\end{rem}

Let us describe how this alternative style of definition is equivalent to that described in definitions \ref{def_infstr} and \ref{def_infmor} beginning with the {\ai} case. Firstly there is a one-to-one correspondence between systems of maps $\check{m}_i:V^{\otimes i} \to V, i\geq 1$ and systems of maps $m_i:\Sigma V^{\otimes i} \to \Sigma V, i\geq 1$ via the following commutative diagram:
\begin{equation} \label{fig_infdef}
\xymatrix{V^{\otimes i} \ar^{\check{m}_i}[r] & V \\ \Sigma V^{\otimes i} \ar^{(\Sigma^{-1})^{\otimes i}}[u] \ar_{m_i}[r] & \Sigma V \ar^{\Sigma^{-1}}[u]}
\end{equation}
Of course the $m_i$'s will inherit additional signs from the Koszul sign rule. It is well-known that a system of maps $m_i:\Sigma V^{\otimes i} \to \Sigma V,i\geq 1$ can be uniquely extended to a coderivation $m$ on the tensor coalgebra $T\Sigma V$ which vanishes on $\gf\subset T\Sigma V$.
 Furthermore all coderivations vanishing on {\gf} are obtained in this way, hence there is a one-to-one correspondence
\[ \Hom_\gf(T\Sigma V/\gf,\Sigma V) \leftrightarrow \{ m \in \Coder(T\Sigma V) : m(\gf) = 0 \}. \]
The condition $m^2=0$ turns out to be equivalent to the higher homotopy associativity axioms for the $\check{m}_i$'s. Now simply observe that the dual of a coderivation on $T\Sigma V$ is a continuous derivation on $(T\Sigma V)^*=\ctalg{V}$. It follows from Proposition \ref{prop_antieq} that our two definitions of an \ai-structure are equivalent.

\ai-morphisms are dealt with in a similar manner. Again there is a one-to-one correspondence between systems of maps $\check{\phi}_i:V^{\otimes i} \to U, i\geq 1$ and systems of maps $\phi_i:\Sigma V^{\otimes i} \to \Sigma U,i\geq 1$. It is well known that systems of such maps are in one-to-one correspondence with coalgebra morphisms $\phi:T\Sigma V \to T\Sigma U$. The dual of a coalgebra morphism $\phi$ is a continuous algebra homomorphism $\phi^*:\ctalg{U} \to \ctalg{V}$. The condition in Definition \ref{def_infmor} which stipulates that $\phi^*$ commutes with the \ai-structures reflects the compatibility conditions placed on the $\check{\phi}_i$'s alluded to in Remark \ref{rem_infdef}. Further details on the formal passage to the dual framework can be found in Appendix \ref{app_todual}.

The {\ci} case is a restriction of the {\ai} case. Certain graded symmetry conditions are placed on the maps $\check{m}_i:V^{\otimes i} \to V$ in addition to the higher homotopy associativity conditions. These symmetry conditions are satisfied if and only if the corresponding \ai-structure $m:\ctalg{V} \to \ctalg{V}$ is a Hopf algebra derivation (i.e. a derivation and a coderivation), where \ctalg{V} is endowed with the cocommutative comultiplication (\emph{shuffle coproduct}) defined in Remark \ref{rem_ccmult}. Similarly certain graded symmetry conditions are placed on the maps $\check{\phi}_i:V^{\otimes i} \to U$ in addition to the conditions requiring that the $\check{\phi}_i$'s are compatible with the \ai-structures. These symmetry conditions are satisfied if and only if the corresponding \ai-morphism $\phi:\ctalg{U} \to \ctalg{V}$ is a Hopf algebra homomorphism. It now follows from Remark \ref{rem_calseq} that our two definitions of a \ci-structure are equivalent.

The {\li} case is slightly different. Systems of maps $m_i:S^i\Sigma V \to \Sigma V, i\geq 1$ are in one-to-one correspondence with coderivations $m:S\Sigma V \to S\Sigma V$ which vanish on $\gf\subset S\Sigma V$, where $S\Sigma V$ is equipped with the cocommutative comultiplication defined in Remark \ref{rem_ccmult}. Similarly systems of maps $\phi_i: S^i\Sigma V \to \Sigma U,i\geq 1$ are in one-to-one correspondence with coalgebra homomorphisms $\phi:S\Sigma V \to S\Sigma U$. The condition $m^2=0$ corresponds exactly to the higher homotopy Lie axioms. The dual of a coderivation on $S\Sigma V$ which vanishes on {\gf} is a continuous derivation $m^*:(S\Sigma V)^* \to (S\Sigma V)^*$ which vanishes at zero, where $(S\Sigma V)^*$ is endowed with the shuffle product. Similarly the dual of a coalgebra morphism $\phi:S\Sigma V \to S\Sigma U$ is a continuous algebra homomorphism $\phi^*:(S\Sigma U)^* \to (S\Sigma V)^*$.

Since coinvariants are dual to invariants we have the following identity:
\[ (S\Sigma V)^*=\gf \times (\Sigma V^*)^{S_1} \times ((\Sigma V^*)^{\cotimes 2})^{S_2} \times ((\Sigma V^*)^{\cotimes 3})^{S_3} \times \ldots. \]
Recall that in equation \eqref{eqn_symiso} we defined the map $i:\csalg{V} \to (S\Sigma V)^*$ which is the canonical map identifying coinvariants with invariants and is given by the formula
\begin{equation} \label{eqn_liduin}
i(x_1 \cotimes \ldots \cotimes x_n):= \sum_{\sigma \in S_n} \sigma \cdot x_1 \cotimes \ldots \cotimes x_n.
\end{equation}
This map has an inverse $\pi:(S\Sigma V)^* \to \csalg{V}$ given by the projection
\begin{equation} \label{eqn_lidupr}
\pi(x_1 \cotimes \ldots \cotimes x_n):= \frac{1}{n!} x_1 \cotimes \ldots \cotimes x_n.
\end{equation}
Using the maps $i$ and $\pi$ to identify the module $(S\Sigma V)^*$ with the module \csalg{V}, the shuffle product on $(S\Sigma V)^*$ is transformed into the canonical commutative multiplication on \csalg{V} which is inherited from the associative multiplication on \ctalg{V}. From this it follows that our two definitions of a \li-structure are equivalent.

Further details on $\infty$-algebras and their definitions can be found in \cite{keller}, \cite[\S 5]{getjon} and \cite{markl}.

\section{Minimal Infinity-algebras}\label{minimal}

In this section we prove that any $C_\infty$-algebra is weakly equivalent to a minimal one. This theorem was first proved in the context of $A_\infty$-algebras by Kadeishvili in \cite{kadvil}, but it had its precursors in the theory of minimal models in rational homotopy theory, cf. \cite{sulvan}. It has since been reproved by many authors, cf. \cite{keller} and the references therein. The proof in the $L_\infty$ case was outlined in \cite{kontpm}. A proof in the $C_\infty$-case was given in the recent paper \cite{getzc}. We will give here a short proof based on the notion of the Maurer-Cartan moduli space associated to a differential graded Lie algebra, cf. \cite{golmil}. With obvious modifications our proof works for the $A_\infty$ and $L_\infty$ cases as well. Martin Markl pointed out to us that the minimality theorem could also
be derived from the main result of \cite{markl1} which states that structures of algebras over a cofibrant operad
(such as the $C_\infty$-operad) could be transferred across quasi-isomorphisms.

In this section we make an assumption that $\gf$ is a field, or, more generally, a graded field (i.e. a graded commutative ring whose homogeneous nonzero elements are invertible).

We first recall some standard facts from the Maurer-Cartan theory following \cite{golmil}. Let $\mathcal{G}$ be a differential graded Lie algebra which we assume to be nilpotent, or, more generally, an inverse limit of nilpotent differential graded Lie algebras. The set $\mathcal{MC}(\mathcal{G})\subset \mathcal{G}^1$ is by definition the set of solutions of the Maurer-Cartan equation
\begin{equation} \label{MC}
d\gamma+\frac{1}{2}[\gamma,\gamma]=0.
\end{equation}
Furthermore, the Lie algebra $\mathcal{G}^0$ acts on $\mathcal{MC}(\mathcal{G})$ by infinitesimal affine transformations: for $\alpha\in \mathcal{G}^0$ we have a vector field $\gamma \mapsto d\alpha+[\gamma,\alpha]$ on  $\mathcal{MC}(\mathcal{G})$. This action exponentiates to an action of the Lie group $\exp(\mathcal{G}^0)$ on $\mathcal{MC}(\mathcal{G})$. We will call the set of orbits with respect to this action \emph{the Maurer-Cartan moduli space} associated to the differential graded Lie algebra $\mathcal{G}$.

Now let $\mathcal{G}_1,\mathcal{G}_2$ be two nilpotent differential graded Lie algebras. We assume that they are endowed with finite filtrations $\{F_p(\mathcal{G}_1)\}$ and $\{F_p(\mathcal{G}_2)\}$, $p=1,2,\ldots, n$ and there is a map of filtered differential graded Lie algebras $f:\mathcal{G}_1\rightarrow \mathcal{G}_2$ which induces a quasi-isomorphism on the associated graded of $\mathcal{G}_1$ and $\mathcal{G}_2$. Under these assumptions we have the following result.

\begin{theorem} \label{getgm}
The map $f$ induces a bijection $\mathcal{MC}(\mathcal{G}_1)\rightarrow \mathcal{MC}(\mathcal{G}_2)$.
\end{theorem}

 This seems to be a well known result and is formulated in this form in \cite{getzlr}, Theorem 2.1. However Goldman-Millson's version \cite{golmil}, Theorem 2.4 does not readily carry over since we are dealing with differential graded Lie algebras which are not necessarily concentrated in nonnegative degrees. We therefore sketch a proof suitable for this more general situation. This proof is modelled on \cite{getzs}, Proposition 4.6. We start with the definition of the simplicial Maurer-Cartan set (also called the simplicial Deligne groupoid) associated to a differential graded Lie algebra $\mathcal{G}$. A more detailed discussion could be found in e.g. \cite{hinich}.  Let $\Omega_n$ be the algebra of polynomial differential forms on the standard $n$-simplex; the collection $\Omega_\bullet=\{\Omega_n\}_{n=0}^\infty$ forms a commutative simplicial differential graded algebra.  Note that $\mathcal{G}\otimes \Omega_n$ has the structure of a differential graded Lie algebra.

\begin{defi}
For $n\geq 0$ set $\mathcal{MC}_n(\mathcal{G}):= \mathcal{MC}(\mathcal{G}\otimes \Omega_n)$. The simplicial structure on $\Omega_\bullet$ determines the structure of a simplicial set on  $\mathcal{MC}_\bullet(\mathcal{G}):=\{\mathcal{MC}_n(\mathcal{G})\}_{n=0}^\infty$ which will be called \emph{the Maurer-Cartan simplicial set} associated with $\mathcal{G}$.
\end{defi}

The main result of \cite{sstas} implies that there is a one-to-one correspondence between $\mathcal{MC}(\mathcal{G})$ and $\pi_0(\mathcal{MC}_\bullet(\mathcal{G}))$, the set of connected components of the simplicial set $\mathcal{MC}_\bullet(\mathcal{G})$. Furthermore, for a surjective map of differential graded Lie algebras $\mathcal{G}\rightarrow \mathcal{G}^\prime$ the induced map $\mathcal{MC}_\bullet(\mathcal{G})\rightarrow \mathcal{MC}_\bullet(\mathcal{G}^\prime)$ is a fibration of simplicial sets. This is proved (in a more general context for $L_\infty$-algebras) in \cite{getzs}.

Finally, using induction up the filtrations of $\mathcal{G}_1$ and $\mathcal{G}_2$ and comparing the associated towers of fibrations of simplicial sets we show that the simplicial sets $\mathcal{MC}(\mathcal{G}_1)$ and $\mathcal{MC}(\mathcal{G}_2)$ are weakly equivalent. In particular, their sets of connected components are in one-to-one correspondence. This finishes our sketch proof of Theorem \ref{getgm}.

Now let $V$ be a graded vector space and $m$ be a \ci-structure on $V$.  Then $m$ is a vector field on the Lie algebra $\widehat{L}(\Sigma V^*)$. We will denote $\Sigma V^*$ by $W$. Choose a topological basis $\{x_i\}_{i \in I}$ in $W$ and denote by $\mathcal{G}$ the Lie algebra of vector fields having the form $\sum_{i\in I} f_i\partial_{x_i}$, where the $f_i$'s are (possibly uncountably infinite) sums of Lie monomials in $x_i$ of order 2 or higher. Clearly, this definition of $\mathcal{G}$ does not depend on the choice of a basis and $\mathcal{G}$ is a formal Lie algebra. The Lie group $\exp(\mathcal{G}^0)$ is the group of formal Lie series whose linear term is the identity linear transformation.

We will view $\mathcal{G}$ as a Lie DGA with respect to the operator $m_1$, the linear part of the vector field $m$. Denote by $Z^\bullet(W)$ and $H^\bullet(W)$ respectively, the cocycles and cohomology of $W$ with respect to $m_1$.

 Note that $\mathcal{G}$ has a filtration given by the order of a vector field. This filtration is in fact a grading and $m_1$ preserves this grading since a commutator of a linear vector filed and a vector field of order $n$ is again a vector field of order $n$. Therefore the Lie algebra $H^\bullet(\mathcal{G})$ is bigraded.

Furthermore, direct inspection shows that $m-m_1$ is a Maurer-Cartan element in $\mathcal{G}$, i.e. it satisfies  equation \eqref{MC}. Moreover, two Maurer-Cartan elements are  equivalent with respect to the action of $\exp(\mathcal{G}^0)$ if and only if the corresponding \ci-structures are weakly equivalent through a formal diffeomorphism $f=(f_1,f_2,\ldots)$ of $\widehat{L}(W)$ whose linear part $f_1$ is the identity map.   Following \cite{polchu} we will call such equivalences \emph{strict $C_\infty$-isomorphisms}.

Let us now fix an integer $n>1$ and denote by $\mathcal{G}_n$ the quotient of $\mathcal{G}$ by the ideal of vector fields having order $>n$. Using the geometric langauge we say that $\mathcal{G}_n$ is the group of germs of formal diffeomorphisms of order $n$. Clearly $\mathcal{G}_n$ is a nilpotent differential graded algebra and $\mathcal{G}=\inlim{n} {\mathcal{G}_n}$.

\begin{prop} \label{form}
The  Lie DGA $\mathcal{G}$  is quasi-isomorphic to its cohomology $H^\bullet(\mathcal{G})$ which is considered as a differential graded Lie algebra with trivial differential.
\end{prop}

\begin{rem}
The above proposition says that $\mathcal{G}$ is \emph{formal} in the sense of rational homotopy theory. We refrain from such a formulation, however, since the term `formal' has a different meaning for us.
\end{rem}

\begin{proof}
Let $Z^\bullet(W)$ be the space of $m_1$-cocycles in $W$ and choose a section of the projection map $Z^\bullet(W)\rightarrow H^\bullet(W)$. Then $H^\bullet(W)$ could be considered as a subspace in $Z^\bullet(W)$ and hence - also in $W$ itself.  Next, choose a section $j:W\rightarrow H^\bullet(W)$
of the embedding $i:H^\bullet(W)\hookrightarrow W$.

Denote by $\mathcal{H}$ the Lie algebra of vector fields on the Lie algebra $\widehat{L}(H^\bullet(W))$ which have order $\geq 2$.
 Observe that we have a natural isomorphism of vector  spaces
\[{\mathcal G}\cong \Hom(W, \hat{L}(W))\]
and similarly
\[{\mathcal H} \cong \Hom(H^\bullet(W),\hat{L}(H^\bullet(W)).\]
Using this identification, we define the maps $\tilde{i}:\mathcal{H} \rightarrow \mathcal{G}$ and $\tilde{j}:\mathcal{G}\rightarrow \mathcal{H}$ as follows: for $f\in \Hom(H^\bullet(W),\hat{L}^n(H^\bullet(W))$, $g \in \Hom(W, \hat{L}^n(W))$, $a\in W$, $b\in H^\bullet(W)$  set
\[[\tilde{i}(f)](a)=\hat{L}^n(i)\circ f\circ j(a),\]
\[[\tilde{j}(g)](b)= \hat{L}^n(j)   \circ g\circ i(b).\]
Here we denoted by $\widehat{L}^n, n=2,3,\ldots$ the functor which associates to a vector space $X$ the space generated by the Lie monomials of length $n$ inside $\widehat{L}(X)$. Note that only $\tilde{i}$ is a map of Lie algebras.

The functor $L$ associating to a graded vector space the free Lie algebra on it commutes with cohomology by \cite{quillen}, Appendix B, Proposition 2.1 (this is a consequence of the Poincar\'e-Birkhoff-Witt theorem). Since the inverse limit functor is exact on the category of finite dimensional vector spaces we conclude that $H^\bullet(\widehat{L}(W))\cong \widehat{L}(H^\bullet(W))$ and it follows that the maps $\tilde{i},\tilde{j}$ are quasi-isomorphisms as required.
\end{proof}

\begin{rem}
Note that our proof of Proposition \ref{form} actually gives slightly more than claimed. First, our proof shows the associated graded to $\mathcal{G}$ is quasi-isomorphic to its own cohomology.
(This holds simply because the filtration on $\mathcal{G}$ is in fact a second grading.)

Second, denote by $\mathcal{H}_n$ the quotient of $\mathcal{H}$ by the ideal of vector fields of order $>n$. Then the restrictions of `formality maps' $\tilde{i}$ and $\tilde{j}$ determine  mutually (quasi-)inverse quasi-isomorphisms between $\mathcal{G}_n$ and $\mathcal{H}_n$. Invoking
again the analogy with rational homotopy theory one can express this by saying that the pronilpotent differential graded algebra $\mathcal{G}$ is \emph{continuously formal}. This is an important technical point which is necessary for the proof of the minimality theorem.
\end{rem}

\begin{cor} \label{kadeish}
(Minimality theorem) Let $(V,m_V)$ be a \ci-algebra. Then there exists a minimal \ci-algebra $(U, m_U)$ and a weak equivalence of \ci-algebras $(U, m_U)\rightarrow (V,m_V)$.
\end{cor}

\begin{proof}
Set $U:=H^\bullet(V)$ and $W:=\Sigma V^*$. Denote, as before, by $\mathcal{G}$ the Lie algebra of vector fields on $\widehat{L}(W)$ of order $\geq 2$ and similarly denote by $\mathcal{H}$ the Lie algebra of vector fields on $\widehat{L}(H^\bullet(W))$ whose order is $\geq 2$. Choosing a basis for representatives of cohomology classes, we will regard $H^\bullet(W)$ as a subspace in $W$. A choice of a complement will determine a map of differential graded Lie algebras $\tilde{i}: \mathcal{H}\rightarrow \mathcal{G}$ which is a quasi-isomorphism by Theorem \ref{form}. It suffices to show that $\tilde{i}$ induces a bijection on the Maurer-Cartan moduli spaces
\[\mathcal{MC}(\mathcal{H})/\exp(\mathcal{H}^0) \overset{\cong}{\to} \mathcal{MC}(\mathcal{G})/\exp(\mathcal{G}^0).\]
Indeed, that would mean that any \ci-structure on $V$ could be reduced to a minimal one using a composition of a strict \ci-isomorphism and a linear projection $V\rightarrow U=H^\bullet(V)$.

To get the desired isomorphism note that the tower of differential graded Lie algebras $\mathcal{H}_2\leftarrow \mathcal{H}_3\leftarrow\ldots$ determines a tower of fibrations of simplicial sets $\mathcal{MC}_\bullet(\mathcal{H}_2)\leftarrow \mathcal{MC}_\bullet(\mathcal{H}_3)\ldots$. Similarly we have the tower of simplicial sets $\mathcal{MC}_\bullet(\mathcal{G}_2)\leftarrow \mathcal{MC}_\bullet(\mathcal{G}_3)\ldots$ associated with the tower of Lie algebras $\mathcal{G}_2\leftarrow \mathcal{G}_3\leftarrow\ldots$. Since $\mathcal{G}_n$ is quasi-isomorphic to $\mathcal{H}_n$ for all $n$ we conclude that these towers of fibrations are level-wise weakly equivalent. There are isomorphisms
\[ \inlim{n}{\mathcal{MC}_\bullet(\mathcal{G}_n)} \cong \mathcal{MC}_\bullet(\mathcal{G}), \]
\[ \inlim{n}{\mathcal{MC}_\bullet(\mathcal{H}_n)} \cong \mathcal{MC}_\bullet(\mathcal{H}) \]
and it follows that $\mathcal{MC}_\bullet(\mathcal{G})$ and $\mathcal{MC}_\bullet(\mathcal{H})$ are weakly equivalent simplicial sets. In particular, their sets of connected components are in one-to-one correspondence.
\end{proof}

\begin{rem}
Theorem \ref{form} and Corollary \ref{kadeish} extend in the context of $A_\infty$ and  \li-algebras. For \li-algebras the reference to the Poincar\'e-Birkhoff-Witt theorem is replaced by the fact that the functor of $S_n$-coinvariants is exact and therefore commutes with cohomology. For $A_\infty$-algebras the corresponding issue never arises and in fact the minimality theorem in the $A_\infty$ context is true in arbitrary characteristic.
\end{rem}

\section{The Cohomology of Infinity-algebras} \label{sec_iachom}

In this section we will define the various cohomology theories for \li, {\ai} and \ci-algebras that we will use throughout the rest of our paper. There will be quite a number of different cohomology theories to be defined. For \ai-algebras it will be useful to define additional quasi-isomorphic complexes computing the cohomology of the \ai-algebras. For example, this will allow us to describe a periodicity exact sequence for \ai-algebra cohomology. We will also prove some other basic facts about \ai-algebra cohomology which are the infinity-analogues of familiar results for strictly associative graded algebras.

\subsection{Hochschild, Harrison and Chevalley-Eilenberg theories for $\infty$-algebras}

We will begin by defining the cohomology theories which do not involve the action of the cyclic groups. These are the theories that control the deformations of $\infty$-algebras, cf. \cite{hamilt}, \cite{pensch} and \cite{penfia}.

\begin{defi}
Let $V$ be a free graded module:
\begin{enumerate}
\item[(a)]
Let $m:\csalg{V} \to \csalg{V}$ be an \li-structure. The Chevalley-Eilenberg complex of the \li-algebra $V$ with coefficients in $V$ is defined on the module consisting of all vector fields on \csalg{V}:
\[ \clac{V}{V}:=\Sigma^{-1}\Der(\csalg{V}). \]
The differential $d:\clac{V}{V} \to \clac{V}{V}$ is given by
\[ d(\xi):=[m,\xi], \quad \xi \in \Der(\csalg{V}). \]
The Chevalley-Eilenberg cohomology of $V$ with coefficients in $V$ is defined as the cohomology of the complex \clac{V}{V} and denoted by \hlac{V}{V}.
\item[(b)]
Let $m:\ctalg{V} \to \ctalg{V}$ be an \ai-structure. The Hochschild complex of the \ai-algebra $V$ with coefficients in $V$ is defined on the module consisting of all vector fields on \ctalg{V}:
\[ \choch{V}{V}:=\Sigma^{-1}\Der(\ctalg{V}). \]
The differential $d:\choch{V}{V} \to \choch{V}{V}$ is given by
\[ d(\xi):=[m,\xi], \quad \xi \in \Der(\ctalg{V}). \]
The Hochschild cohomology of $V$ with coefficients in $V$ is defined as the cohomology of the complex \choch{V}{V} and denoted by \hhoch{V}{V}.

\item[(c)]
Let $m:\clalg{V} \to \clalg{V}$ be a \ci-structure. The Harrison complex of the \ci-algebra $V$ with coefficients in $V$ is defined on the module consisting of all vector fields on \clalg{V}:
\[ \caq{V}{V}:=\Sigma^{-1}\Der(\clalg{V}). \]
The differential $d:\caq{V}{V} \to \caq{V}{V}$ is given by
\[ d(\xi):=[m,\xi], \quad \xi \in \Der(\clalg{V}). \]
The Harrison cohomology of $V$ with coefficients in $V$ is defined as the cohomology of the complex \caq{V}{V} and denoted by \haq{V}{V}.
\end{enumerate}
\end{defi}

\begin{rem}
The purpose of the desuspensions in the above definition is to make the grading consistent with the classical grading on these cohomology theories.
\end{rem}

\begin{rem}
The Chevalley-Eilenberg, Hochschild and Harrison complexes all have the structure of a differential graded Lie algebra under the commutator bracket $[-,-]$. It was Stasheff who realised in \cite{stasgb} that on the Hochschild cohomology of a strictly associative algebra this was the Gerstenhaber bracket introduced by Gerstenhaber in \cite{gerbra}. Stasheff considered the \ai-structure as a coderivation $m:T\Sigma V \to T\Sigma V$ and defined his bracket on $\Coder(T\Sigma V)$ as the commutator of coderivations. If we were to translate this structure directly to $\Der(\ctalg{V})$ via the contravariant functor of Proposition \ref{prop_antieq} then the induced Lie bracket $\{-,-\}$ would be the reverse of the commutator bracket
\[ \{\xi,\gamma\} = [\gamma,\xi]; \quad \xi,\gamma \in \Der(\ctalg{V}) \]
and the differential would be the map
\[ \xi \mapsto [\xi,m]; \quad \xi \in \Der(\ctalg{V}), \]
however, there is an isomorphism between this differential graded Lie structure and the differential graded Lie structure we defined on \choch{V}{V} above given by the map
\[ \xi \mapsto (-1)^{|\xi|(|\xi|+1)(2|\xi|+1)/6 + 1}\xi, \quad \xi \in \Der(\ctalg{V}). \]
\end{rem}

Now let us describe the cohomology theories which are dual to the corresponding cohomology theories defined above.

\begin{defi}
Let $V$ be a free graded module:
\begin{enumerate}
\item[(a)]
Let $m:\csalg{V} \to \csalg{V}$ be an \li-structure. The Chevalley-Eilenberg complex of the \li-algebra $V$ with coefficients in $V^*$ is defined on the module consisting of all 1-forms:
\[ \clac{V}{V^*}:=\Sigma\drof[Com]{\csalg{V}}. \]
The differential on this complex is the (suspension of the) Lie operator of the vector field $m$;
\[ L_m:\drof[Com]{\csalg{V}} \to \drof[Com]{\csalg{V}}. \]
The Chevalley-Eilenberg cohomology of $V$ with coefficients in $V^*$ is defined as the cohomology of the complex \clac{V}{V^*} and denoted by \hlac{V}{V^*}.
\item[(b)]
Let $m:\ctalg{V} \to \ctalg{V}$ be an \ai-structure. The Hochschild complex of the \ai-algebra $V$ with coefficients in $V^*$ is defined on the module consisting of all 1-forms:
\[ \choch{V}{V^*}:=\Sigma\drof[Ass]{\ctalg{V}}. \]
The differential on this complex is the (suspension of the) Lie operator of the vector field $m$;
\[ L_m:\drof[Ass]{\ctalg{V}} \to \drof[Ass]{\ctalg{V}}. \]
The Hochschild cohomology of $V$ with coefficients in $V^*$ is defined as the cohomology of the complex \choch{V}{V^*} and denoted by \hhoch{V}{V^*}.
\item[(c)]
Let $m:\clalg{V} \to \clalg{V}$ be a \ci-structure. The Harrison complex of the \ci-algebra $V$ with coefficients in $V^*$ is defined on the module consisting of all 1-forms:
\[ \caq{V}{V^*}:=\Sigma\drof[Lie]{\clalg{V}}. \]
The differential on this complex is the (suspension of the) Lie operator of the vector field $m$;
\[ L_m:\drof[Lie]{\clalg{V}} \to \drof[Lie]{\clalg{V}}. \]
The Harrison cohomology of $V$ with coefficients in $V^*$ is defined as the cohomology of the complex \caq{V}{V^*} and denoted by \haq{V}{V^*}.
\end{enumerate}
\end{defi}

\begin{rem}
Again the suspensions appear in order to keep the grading consistent with the classical grading on these cohomology theories.
\end{rem}

\begin{rem}
We did not prove that $L_m$ is indeed a differential, however this is obvious from Lemma \ref{lem_schf}:
\[ L_m^2=\frac{1}{2}[L_m,L_m]=L_{\frac{1}{2}[m,m]}=L_{m^2}=0. \]
\end{rem}

\begin{rem}
Our definition of the Harrison complex implies that in the case of an ungraded commutative algebra the
zeroth term of the complex is not present. This terminology differs slightly from what seems to be adopted in the modern literature (e.g. \cite{loday}) in that in the latter the algebra itself is considered as the zeroth term (note, however, that in his original paper \cite{harrison} Harrison only defined the 1st, 2nd and 3rd cohomology  groups). This distinction is unimportant since the zeroth term always splits off as a direct summand.
\end{rem}

Next we will define the bar cohomology of an \ai-algebra:

\begin{defi} \label{def_barhom}
Let $V$ be an \ai-algebra with \ai-structure $m:\ctalg{V} \to \ctalg{V}$. We define the map $b'$ as the restriction of (the suspension of) $m$ to $\cbr{V}:=\Sigma\left[\col{V}\right]$:
\begin{equation} \label{fig_barhom}
\xymatrix{\col{V} \ar^{b'}[d] \ar@<-0.5ex>@{^(->}[r] & \ctalg{V} \ar^{m}[d] \\ \col{V} \ar@<-0.5ex>@{^(->}[r] & \ctalg{V}}
\end{equation}
The map $b'$ is a differential on \cbr{V} and we define the bar cohomology of the \ai-algebra $V$ as the cohomology of the bar complex \cbr{V} and denote it by \hbr{V}.
\end{defi}

\begin{rem} \label{rem_chomid}
We can also define a map $b:\col{V} \to \col{V}$ by identifying \col{V} with the underlying module of the complex \choch{V}{V^*}. We do this via the module isomorphism \thiso[Ass] defined by Lemma \ref{lem_ofisom}:
\begin{equation} \label{fig_chomid}
\xymatrix{\col{V} \ar^{b}[d] \ar^{\thiso[Ass]}[r] & \drof[Ass]{\ctalg{V}} \ar^{L_m}[d] \\ \col{V} \ar^{\thiso[Ass]}[r] & \drof[Ass]{\ctalg{V}}}
\end{equation}
The complex whose underlying module is $\Sigma\left[\col{V}\right]$ and whose differential is (the suspension of) $b$ is by definition isomorphic to the Hochschild complex of $V$ with coefficients in $V^*$.
\end{rem}

Now we shall show that the bar cohomology of a unital \ai-algebra is trivial.

\begin{lemma} \label{lem_unitch}
Let $V$ be a unital \ai-algebra with unit $1 \in V$. There is a contracting homotopy $h:\cbr{V} \to \cbr{V}$ which is the dual of the map
\begin{displaymath}
\begin{array}{ccc}
\bigoplus_{i=1}^\infty \Sigma V^{\otimes i} & \to & \bigoplus_{i=1}^\infty \Sigma V^{\otimes i}, \\
x & \mapsto & 1 \otimes x. \\
\end{array}
\end{displaymath}
\end{lemma}

\begin{proof}
By Definition \ref{def_infstr} part (b) the \ai-structure $m:\fpsa \to \fpsa$ has the form
\[ m=A(\boldsymbol{t})\partial_\tau+\sum_{i\in I}B_i(\boldsymbol{t})\partial_{t_i} + \ad\tau - \tau^2\partial_\tau, \]
where $\tau$ is dual to the unit $1 \in V$. We have the following formula for the map $h$;
\begin{displaymath}
\begin{array}{cc}
h(\tau x) = x, & x \in \fpsa; \\
h(t_i x) = 0, & x \in \fpsa. \\
\end{array}
\end{displaymath}
We calculate
\begin{displaymath}
\begin{split}
b'h(t_i x)+hb'(t_i x) & = 0 + h[([\tau,t_i]+B_i(\boldsymbol{t}))x] \pm h(t_ib'(x)), \\
& = h(\tau t_i x)=t_i x. \\
b'h(\tau x)+hb'(\tau x) & = b'(x)+h(\tau^2 x) - h(\tau b'(x)), \\
& =b'(x)+ \tau x -b'(x) = \tau x. \\
\end{split}
\end{displaymath}
Similar calculations show that
\[ b'h(t_i)+hb'(t_i)= t_i \quad \text{and} \quad b'h(\tau)+hb'(\tau)= \tau. \]
therefore
\[ b'h+hb'=\id. \]
\end{proof}

Let $V$ be a unital \ai-algebra and let $\tau,\boldsymbol{t}$ be a topological basis of $\Sigma V^*$ where $\tau$ is dual to the unit $1 \in V$, then $\ctalg{V} = \fpsa$. We say a 1-form $\alpha \in \drof[Ass]{\ctalg{V}}$ is normalised if it is a linear combination of elements of the form $q\cdot dv$ where $v \in \Sigma V^*$ and $q \in \gf\langle\langle \boldsymbol{t} \rangle \rangle$;
\begin{equation} \label{eqn_normform}
\alpha = A(\boldsymbol{t})d\tau + \sum_{i \in I} B_i(\boldsymbol{t})dt_i; \quad A(\boldsymbol{t}), B_i(\boldsymbol{t}) \in \gf\langle\langle \boldsymbol{t} \rangle \rangle.
\end{equation}
We will denote the module of normalised 1-forms by \drnof[Ass]{\ctalg{V}}. It is clear to see that the map \thiso[Ass] defined by Lemma \ref{lem_ofisom} identifies the module of normalised 1-forms with the module $\Sigma V^* \cotimes \widehat{T}(\Sigma V/\gf)^*$.

\begin{prop} \label{prop_htpret}
Let $V$ be a unital \ai-algebra with \ai-structure $m:\ctalg{V} \to \ctalg{V}$. The normalised 1-forms $\Sigma\drnof[Ass]{\ctalg{V}}$ form a subcomplex of $\choch{V}{V^*}:= \left(\Sigma\drof[Ass]{\ctalg{V}}, L_m\right)$. Furthermore the subcomplex of normalised 1-forms is a chain deformation retract of \choch{V}{V^*}.
\end{prop}

\begin{proof}
By Definition \ref{def_infstr} the \ai-structure $m$ has the form
\[ m=A(\boldsymbol{t})\partial_\tau+\sum_{i\in I}B_i(\boldsymbol{t})\partial_{t_i} + \ad\tau - \tau^2\partial_\tau. \]
We calculate that for all $q \in \gf\langle\langle \boldsymbol{t} \rangle \rangle$ and $i \in I$;
\begin{displaymath}
\begin{array}{rcll}
L_m(q\cdot dt_i) & = & [\tau,q]\cdot dt_i +(-1)^{|q|+1}q\cdot d[\tau,t_i] & \mod \drnof[Ass]{\ctalg{V}}, \\
& = & [\tau,q]\cdot dt_i +(-1)^{|q|+1}q\cdot([d\tau,t_i]-[\tau,dt_i]) & \mod \drnof[Ass]{\ctalg{V}}, \\
& = & [\tau,q]\cdot dt_i +(-1)^{|q|}[q,t_i]\cdot d\tau +(-1)^{|q|}[q,\tau]\cdot dt_i & \mod \drnof[Ass]{\ctalg{V}}, \\
& = & (-1)^{|q|}[q,t_i]\cdot d\tau = 0 & \mod \drnof[Ass]{\ctalg{V}}. \\
\end{array}
\end{displaymath}
Similar calculations show that $L_m(q\cdot d\tau) \in \drnof[Ass]{\ctalg{V}}$ for all $q \in \gf\langle\langle \boldsymbol{t} \rangle \rangle$, hence \drnof[Ass]{\ctalg{V}} is a subcomplex of \choch{V}{V^*}.

Let $x_0,\ldots,x_n \in \Sigma V^*$, we say that the 1-form
\begin{equation} \label{eqn_htpretdummy}
x_1\ldots x_n\cdot dx_0
\end{equation}
is $i$-normalised if $x_1\ldots x_i \in \gf\langle\langle \boldsymbol{t} \rangle\rangle$. A generic $i$-normalised 1-form is a linear combination of $i$-normalised elements of the form \eqref{eqn_htpretdummy}. Obviously \eqref{eqn_htpretdummy} is $n$-normalised if and only if it is normalised.

Let $\gamma$ be the constant vector field $\gamma := 1\partial_\tau$. Define the map
\[ s_i:\drof[Ass]{\ctalg{V}} \to \drof[Ass]{\ctalg{V}} \]
by the formula
\[ s_i(x_1\ldots x_i\cdot x_{i+1} \cdot q\cdot dx_0) := (-1)^{|x_1|+\ldots+|x_i|}\gamma(x_{i+1})x_1\ldots x_i \cdot q\cdot dx_0, \]
where $x_0,\ldots,x_{i+1} \in \Sigma V^*$ and $q \in \fpsa$. We define the map
\[ h_i:\drof[Ass]{\ctalg{V}} \to \drof[Ass]{\ctalg{V}} \]
as $h_i:= \id + L_m s_i + s_i L_m$ which is of course chain homotopic to the identity.

We claim that $h_i$ takes $i$-normalised 1-forms to $i+1$-normalised 1-forms. Let $x,v \in \Sigma V^*,q \in \fpsa$ and let $p \in \gf\langle\langle \boldsymbol{t} \rangle\rangle$ be a power series of order $i$ so that $\alpha:=pxq\cdot dv$ is an $i$-normalised 1-form. Let $m_1$ be the linear part of the \ai-structure $m$. We will now calculate $h_i(\alpha)$:
\begin{equation} \label{eqn_htpretdummya}
\begin{split}
h_i(\alpha) = & \alpha + L_m((-1)^{|p|}\gamma(x)pq\cdot dv) + s_i(m(p)xq\cdot dv + (-1)^{|p|}pm(x)q\cdot dv \\
& + (-1)^{|p|+|x|}pxm(q)\cdot dv +(-1)^{|p|+|x|+|q|+1}pxq\cdot dm(v)), \\
= & \alpha + (-1)^{|p|}\gamma(x)[m(p)q\cdot dv + (-1)^{|p|}pm(q)\cdot dv + (-1)^{|p|+|q|+1}pq\cdot dm(v)] \\
& + s_i([\tau,p]q\cdot dv) + (-1)^{|p|+1}\gamma(x)m_1(p)q\cdot dv + (-1)^{|p|}s_i(pm(x)q\cdot dv) \\
& + (-1)^{|x|}\gamma(x)pm(q)\cdot dv + (-1)^{|p|+|x|+|q|+1}s_i(pxq\cdot dm(v)), \\
= & \alpha + (-1)^{|p|}\gamma(x)m(p)q\cdot dv + \gamma(x)pm(q)\cdot dv +(-1)^{|q|+1}\gamma(x)pq\cdot dm(v) - pxq\cdot dv \\
& + (-1)^{|p|+1}\gamma(x)m_1(p)q\cdot dv + (-1)^{|p|}s_i(pm(x)q\cdot dv) \\
& + (-1)^{|x|}\gamma(x)pm(q)\cdot dv + (-1)^{|p|+|x|+|q|+1}s_i(pxq\cdot dm(v)), \\
= & (-1)^{|p|}\gamma(x)m(p)q\cdot dv + (-1)^{|q|+1}\gamma(x)pq\cdot dm(v) + (-1)^{|p|+1}\gamma(x)m_1(p)q\cdot dv \\
& + (-1)^{|p|}s_i(pm(x)q\cdot dv) + (-1)^{|p|+|x|+|q|+1}s_i(pxq\cdot dm(v)). \\
\end{split}
\end{equation}

Suppose that $x,v \in \boldsymbol{t}:=\{t_i\}_{i \in I}$, then \eqref{eqn_htpretdummya} implies that
\begin{equation} \label{eqn_htpretdummyb}
\begin{split}
h_i(\alpha) & = (-1)^{|p|}s_i(p[\tau,x]q\cdot dv) + (-1)^{|p|+|x|+|q|+1}s_i(pxq\cdot d[\tau,v]), \\
& = \alpha + (-1)^{|p|+|x|+|q|}s_i([pxq,v]\cdot d\tau + [pxq,\tau]\cdot dv), \\
& = \alpha. \\
\end{split}
\end{equation}
A similar calculation shows that $h_i(\alpha)=\alpha$ when $x \in \boldsymbol{t}$ and $v = \tau$. This means that $h_i$ acts as the identity on $i+1$-normalised 1-forms.

Now suppose that $x = \tau$ and $v \in \boldsymbol{t}$, then
\begin{equation} \label{eqn_htpretdummyc}
\begin{split}
h_i(\alpha) = & (-1)^{|p|}[\tau,p]q\cdot dv + (-1)^{|p|}m_1(p)q\cdot dv + (-1)^{|q|+1}pq\cdot dm(v) + (-1)^{|p|+1}m_1(p)q\cdot dv \\
& + (-1)^{|p|}s_i(p\tau^2q\cdot dv) + (-1)^{|p|+|q|}s_i(p\tau q\cdot dm(v)) + (\text{$i+1$-normalised 1-forms}), \\
= & (-1)^{|p|}\tau pq\cdot dv - \alpha + (-1)^{|q|+1}pq\cdot dm(v) + \alpha \\
& + (-1)^{|p|+|q|}s_i(p\tau q\cdot dm(v)) + (\text{$i+1$-normalised 1-forms}), \\
= & (-1)^{|p|}\tau pq\cdot dv + (-1)^{|q|+1}pq\cdot d[\tau,v] \\
& + (-1)^{|p|+|q|}s_i(p\tau q\cdot d[\tau,v]) + (\text{$i+1$-normalised 1-forms}), \\
= & (-1)^{|p|}\tau pq\cdot dv + (-1)^{|q|}pqv\cdot d\tau + (-1)^{|q|}[pq,\tau]\cdot dv \\
& + (-1)^{|q|+1}pqv\cdot d\tau + (-1)^{|q|+1}pq\tau\cdot dv + (\text{$i+1$-normalised 1-forms}), \\
= & 0 + (\text{$i+1$-normalised 1-forms}). \\
\end{split}
\end{equation}
A similar calculation shows that
\[ h_i(\alpha) = 0 \mod (\text{$i+1$-normalised 1-forms}) \]
when $x = \tau$ and $v = \tau$. This means that $h_i$ takes $i$-normalised forms to $i+1$-normalised forms as claimed.

We define the map $H:\drof[Ass]{\ctalg{V}} \to \drnof[Ass]{\ctalg{V}}$ as
\[ H:= \ldots \circ h_n\circ\ldots\circ h_2\circ h_1. \]
By the definition of the $h_n$'s the map $H$ is homotopic to the identity. Equations \eqref{eqn_htpretdummya}, \eqref{eqn_htpretdummyb} and \eqref{eqn_htpretdummyc} show that $H$ takes 1-forms to normalised 1-forms so $H$ is a well defined map. Lastly equation \eqref{eqn_htpretdummyb} shows that $H$ splits the inclusion
\[ \drnof[Ass]{\ctalg{V}} \hookrightarrow \drof[Ass]{\ctalg{V}}, \]
so $H$ is the chain homotopy retraction exhibiting \drnof[Ass]{\ctalg{V}} as a chain homotopy retract of \drof[Ass]{\ctalg{V}} that we sought.
\end{proof}

\subsection{Cyclic cohomology theories}

We will now define the cyclic cohomology theories for \li, {\ai} and \ci-algebras. The cyclic
theory for \ai-algebras was first defined in \cite{getjo}. In \cite{pensch} Penkava and Schwarz showed that cyclic cohomology controlled the deformations of the \ai-algebra which preserve a fixed invariant inner product. Our approach to defining cyclic cohomology will be different however, in that we will use the framework of noncommutative geometry.

Let $V$ be a profinite graded module. We say $q(\boldsymbol{t}) \in \ctalg{V}$ vanishes at zero if it belongs to the ideal $\Sigma V^* \cdot \ctalg{V}$ of \ctalg{V}, i.e. the power series $q(\boldsymbol{t})$ has vanishing constant term.

\begin{defi}
Let $V$ be a free graded module:
\begin{enumerate}
\item[(a)]
Let $m:\csalg{V} \to \csalg{V}$ be an \li-structure. The cyclic Chevalley-Eilenberg complex of the \li-algebra $V$ is defined on the module of 0-forms vanishing at zero;
\[ \cclac{V}:=\Sigma\left[\{q \in \drzf[Com]{\csalg{V}}:\ q \text{ vanishes at zero} \}\right]. \]
The differential on this complex is the restriction of the (suspension of the) Lie operator of the vector field $m$;
\[ L_m:\drzf[Com]{\csalg{V}} \to \drzf[Com]{\csalg{V}} \]
to \cclac{V}. The cyclic Chevalley-Eilenberg cohomology of $V$ is defined as the cohomology of the complex \cclac{V} and is denoted by \hclac{V}.
\item[(b)]
Let $m:\ctalg{V} \to \ctalg{V}$ be an \ai-structure. The cyclic Hochschild complex of the \ai-algebra $V$ is defined on the module of 0-forms vanishing at zero:
\[ \cchoch{V}:=\Sigma\left[\{q \in \drzf[Ass]{\ctalg{V}}:\ q \text{ vanishes at zero} \}\right]. \]
The differential on this complex is the restriction of the (suspension of the) Lie operator of the vector field $m$;
\[ L_m:\drzf[Ass]{\ctalg{V}} \to \drzf[Ass]{\ctalg{V}} \]
to \cchoch{V}. The cyclic Hochschild cohomology of $V$ is defined as the cohomology of the complex \cchoch{V} and is denoted by \hchoch{V}.
\item[(c)]
Let $m:\clalg{V} \to \clalg{V}$ be a \ci-structure. The cyclic Harrison complex of the \ci-algebra $V$ is defined on the module consisting of all 0-forms:
\[ \ccaq{V}:=\Sigma\drzf[Lie]{\clalg{V}}. \]
The differential on this complex is the (suspension of the) Lie operator of the vector field $m$;
\[ L_m:\drzf[Lie]{\clalg{V}} \to \drzf[Lie]{\clalg{V}}. \]
The cyclic Harrison cohomology of $V$ is defined as the cohomology of the complex \ccaq{V} and is denoted by \hcaq{V}.
\end{enumerate}
\end{defi}

\begin{rem}
Note that a consequence of the above definition is that the cyclic Chevalley-Eilenberg complex is just the Chevalley-Eilenberg complex with \emph{trivial coefficients}.
\end{rem}

It will be useful for the purposes of the next section to describe equivalent formulations for the cyclic cohomology of an \ai-algebra. For this we will need the following lemmas:

\begin{lemma} \label{lem_difcom}
Let $V$ be an \ai-algebra and consider the maps $b'$ and $b$ defined by diagrams \eqref{fig_barhom} and \eqref{fig_chomid} respectively:
\begin{enumerate}
\item[(i)]
\[ bN=-Nb'. \]
\item[(ii)]
\[ b'(1-z)=-(1-z)b. \]
\end{enumerate}
\end{lemma}

\begin{proof}
\
\begin{enumerate}
\item[(i)]
This is a tautological consequence of Definition \ref{def_barhom}, Lemma \ref{lem_nrmdif} and Lemma \ref{lem_schf} part (v).
\item[(ii)]
Let $x \in \Sigma V^*$ and $y \in (\Sigma V^*)^{\cotimes n}$ for $n \geq 0$:
\begin{displaymath}
\begin{split}
b'(1-z)xy & = b'([x,y]) = [m(x),y] + (-1)^{|x|}[x,m(y)]. \\
(1-z)b(xy) & =(1-z)\thiso[Ass]^{-1}L_m(dx\cdot y), \\
& = -(1-z)\thiso[Ass]^{-1}(dm(x)\cdot y) - (-1)^{|x|}(1-z)\thiso[Ass]^{-1}(dx\cdot m(y)). \\
\end{split}
\end{displaymath}
By a calculation similar to that of equation \eqref{eqn_nrmdif} we could deduce that
\[ \thiso[Ass]^{-1}(d(x_1\cotimes\ldots\cotimes x_i)\cdot y) = (1+z+\ldots+z^{i-1}) \cdot x_1\cotimes\ldots\cotimes x_i \cdot y \]
and since $(1-z)(1+z+\ldots+z^{i-1})=1-z^i$ we use equation \eqref{eqn_cyccom} to conclude that
\[ (1-z)b(xy) = -\left([m(x),y] + (-1)^{|x|}[x,m(y)]\right), \]
therefore $b'(1-z)=-(1-z)b$ as claimed.
\end{enumerate}
\end{proof}

We can now define a new bicomplex computing the cyclic cohomology of an \ai-algebra:

\begin{defi} \label{def_tsygan}
Let $V$ be an \ai-algebra and consider the maps $b'$ and $b$ defined by diagrams \eqref{fig_barhom} and \eqref{fig_chomid} respectively. The Tsygan bicomplex \ctsygan{V} is the bicomplex;
\[ \xymatrix{\left(\Sigma\left[\col{V}\right],b\right) \ar^{1-z}[r] & \left(\Sigma\left[\col{V}\right],b'\right) \ar^{N}[r] & \left(\Sigma\left[\col{V}\right],b\right) \ar^-{1-z}[r] & \ldots}. \]
An element has homogeneous bidegree $(i,j)$ if it is in the $i$th column from the left (where the leftmost column has bidegree $(0,\bullet)$) and has degree $j$ in the graded profinite module $\Sigma\left[\col{V}\right]$. The Tsygan cohomology of the \ai-algebra $V$ is defined as the total cohomology of the bicomplex \ctsygan{V} \emph{formed by taking direct sums}, that is to say that it is the cohomology of the complex
\[ C:= \prod_{n \in \mathbb{Z}}\left(\bigoplus_{i+j=n}C^{ij}\right) \]
where $C^{ij}$ is the component of \ctsygan{V} of bidegree $(i,j)$.
\end{defi}

\begin{rem}
Note that the odd numbered columns are copies of the Bar complex \cbr{V} whilst the even numbered columns are isomorphic to the Hochschild complex \choch{V}{V^*}.
\end{rem}

\begin{rem}
Since the underlying module $V$ is $\mathbb{Z}$-graded, the direct product totalisation of the Tsygan complex might a priori lead to different cohomology. This issue does not arise if we assume that the underlying module $V$ is connected ($\mathbb{Z}_{\geq 0}$-graded) since in this case, the direct product and the direct sum totalisation are in agreement.
\end{rem}

Let us now show that Tsygan cohomology computes the same cohomology as cyclic Hochschild cohomology which we defined earlier:

\begin{lemma} \label{lem_cyciso}
Let $V$ be an \ai-algebra:
\begin{enumerate}
\item[(i)]
The map $b':\col{V} \to \col{V}$ defined by diagram \eqref{fig_barhom} lifts uniquely to a map
\[ b' : \prod_{i=1}^\infty \left((\Sigma V^*)^{\cotimes i}\right)_{Z_i} \to \prod_{i=1}^\infty \left((\Sigma V^*)^{\cotimes i}\right)_{Z_i}. \]
This gives rise to the following identity:
\[ \cchoch{V} = \left(\Sigma\left[\prod_{i=1}^\infty \left((\Sigma V^*)^{\cotimes i}\right)_{Z_i}\right],b'\right). \]
\item[(ii)]
There is a quasi-isomorphism of complexes
\[ \left(\Sigma\left[\prod_{i=1}^\infty \left((\Sigma V^*)^{\cotimes i}\right)_{Z_i}\right],b'\right) \overset{\sim}{\to} \ctsygan{V} \]
which is given by mapping $\Sigma\left[\prod_{i=1}^\infty \left((\Sigma V^*)^{\cotimes i}\right)_{Z_i}\right]$ onto the leftmost column of \ctsygan{V} by
\[ \Sigma\left[\prod_{i=1}^\infty \left((\Sigma V^*)^{\cotimes i}\right)_{Z_i}\right] \overset{N}{\longrightarrow} \Sigma\left[\col{V}\right] = CC^{0,\bullet}_\mathrm{Tsygan}(V) \subset \ctsygan{V}. \]
\end{enumerate}
\end{lemma}

\begin{proof}
\
\begin{enumerate}
\item[(i)]
Recall that
\[ \drzf[Ass]{\ctalg{V}} = \ctalg{V}/[\ctalg{V},\ctalg{V}]. \]
The claim then follows tautologically from equation \eqref{eqn_cyccom} and the definitions.
\item[(ii)]
This map is a quasi-isomorphism because
\[ \xymatrix{ \left(\Sigma\left[\prod_{i=1}^\infty \left((\Sigma V^*)^{\cotimes i}\right)_{Z_i}\right],b'\right) \ar^-{N}[r] & \ctsygan{V}} \]
is a resolution of $\left(\Sigma\left[\prod_{i=1}^\infty \left((\Sigma V^*)^{\cotimes i}\right)_{Z_i}\right],b'\right)$.
\end{enumerate}
\end{proof}

Following Lemma \ref{lem_cyciso} we will denote the cohomology of the Tsygan complex by \hchoch{V}. We will now describe the periodicity long exact sequence linking Hochschild cohomology with cyclic cohomology:

\begin{prop} \label{prop_sbiseq}
Let $V$ be a unital \ai-algebra, then we have the following long exact sequence in cohomology:
\[ \xymatrix{ \ldots \ar^-{B}[r] & \hchoch[n-2]{V} \ar^{S}[r] & \hchoch[n]{V} \ar^{I}[r] & \hhoch[n]{V}{V^*} \ar^{B}[r] & \hchoch[n-1]{V} \ar^-{S}[r] & \ldots \\ }. \]
\end{prop}

\begin{proof}
There is a short exact sequence of complexes
\[ \xymatrix{ 0 \ar[r] & \Sigma^{-2}\ctsygan{V} \ar^{S}[r] & \ctsygan{V} \ar^{I}[r] & CC\{2\}^{\bullet\bullet}_{\mathrm{Tsygan}}(V) \ar[r] & 0 \\ } \]
where $CC\{2\}^{\bullet\bullet}_{\mathrm{Tsygan}}(V)$ is the total complex of the first two columns of \ctsygan{V};
\[ CC\{2\}^{\bullet\bullet}_{\mathrm{Tsygan}}(V) := \xymatrix{\left(\Sigma\left[\col{V}\right],b\right) \ar^{1-z}[r] & \left(\Sigma\left[\col{V}\right],b'\right)}. \]

By Lemma \ref{lem_unitch} the second column of $CC\{2\}^{\bullet\bullet}_{\mathrm{Tsygan}}(V)$ is acyclic, so $CC\{2\}^{\bullet\bullet}_{\mathrm{Tsygan}}(V)$ is quasi-isomorphic to its first column which is just the complex \choch{V}{V^*} by Remark \ref{rem_chomid}. The long exact sequence in cohomology of the proposition is derived from the above short exact sequence of complexes.
\end{proof}

\begin{rem} \label{rem_sbiseq}
We can describe the connecting map
\[ B: \hhoch[n]{V}{V^*} \to \hchoch[n-1]{V} \]
more explicitly as follows: first of all as mentioned in Proposition \ref{prop_sbiseq}, the projection
\[ \pi: CC\{2\}^{\bullet\bullet}_{\mathrm{Tsygan}}(V) \to \left(\Sigma\left[\col{V}\right],b\right) \]
of $CC\{2\}^{\bullet\bullet}_{\mathrm{Tsygan}}(V)$ onto the first column is a quasi-isomorphism. A simple check reveals that this map has a section (and hence a quasi-inverse),
\[ i:\left(\Sigma\left[\col{V}\right],b\right) \to CC\{2\}^{\bullet\bullet}_{\mathrm{Tsygan}}(V) \]
which is defined as follows;
\begin{displaymath}
\begin{array}{ccc}
\col{V} & \to & \col{V} \oplus \col{V}, \\
x & \mapsto & x \oplus -h(1-z)[x]; \\
\end{array}
\end{displaymath}
where $h:\col{V} \to \col{V}$ is the contracting homotopy of Lemma \ref{lem_unitch}. It follows from the definitions and a simple check that the connecting map $B$ is the map induced by the following map of complexes:
\begin{equation} \label{eqn_sbiseq}
\begin{array}{ccc}
\left(\Sigma\left[\col{V}\right],b\right) & \to & \left(\Sigma\left[\prod_{i=1}^\infty\left((\Sigma V^*)^{\cotimes i}\right)_{Z_i}\right],b'\right), \\
x & \mapsto & -h(1-z)[x]. \\
\end{array}
\end{equation}

A similar check also shows that the map $I:\hchoch[n]{V} \to \hhoch[n]{V}{V^*}$ is the map induced by the following map of complexes:
\begin{displaymath}
\begin{array}{ccc}
\left(\Sigma\left[\prod_{i=1}^\infty\left((\Sigma V^*)^{\cotimes i}\right)_{Z_i}\right],b'\right) & \to & \left(\Sigma\left[\col{V}\right],b\right), \\
x & \mapsto & N\cdot x. \\
\end{array}
\end{displaymath}
\end{rem}

\subsection{Connes complex and normalised cyclic cohomology}
In the case of a unital $\ai$-algebra (just as in the case of a unital strictly associative algebra) the Tsygan complex could be simplified, excluding the acyclic columns.  The obtained bicomplex is called the Connes complex. It also makes sense to consider the normalised cyclic complex. In contrast to the Hochschild cohomology, the normalised  cyclic cohomology do not agree with cyclic cohomology, the discrepancy being the cyclic cohomology of the ground field. In the strictly associative case this is proved in  \cite{loday}.

\begin{defi} \label{def_connes}
Let $V$ be a unital $\ai$-algebra and let $\tau,\boldsymbol{t}$ be a topological basis of $\Sigma V^*$, where $\tau$ is dual to the unit $1\in V$ so that $\Sigma V^*$ splits (as a module)
\begin{equation} \label{eqn_unisplit}
\Sigma V^*:=\Sigma\gf^*\oplus \Sigma\bar{V}^*,
\end{equation}
where $\Sigma\bar{V}^*$ is the module linearly generated by the $t_i$ and $\gf^*$ is generated by $\tau$.
\begin{enumerate}
\item[(i)]
Define the map
\[ B':\col{V}\to\col{V} \]
of degree $-1$ by the formula $B':=Nh(1-z)$, where $h$ is the contracting homotopy of Lemma \ref{lem_unitch}. It follows from Lemmas \ref{lem_unitch} and \ref{lem_difcom} that the following diagram is a bicomplex.
\[ \xymatrix{\left(\Sigma\left[\col{V}\right],b\right) \ar^{\Sigma^{-1} B'}[r] & \left(\left[\col{V}\right],b\right) \ar^{\Sigma B'}[r] & \left(\Sigma^{-1}\left[\col{V}\right],b\right) \ar^-{\Sigma^{-1} B'}[r] & \ldots}. \]
An element has homogeneous bidegree $(i,j)$ if it is in the $i$th column from the left (where the leftmost column has bidegree $(0,\bullet)$) and has degree $j$ in the graded profinite module $\Sigma^{1-i}\left[\col[r]{V}\right]$. We call the direct sum totalisation of this bicomplex the Connes complex and denote it by $\connes{V}$.
\item[(ii)]
Define the map
\[ \bar{B}':\col{V}\to\col{V} \]
of degree $-1$ by the formula $\bar{B}':=Nh$. The subcomplex
\[ \xymatrix{\left(\Sigma\left[\prod_{i=0}^\infty \Sigma V^*\cotimes(\Sigma\bar{V}^*)^{\cotimes i}\right],b\right) \ar^{\Sigma^{-1}\bar{B}'}[r] & \left(\left[\prod_{i=0}^\infty \Sigma V^*\cotimes(\Sigma\bar{V}^*)^{\cotimes i}\right],b\right)  \ar^-{\Sigma\bar{B}'}[r] & \ldots}. \]
of the Connes complex is called the normalised Connes compex and denoted by $\nconnes{V}$.
\end{enumerate}
\end{defi}

\begin{lemma}
Let $V$ be a non-negatively graded unital $\ai$-algebra.
\begin{enumerate}
\item[(i)]
The map from $\connes{V}$ to $\ctsygan{V}$ given by mapping the $i$th column of $\connes{V}$ onto the $(2i)$th and $(2i+1)$th column of $\ctsygan{V}$ according to the formula
\begin{equation} \label{eqn_colctr}
\begin{array}{ccc}
\Sigma^{1-i}\left[\col[r]{V}\right] & \to & \Sigma\left[\col[r]{V}\right] \oplus \Sigma\left[\col[r]{V}\right], \\
x & \mapsto & (-1)^i x \oplus (-1)^{i+1} h(1-z)[x]; \\
\end{array}
\end{equation}
is a quasi-isomorphism.
\item[(ii)]
The map from $\cchoch{V}$ to $\connes{V}$ which is given by mapping $\Sigma\left[\prod_{i=1}^\infty \left((\Sigma V^*)^{\cotimes i}\right)_{Z_i}\right]$ onto the leftmost column of $\connes{V}$ by
\[ \Sigma\left[\prod_{i=1}^\infty \left((\Sigma V^*)^{\cotimes i}\right)_{Z_i}\right] \overset{N}{\longrightarrow} \Sigma\left[\col{V}\right] = CC^{0,\bullet}_\mathrm{Connes}(V) \subset \connes{V} \]
is a quasi-isomorphism.
\item[(iii)]
The canonical inclusion of the normalised Connes complex $\nconnes{V}$ into the Connes complex $\connes{V}$ is a quasi-isomorphism.
\end{enumerate}
\end{lemma}

\begin{proof}
It follows from Lemmas \ref{lem_unitch} and \ref{lem_difcom} (ii) that the map described by \eqref{eqn_colctr} induces a map on the totalisations of the bicomplexes respecting the differentials. Consider the filtration on $\connes{V}$ given by columns and the filtration on $\ctsygan{V}$ given by filtering by double columns. The map \eqref{eqn_colctr} respects this filtration and it follows from Remark \ref{rem_sbiseq} that it induces an isomorphism between the spectral sequences at the $E^1$ term. Therefore it is a quasi-isomorphism. Part (ii) is an immediate consequence of part (i) and Lemma \ref{lem_cyciso} (ii). Part (iii) follows from a similar spectral sequence argument by applying Proposition \ref{prop_htpret}.
\end{proof}

One can also consider the appropriate normalised version of cyclic cohomology for unital $\ai$-algebras. Unlike the noncyclic case, its cohomology is not the same as that of the unnormalised complex, but the discrepancy is measured by the cyclic cohomology of the field (cf. \cite{loday}).

\begin{prop} \label{lem_nrmles}
Let $V$ be a minimal unital non-negatively graded $\ai$-algebra and let $\tau,\boldsymbol{t}$ be a topological basis of $\Sigma V^*$ so that $\Sigma V^*$ splits (as a module)
\[ \Sigma V^*:=\Sigma\gf^*\oplus \Sigma\bar{V}^* \]
in the same way as before (cf. \eqref{eqn_unisplit}). The subspace of $\cchoch{V}$ of normalised noncommutative 0-forms defined by the formula:
\[ \ccnhoch{V}:=\Sigma\left[\{q \in \drzf[Ass]{\widehat{T}[\Sigma\bar{V}^*]}:\ q \text{ vanishes at zero} \}\right] \]
forms a subcomplex of $\cchoch{V}$. It fits into a long exact sequence in cohomology:
\begin{equation} \label{eqn_nrmles}
\xymatrix{ \ldots \ar[r]^-{\partial} & \hcnhoch[n]{V} \ar[r]^{\iota} & \hchoch[n]{V} \ar[r]^{\pi} & \hchoch[n]{\gf} \ar[r]^{\partial} & \hcnhoch[n+1]{V} \ar[r]^-{i} & \ldots}
\end{equation}
where $\pi$ is dual to the inclusion $\gf\subset V$ and $\iota$ is the canonical inclusion of the normalised cochains into $\cchoch{V}$. Moreover, if $V$ is a minimal unital non-negatively graded $\ci$-algebra, then the connecting map $\partial$ is equal to zero.
\end{prop}

\begin{proof}
By definition, any unital $\ai$-structure has the form \eqref{eqn_unitform}; therefore it sends a normalised 0-form to another normalised 0-form plus a term given by the derivation $\ad\tau$. This term is a commutator and therefore vanishes. It follows that $\ccnhoch{V}$ is a subcomplex as claimed.

Consider the following commutative diagram:
\begin{equation} \label{eqn_normseq}
\xymatrix{ \ccnhoch[n]{V} \ar[dd]_{N} \ar[r]^{\iota} & \cchoch[n]{V} \ar[r]^{\pi} \ar[d]_{N} & \cchoch[n]{\gf} \ar[d]_{N} \\ & \connes[n]{V} \ar[r]^{\pi} & \connes[n]{\gf} \\ \ker(\pi) \ar@{^{(}->}[r] & \nconnes[n]{V} \ar@{^{(}->}[u] \ar[r]^{\pi} & \nconnes[n]{\gf} \ar@{^{(}->}[u] }
\end{equation}
Since the central and rightmost vertical arrows are quasi-isomorphisms and the bottom row is a short exact sequence, we will have the claimed long exact sequence in cohomology if we can prove that the lefthand vertical arrow is a quasi-isomorphism.

To this end, consider the filtration $\{F_p\}_{p=1}^\infty$ of $K:=\ker(\pi)$ whose $p$th term is the subcomplex
\[\xymatrix{\left(\Sigma\left[{\begin{subarray}{c} \prod_{r=p}^\infty \Sigma V^*\cotimes(\Sigma\bar{V}^*)^{\cotimes r} \\ \oplus \\ (\Sigma\bar{V}^*)^{\cotimes p} \end{subarray}}\right],b\right) \ar^{\Sigma^{-1}\bar{B}'}[r] & \left(\left[{\begin{subarray}{c} \prod_{r=p}^\infty \Sigma V^*\cotimes(\Sigma\bar{V}^*)^{\cotimes r} \\ \oplus \\ (\Sigma\bar{V}^*)^{\cotimes p} \end{subarray}} \right],b\right) \ar^-{\Sigma^{-1}\bar{B}'}[r] & \ldots}\]
and the filtration of $\ccnhoch{V}$ whose $p$th term is the subcomplex
\[ \left(\Sigma\left[\prod_{r=p}^\infty \left((\Sigma \bar{V}^*)^{\cotimes r}\right)_{Z_r}\right],b'\right). \]
Both filtrations are exhaustive and bounded below. The lefthand vertical arrow of \eqref{eqn_normseq} respects these filtrations.

Let $v\in\Sigma\bar{V}^*$ and $q\in\Sigma\bar{V}^{*\cotimes p-1}$. It follows from the definition of a unital $\ai$-structure \eqref{eqn_unitform} and a simple calculation that
\[L_m(dv\cdot q)= d\tau\cdot[v,q] + dv\cdot(\text{terms of order $\geq p$}).\]
From this calculation we can deduce the shape of the complex $F_p K/F_{p+1} K$:
\[\xymatrix{\Sigma\bar{V}^{*\otimes p} \ar[rd]_{N} & \Sigma\bar{V}^{*\otimes p} \ar[rd]_{N} & \Sigma\bar{V}^{*\otimes p} \ar[rd]_{N} \\ \Sigma\bar{V}^{*\otimes p} \ar[u]^{1-z_p} & \Sigma\bar{V}^{*\otimes p} \ar[u]^{1-z_p} & \Sigma\bar{V}^{*\otimes p} \ar[u]^{1-z_p} & \ldots } \]
It follows that the lefthand vertical arrow of \eqref{eqn_normseq} induces an isomorphism between the spectral sequences associated to these filtrations at the $E^1$ term and hence it is a quasi-isomorphism as claimed.

Lastly, we must verify that if $V$ is a $\ci$-algebra then the connecting map $\partial$ is equal to zero. Since $V$ is a (minimal) $\ci$-algebra with $\ci$-structure $m$ it follows that given any $v\in\Sigma V^*$, $m(v)$ is a sum of commutators and therefore
\[ b(v) = \thiso[Ass]^{-1}L_m d(v) = \thiso[Ass]^{-1} d(L_m(v)) = 0. \]
Using this fact one can check that the connecting map is indeed equal to zero.
\end{proof}

\section{The Hodge Decomposition of Hochschild and Bar Cohomology} \label{sec_hdgdec}

In this section we will be concerned with constructing the Hodge decomposition of the Hochschild and bar cohomology of a \ci-algebra. Given a \ci-algebra $V$ we will construct both the Hodge decomposition of the Hochschild cohomology of $V$ with coefficients in $V^*$ as well as the Hodge decomposition of the bar cohomology of $V$.

The Hodge decomposition of the Hochschild (co)homology of a commutative algebra has been described by many authors such as \cite{gersch}, \cite[\S4.5,\S4.6]{loday} and \cite{natsch}. Our approach will be to determine a Hodge decomposition for \ci-algebras and this will naturally include a Hodge decomposition of the Hochschild cohomology of a commutative algebra. Our results however are more than just a mere generalisation of the results contained in \cite{loday} and \cite{natsch}. Considering the Hodge decomposition in the broader perspective of \ci-algebras leads us to use the framework of noncommutative geometry. This is a very natural setting in which to construct the Hodge decomposition and leads to a much more conceptual approach, resulting in a streamlining of the calculations. The pay off is that we are able to obtain new results, even for the Hodge decomposition of commutative algebras.

Let $W$ be a profinite graded module. For the rest of this section we will denote the canonical associative multiplication on $\widehat{T}W$ by $\mu$. The cocommutative comultiplication (\emph{shuffle product}) on $\widehat{T}W$ defined by Remark \ref{rem_ccmult} will be denoted by $\Delta$. This gives $\widehat{T}W$ the structure of a Hopf algebra.

The Hodge decomposition will be constructed from a spectral decomposition of an operator which we will refer to as the \emph{modified shuffle operator} and which we will now define. The original `shuffle operator' was used in \cite{gersch} and \cite{natsch} and seems to have originally been introduced by Barr in \cite{barr}.

\begin{defi} \label{def_shufop}
Let $W$ be a profinite graded module. The \emph{modified shuffle operator} $s$ is defined as
\[ s:= \mu\Delta:\widehat{T}W \to \widehat{T}W. \]
We denote its components by
\[ s_n:W^{\cotimes n} \to W^{\cotimes n}. \]

We also define a second operator $\tilde{s}:\prod_{i=1}^\infty W^{\cotimes i} \to \prod_{i=1}^\infty W^{\cotimes i}$ by the commutative diagram;
\begin{displaymath}
\xymatrix{ W\cotimes\widehat{T}W \ar^{1\cotimes s}[r] \ar@{=}[d] & W\cotimes\widehat{T}W \ar@{=}[d] \\ \prod_{i=1}^\infty W^{\cotimes i} \ar^{\tilde{s}}[r] & \prod_{i=1}^\infty W^{\cotimes i}}
\end{displaymath}
We denote its components by
\[ \tilde{s}_n:W^{\cotimes n+1} \to W^{\cotimes n+1}. \]
\end{defi}

Next we will define a filtration of $\widehat{T}W$ which comes from the observation $\widehat{T}W = \widehat{\mathcal{U}}(\widehat{L}W)$. This will assist us in our calculations later:

\begin{defi} \label{def_pbwfil}
Let $W$ be a profinite graded module. The increasing filtration $\{F_p(\widehat{T}W)\}_{p=0}^\infty$ of $\widehat{T}W$ is defined as follows: $F_p(\widehat{T}W)$ is defined as the submodule of $\widehat{T}W$ which consists of all (possibly uncountably infinite) linear combinations of elements of the form
\[ g_1 \ldots g_i, \]
where $g_1,\ldots,g_i \in \widehat{L}W$ and $0 \leq i \leq p$. By convention $F_0(\widehat{T}W)= \gf$.
\end{defi}

\begin{rem}
Note that this filtration is not exhaustive but that $W^{\cotimes n} \subset F_n(\widehat{T}W)$. Also note that this filtration is preserved by the modified shuffle operator $s$. This follows from the fact that $\widehat{L}W$ consists of primitive elements of the Hopf algebra $(\widehat{T}W,\mu,\Delta)$.
\end{rem}

Now we will determine the spectral decomposition of the modified shuffle operator $s$. We will do this by exhibiting a certain polynomial that annihilates it. Let us define the polynomial $\nu_n(X) \in \mathbb{Z}(X)$  for $n \geq 0$ by the formula
\begin{equation} \label{eqn_poldef}
\nu_n(X):=\prod_{i=0}^n (X-\lambda_i), \quad \lambda_i:=2^i.
\end{equation}

\begin{lemma} \label{lem_shuann}
Let $W$ be a profinite graded module. For all $n \geq 0$;
\[ \nu_n(s_n)=\prod_{i=0}^n (s_n-\lambda_i\id)=0. \]
\end{lemma}

\begin{proof}
Let us prove the following equation:
\begin{equation} \label{eqn_shufev}
s_n(x)=\lambda_p x \mod F_{p-1}(\widehat{T}W), \quad \text{for all } x \in F_p(\widehat{T}W)\cap W^{\cotimes n}.
\end{equation}
Since the modified shuffle operator $s$ preserves the filtration, we may assume that $x$ is a linear combination of elements of the form $g_1\ldots g_p$ where $g_i \in \widehat{L}W$.

We calculate that for all $g_1,\ldots,g_p \in \widehat{L}W$;
\begin{displaymath}
\begin{split}
s_n(g_1\ldots g_p) & = \mu\Delta(g_1\ldots g_p), \\
& =\mu\left[(g_1\cotimes 1 + 1 \cotimes g_1)\ldots(g_p\cotimes 1 + 1 \cotimes g_p)\right], \\
& =2^p g_1\ldots g_p \mod F_{p-1}(\widehat{T}W). \\
\end{split}
\end{displaymath}
The last equality follows since by using commutators we can transform
\[ \mu\left[(g_1\cotimes 1 + 1 \cotimes g_1)\ldots(g_p\cotimes 1 + 1 \cotimes g_p)\right] \]
into $2^p g_1\ldots g_p$. By definition, these commutating elements are in $F_{p-1}(\widehat{T}W)$. This calculation implies equation \eqref{eqn_shufev} from which the lemma follows as a trivial consequence.
\end{proof}

\begin{rem}
Note $\nu_n(X)$ is not the minimal polynomial for $s_n$. The minimal polynomial of $s_n$ was constructed in \cite{gersch} and is given by the formula:
\begin{displaymath}
p_n(X):=\left\{
\begin{array}{ccc}
\prod_{i=1}^n (X-\lambda_i) & , & n \geq 1 \\
X-\lambda_0 & , & n = 0 \\
\end{array}\right..
\end{displaymath}
\end{rem}

\begin{defi} \label{def_idemdf}
Let $W$ be a profinite graded module. We define a family of operators
\begin{displaymath}
\begin{array}{ll}
e(j):\widehat{T}W \to \widehat{T}W, & i \geq 0; \\
e_n(j):W^{\cotimes n} \to W^{\cotimes n}; \\
\end{array}
\end{displaymath}
as the Lagrange interpolation polynomials of the operator $s_n$:
\begin{displaymath}
e_n(j):=\left\{
\begin{array}{ccc}
\left[ \underset{\begin{subarray}{c} 0 \leq r \leq n \\ r \neq j \end{subarray}}{\prod} (\lambda_j - \lambda_r) \right]^{-1} \underset{\begin{subarray}{c} 0 \leq r \leq n \\ r \neq j \end{subarray}}{\prod} (s_n - \lambda_r\id) & , & 0 \leq j \leq n; \\
0 & , & j > n. \\
\end{array}
\right.
\end{displaymath}

We also define a family of operators
\[ \tilde{e}(j):\prod_{i=1}^\infty W^{\cotimes i} \to \prod_{i=1}^\infty W^{\cotimes i}, \quad j \geq 0 \]
by the commutative diagram
\begin{displaymath}
\xymatrix{W\cotimes\widehat{T}W \ar^{1\cotimes e(j)}[r] \ar@{=}[d] & W\cotimes\widehat{T}W \ar@{=}[d] \\ \prod_{i=1}^\infty W^{\cotimes i} \ar^{\tilde{e}(j)}[r] & \prod_{i=1}^\infty W^{\cotimes i}}
\end{displaymath}
and denote their components by
\[ \tilde{e}_n(j):W^{\cotimes n+1} \to W^{\cotimes n+1}. \]
\end{defi}

\begin{lemma} \label{lem_shidem}
Let $W$ be a profinite graded module. We have the following identities:
\begin{displaymath}
\begin{array}{rl}
\textnormal{(a)} & s_n=\lambda_0 e_n(0) + \ldots + \lambda_n e_n(n). \\
\textnormal{(a*)} & s=\sum_{j=0}^\infty \lambda_j e(j). \\
\textnormal{(b)} & \id_n=e_n(0) + \ldots + e_n(n). \\
\textnormal{(b*)} & \id = \sum_{j=0}^\infty e(j). \\
\textnormal{(c)} & e_n(i) \circ e_n(j)=\left\{
\begin{array}{ccc}
e_n(i) & , & i=j; \\
0 & , & i \neq j. \\
\end{array}
\right.\\
\textnormal{(c*)} & e(i) \circ e(j)=\left\{
\begin{array}{ccc}
e(i) & , & i=j; \\
0 & , & i \neq j. \\
\end{array}
\right. \\
\end{array}
\end{displaymath}
\end{lemma}

\begin{proof}
(a), (b) and (c) are a formal consequence of Lemma \ref{lem_shuann}. Since we assume that our ground ring {\gf} contains the field $\mathbb{Q}$ and since the polynomial $\nu_n(X) \in \mathbb{Z}(X)$ defined by equation \eqref{eqn_poldef} annihilates $s_n$ and has no repeated roots, elementary linear algebra implies that $s_n$ is diagonalisable. The map $e(j)$ is the projection onto the eigenspace
\[ \{x \in \widehat{T}W : s(x)=\lambda_j x \}. \]
(a*), (b*) and (c*) are trivial consequences of (a), (b) and (c) respectively.
\end{proof}

\begin{rem}
Since $\tilde{e}(j):=1\cotimes e(j)$, the same identities hold when we replace $e(j)$ with $\tilde{e}(j)$, $s$ with $\tilde{s}$ and $\id_n$ with $\id_{n+1}$ in the above Lemma.
\end{rem}

We want to use the spectral decomposition of the modified shuffle operator $s$ to construct a decomposition of the relevant cohomology theories. For this we need the following lemma:

\begin{lemma} \label{lem_shucom}
Let $V$ be a \ci-algebra and consider the maps $b'$ and $b$ defined by diagrams \eqref{fig_barhom} and \eqref{fig_chomid} respectively. Also consider the maps $s,\tilde{s}:\col{V} \to \col{V}$ defined in Definition \ref{def_shufop}. We have the following identities:
\begin{enumerate}
\item[(i)]
\[ s \circ b'=b' \circ s. \]
\item[(ii)]
\[ \tilde{s} \circ b=b \circ \tilde{s}. \]
\end{enumerate}
\end{lemma}

\begin{proof}
Let $m:\clalg{V} \to \clalg{V}$ be the \ci-structure on $V$. Recall from Remark \ref{rem_calseq} that this could be considered to be an \ai-structure $m:\ctalg{V} \to \ctalg{V}$ which is also a Hopf algebra derivation of the Hopf algebra $(\ctalg{V},\mu,\Delta)$:
\begin{enumerate}
\item[(i)]
This is just a trivial consequence of the definition of $b'$ and the fact that $m$ is a derivation of the Hopf algebra $(\ctalg{V},\mu,\Delta)$:
\begin{equation} \label{eqn_shucom}
\mu\Delta m=\mu(m\cotimes 1+1\cotimes m)\Delta=m\mu\Delta.
\end{equation}
Since $b'$ is the restriction of $m$ to \col{V} we obtain $s \circ b'=b' \circ s$
\item[(ii)]
Let $\thiso[Ass]:\col{V} \to \drof[Ass]{\ctalg{V}}$ be the module isomorphism we defined in Lemma \ref{lem_ofisom}. Let $x \in \Sigma V^*$ and $y \in (\Sigma V^*)^{\cotimes n}$ for $n \geq 0$;
\begin{displaymath}
\begin{split}
b\tilde{s}(xy) & = b(xs(y)) = \thiso[Ass]^{-1}L_m(dx\cdot s(y)), \\
& = -\thiso[Ass]^{-1}(dm(x)\cdot s(y))-(-1)^{|x|}x\cdot sm(y); \\
\end{split}
\end{displaymath}
where the computation of the last term follows from equation \eqref{eqn_shucom}. Furthermore;
\begin{displaymath}
\begin{split}
\tilde{s}b(xy) & = \tilde{s}\thiso[Ass]^{-1}L_m(dx\cdot y), \\
& = -\tilde{s}\thiso[Ass]^{-1}(dm(x)\cdot y)-(-1)^{|x|}x\cdot sm(y). \\
\end{split}
\end{displaymath}

We would like to show that given any $u \in \clalg{V}$ and $w \in \ctalg{V}$ we have the following identity:
\begin{equation} \label{eqn_shucomdummy}
\thiso[Ass]^{-1}(du\cdot s(w))=\tilde{s}\thiso[Ass]^{-1}(du\cdot w).
\end{equation}
Since $m$ maps $\Sigma V^*$ to \clalg{V} it will follow from equation \eqref{eqn_shucomdummy} that $\tilde{s} \circ b=b \circ \tilde{s}$. Since \clalg{V} is generated by Lie monomials and \eqref{eqn_shucomdummy} is a tautology for a Lie monomial $u$ of order one (i.e. $u \in \Sigma V^*$), we can proceed by induction and assume there exists $u_1,u_2 \in \clalg{V}$ such that $u=[u_1,u_2]$ and such that \eqref{eqn_shucomdummy} holds for $u_1$ and $u_2$:
\begin{displaymath}
\begin{split}
du\cdot s(w) & = d[u_1,u_2]\cdot s(w) =[du_1,u_2]\cdot s(w) - (-1)^{|u_1||u_2|}[du_2,u_1] \cdot s(w), \\
& = du_1\cdot [u_2,s(w)] - (-1)^{|u_1||u_2|}du_2\cdot [u_1,s(w)]. \\
du \cdot w & = d[u_1,u_2]\cdot w = du_1\cdot [u_2,w] - (-1)^{|u_1||u_2|}du_2\cdot [u_1,w]. \\
\end{split}
\end{displaymath}
By the inductive hypothesis we obtain;
\begin{displaymath}
\begin{split}
\tilde{s}\thiso[Ass]^{-1}(du\cdot w) & = \tilde{s}\thiso[Ass]^{-1}(du_1\cdot [u_2,w]) - (-1)^{|u_1||u_2|}\tilde{s}\thiso[Ass]^{-1}(du_2\cdot [u_1,w]), \\
& = \thiso[Ass]^{-1}(du_1\cdot s([u_2,w])) - (-1)^{|u_1||u_2|}\thiso[Ass]^{-1}(du_2\cdot s([u_1,w])). \\
\end{split}
\end{displaymath}

In order to complete the proof we will need one final auxiliary calculation. Recall that \clalg{V} coincides with the Lie subalgebra of primitive elements of the Hopf algebra $(\ctalg{V},\mu,\Delta)$. Let $a \in \clalg{V}$ and $b \in \ctalg{V}$;
\begin{equation} \label{eqn_shucomdummya}
\begin{split}
s([a,b]) & = \mu[a\cotimes 1+1\cotimes a,\Delta b], \\
& = \mu(a\cotimes 1\cdot\Delta b -(-1)^{|a||b|}\Delta b \cdot 1\cotimes a) + \mu(1\cotimes a\cdot\Delta b - (-1)^{|a||b|}\Delta b \cdot a\cotimes 1), \\
& = a\cdot\mu\Delta b - (-1)^{|a||b|}\mu\Delta b \cdot a - 0, \\
& = [a,s(b)]. \\
\end{split}
\end{equation}
Finally we can establish equation \eqref{eqn_shucomdummy} and hence conclude the proof:
\begin{displaymath}
\begin{split}
\thiso[Ass]^{-1}(du\cdot s(w)) & = \thiso[Ass]^{-1}(du_1\cdot [u_2,s(w)]) - (-1)^{|u_1||u_2|}\thiso[Ass]^{-1}(du_2\cdot [u_1,s(w)]), \\
& = \thiso[Ass]^{-1}(du_1\cdot s([u_2,w])) - (-1)^{|u_1||u_2|}\thiso[Ass]^{-1}(du_2\cdot s([u_1,w])), \\
& = \tilde{s}\thiso[Ass]^{-1}(du\cdot w). \\
\end{split}
\end{displaymath}
\end{enumerate}
\end{proof}

\begin{cor} \label{cor_shucom}
Let $V$ be a \ci-algebra and consider the maps
\[ e(j),\tilde{e}(j):\col{V}\to\col{V} \]
defined by Definition \ref{def_idemdf}. We have the following identities:
\begin{enumerate}
\item[(i)]
For all $j \geq 1$,
\[ e(j) \circ b' = b' \circ e(j). \]
\item[(ii)]
For all $j \geq 0$,
\[ \tilde{e}(j) \circ b = b \circ \tilde{e}(j). \]
\end{enumerate}
\end{cor}

\begin{proof}
This is just a formal consequence of Lemma \ref{lem_shucom} and Lemma \ref{lem_shidem}. $e(j)$ is just the projection onto the eigenspace
\[ \{ x \in \col{V} : s(x)=\lambda_j x \} \]
whilst $\tilde{e}(j)$ is the projection onto the eigenspace
\[ \{ x \in \col{V} : \tilde{s}(x)=\lambda_j x \}. \]
Lemma \ref{lem_shucom} then tells us that $b'$ and $b$ preserve these eigenspaces respectively and hence commute with the projections $e(j)$ and $\tilde{e}(j)$ respectively.
\end{proof}

We are now in a position to state the main theorem of this section; the Hodge decomposition of the Hochschild and bar cohomology of a \ci-algebra:

\begin{theorem}
Let $V$ be a \ci-algebra:
\begin{enumerate}
\item[(i)]
The bar complex \cbr{V} of the \ci-algebra $V$ splits as the direct product of the subcomplexes $\Sigma\left(e(j)\left[\col{V}\right]\right)$:
\[ \cbr{V}=\prod_{j=1}^{\infty}\left(\Sigma\left(e(j)\left[\col{V}\right]\right),b'\right). \]
\item[(ii)]
The Hochschild complex of the \ci-algebra $V$ with coefficients in $V^*$ splits as the direct product of the subcomplexes $\Sigma\left(\tilde{e}(j)\left[\col{V}\right]\right)$:
\[ \choch{V}{V^*}=\prod_{j=0}^\infty\left(\Sigma\left(\tilde{e}(j)\left[\col{V}\right]\right),b\right). \]
\end{enumerate}
\end{theorem}

\begin{proof}
\
\begin{enumerate}
\item[(i)]
This is just a trivial consequence of Lemma \ref{lem_shidem} and Corollary \ref{cor_shucom}. Lemma \ref{lem_shidem} tells us that the module \col{V} splits as a product of submodules;
\[ \col{V}=\prod_{j=1}^\infty e(j)\left[\col{V}\right]. \]
Corollary \ref{cor_shucom} part (i) tells us that when we equip \col{V} with the differential $b'$, the modules $e(j)\left[\col{V}\right]$ are actually subcomplexes of $\left(\col{V},b'\right)$.
\item[(ii)]
Lemma \ref{lem_shidem} tells us that the module \col{V} splits as a product of submodules;
\[ \col{V}=\prod_{j=0}^\infty \tilde{e}(j)\left[\col{V}\right]. \]
Corollary \ref{cor_shucom} part (ii) tells us that when we equip \col{V} with the differential $b$ the modules $\tilde{e}(j)\left[\col{V}\right]$ are actually subcomplexes of $\left(\col{V},b\right)$. By Remark \ref{rem_chomid} the complex \choch{V}{V^*} is isomorphic to $\left(\Sigma\left[\col{V}\right],b\right)$.
\end{enumerate}
\end{proof}

\begin{defi}
Let $V$ be a \ci-algebra. We will define the complex \echoch{V}{V^*}{j} as the subcomplex of \choch{V}{V^*} given by the formula;
\[ \echoch{V}{V^*}{j}:= \left(\Sigma\left(\tilde{e}(j)\left[\col{V}\right]\right),b\right), \quad j \geq 0 \]
and denote its cohomology by \ehhoch{V}{V^*}{j}.
\end{defi}

Let $l: \drof[Lie]{\clalg{V}} \hookrightarrow \drof[Ass]{\ctalg{V}}$ be the map defined by equation \eqref{eqn_liasmp}. Suppose that $V$ is a \ci-algebra, then by Lemma \ref{lem_plcomm} we see that $l$ is a map of complexes;
\[ l:\caq{V}{V^*} \to \choch{V}{V^*}. \]
We shall now show that this map is in fact split:

\begin{prop} \label{prop_aqsplt}
Let $V$ be a \ci-algebra. There is an isomorphism
\[ \caq{V}{V^*} \cong \echoch{V}{V^*}{1}. \]
In particular the map $l:\caq{V}{V^*} \to \choch{V}{V^*}$ is split.
\end{prop}

\begin{proof}
We shall show that
\begin{equation} \label{eqn_idelie}
e(1)\left[\col{V}\right] = \clalg{V}.
\end{equation}
The proposition will then follow as a simple consequence of this. Since the Lie monomials generate the Lie subalgebra of primitive elements of the Hopf algebra $(\ctalg{V},\mu,\Delta)$ we have;
\[ \clalg{V} \subset \{ x \in \col{V}:s(x)=2x \} = e(1)\left[\col{V}\right]. \]

To prove the converse we consider the filtration $\{F_p(\ctalg{V})\}_{p=0}^\infty$ of \ctalg{V} defined by Definition \ref{def_pbwfil}. Let $x \in F_p(\ctalg{V})$ and suppose that $s(x)=2x$, then by equation \eqref{eqn_shufev};
\[ 2x = s(x) = 2^p x \mod F_{p-1}(\ctalg{V}). \]
It follows that if $p>1$ then $x \in F_{p-1}(\ctalg{V})$. By induction on $p$ we conclude that $x \in F_1(\ctalg{V})$. Equation \eqref{eqn_idelie} now follows from the fact that $(\Sigma V^*)^{\cotimes n} \subset F_n(\ctalg{V})$ and from the identity
\[ F_1(\ctalg{V}) = \gf \oplus \clalg{V}. \]

Equation \eqref{eqn_idelie} implies that
\[ \tilde{e}(1)\left[\col{V}\right] = \Sigma V^* \cotimes \clalg{V}. \]
We see from diagram \eqref{fig_clsdtf} that we have the following identity:
\begin{equation} \label{eqn_liespl}
l\left[\drof[Lie]{\clalg{V}}\right] = \thiso[Ass]\tilde{e}(1)\left[\col{V}\right].
\end{equation}
It follows that $\thiso[Ass]^{-1} \circ l$ is an isomorphism between \caq{V}{V^*} and \echoch{V}{V^*}{1} and that $l$ is split.
\end{proof}

Proposition \ref{prop_aqsplt} allows us to establish the analogue of Proposition \ref{prop_htpret} for $\ci$-algebras. Let $V$ be a unital $\ci$-algebra and choose a topological basis $\tau,\boldsymbol{t}$ of $\Sigma V^*$, so that $\Sigma V^*$ splits (as a module)
\[ \Sigma V^*:=\Sigma\gf^*\oplus \Sigma\bar{V}^*,\]
where $\Sigma\bar{V}^*$ is the module linearly generated by the $t_i$ and $\gf^*$ is generated by $\tau$, where $\tau$ is dual to the unit $1 \in V$. Recall that the module of \emph{normalised} noncommutative 1-forms $\drnof[Ass]{\ctalg{V}}$ is the module generated by elements of the form $q\cdot dv$ for $q\in\gf\langle\langle\boldsymbol{t}\rangle\rangle$ and $v\in\Sigma V^*$, cf. equation \eqref{eqn_normform}. Likewise, we can define the module of \emph{normalised} Lie 1-forms $\drnof[Lie]{\clalg{V}}$ as the module which is generated by elements of the form $q\cdot dv$ for $q$ a formal Lie power series in the variables $\boldsymbol{t}=\{t_i\}_{i\in I}$ and $v\in\Sigma V^*$; i.e. $\drnof[Lie]{\clalg{V}}$ is the inverse image of $\drnof[Ass]{\ctalg{V}}$ under the map $l$.

\begin{cor}
The normalised 1-forms $\Sigma\drnof[Lie]{\clalg{V}}$ form a subcomplex of
\[ \caq{V}{V^*}:=\left(\Sigma\drof[Lie]{\clalg{V}},L_m\right) \]
which is quasi-isomorphic to $\caq{V}{V^*}$ under the canonical inclusion.
\end{cor}

\begin{proof}
It is clear that the operator $\tilde{s}$ on $\choch{V}{V^*}$ descends to an operator on the normalised cochains
\[ \drnof[Ass]{\ctalg{V}}\overset{\thiso[Ass]}{\cong}\prod_{i=0}^\infty\Sigma V^*\cotimes(\Sigma\bar{V}^*)^{\cotimes i}. \]
Consequently we obtain a Hodge decomposition
\[\left(\prod_{i=0}^\infty\Sigma V^*\cotimes(\Sigma\bar{V}^*)^{\cotimes i},b\right)=\prod_{j=0}^{\infty}\left(\tilde{e}(j)\left[\prod_{i=0}^\infty\Sigma V^*\cotimes(\Sigma\bar{V}^*)^{\cotimes i}\right],b\right)\]
of the complex of normalised cochains which is compatible with the Hodge decomposition of $\choch{V}{V^*}$ under the inclusion. Using the same arguments as above, we can identify the summand $\tilde{e}(1)\left[\prod_{i=0}^\infty\Sigma V^*\cotimes(\Sigma\bar{V}^*)^{\cotimes i}\right]$ with the image of $\drnof[Lie]{\clalg{V}}$ under the map $l$; hence we have the following commutative diagram:
\begin{displaymath}
\xymatrix{ \tilde{e}(1)\left[\prod_{i=0}^\infty\Sigma V^*\cotimes(\Sigma\bar{V}^*)^{\cotimes i}\right] \ar@{^{(}->}[r] &  \tilde{e}(1)\left[\col{V}\right] \\ \drnof[Lie]{\clalg{V}} \ar@{^{(}->}[r] \ar[u]^{\thiso[Ass]^{-1} \circ l} & \drof[Lie]{\clalg{V}} \ar[u]^{\thiso[Ass]^{-1}\circ l} }
\end{displaymath}
We have just seen that the vertical arrows are isomorphisms and the top horizontal arrow is a quasi-isomorphism by Proposition \ref{prop_htpret}, therefore the bottom arrow is also a quasi-isomorphism.
\end{proof}

\begin{rem}
Given a \ci-algebra $V$ with \ci-structure $m$ it is also possible to construct a Hodge decomposition of the Hochschild cohomology of $V$ with coefficients in $V$. We will describe it briefly here. Choose a topological basis $\{t_i\}_{i\in I}$ of $\Sigma V^*$. One introduces the operator $\bar{s}:\choch{V}{V} \to \choch{V}{V}$ as
\[ \bar{s}(\xi):= \sum_{i \in I} \mu\Delta\xi(t_i)\partial_{t_i}, \quad \xi \in \Der(\ctalg{V}); \]
which is obviously defined independently of the choice of basis. It is then possible to show, using equations \ref{eqn_shucom} and \ref{eqn_shucomdummya}, that this operator commutes with the differential $d:=\ad m$. It follows that the differential $d$ preserves the eigenspaces of the operator $\bar{s}$ and these eigenspaces can be identified using Lemma \ref{lem_shidem} in the usual way to give the Hodge decomposition of \choch{V}{V}. In particular it is possible to show that the Harrison complex \caq{V}{V} splits off of the Hochschild complex \choch{V}{V}.

Note that the Hodge decomposition of \choch{V}{V} was also obtained using different methods
in \cite{markl2}.
\end{rem}

\section{The Hodge Decomposition of Cyclic Hochschild Cohomology}\label{cyclichodge}

In this section we will build on our results from the last section to construct the Hodge decomposition of the cyclic Hochschild cohomology of a \ci-algebra. Recall that in section \ref{sec_iachom} we defined two quasi-isomorphic complexes \cchoch{V} and \ctsygan{V} which compute the cyclic Hochschild cohomology of an \ai-algebra. We will describe a Hodge decomposition for both of these complexes based on the spectral decomposition of the modified shuffle operator $s$ and show that the quasi-isomorphism defined in Lemma \ref{lem_cyciso} respects this decomposition.

In order to apply the work we carried out in section \ref{sec_hdgdec} to the cyclic complexes, we will need to prove the following lemma:

\begin{lemma} \label{lem_nrmshu}
Let $W$ be a profinite graded module and consider the maps
\[ s,\tilde{s}:\prod_{i=1}^\infty W^{\cotimes i} \to \prod_{i=1}^\infty W^{\cotimes i} \]
defined in Definition \ref{def_shufop}. We have the following identities:
\begin{enumerate}
\item[(i)]
\[ 2\tilde{s} \circ N = N \circ s. \]
\item[(ii)]
\[ s\circ(1-z) = (1-z) \circ \tilde{s}. \]
\end{enumerate}
\end{lemma}

\begin{proof}
\
\begin{enumerate}
\item[(i)]
Let $\thiso[Ass]:\prod_{i=1}^\infty W^{\cotimes i} \to \drof[Ass]{\widehat{T}W}$ be the map defined by Lemma \ref{lem_ofisom} and let
\[ \mu:\de[Ass]{\widehat{T}W}\cotimes\de[Ass]{\widehat{T}W} \to \de[Ass]{\widehat{T}W} \]
be the multiplication map $\mu(x \cotimes y):= x \cdot y$. Let $\alpha_i := (-1)^{(|x_1|+\ldots+|x_{i-1}|)(|x_i|+\ldots+|x_n|)}$ and let $x:=x_1\cotimes\ldots\cotimes x_n \in W^{\cotimes n}$;
\begin{displaymath}
\begin{split}
\tilde{s} N(x) & = \tilde{s}\left(\sum_{i=1}^{n} \alpha_i x_i\cotimes\ldots\cotimes x_n\cotimes x_1 \cotimes \ldots \cotimes x_{i-1}\right), \\
& = \sum_{i=1}^{n} \alpha_i x_i\cotimes \mu\Delta(x_{i+1}\cotimes\ldots\cotimes x_n\cotimes x_1 \cotimes \ldots \cotimes x_{i-1}). \\
N s(x) & = \thiso[Ass]^{-1}d\mu\Delta(x), \\
& = \thiso[Ass]^{-1}d\mu(\Delta(x_1)\ldots\Delta(x_n)), \\
& = \thiso[Ass]^{-1}\mu(d\cotimes 1+1\cotimes d)[(x_1\cotimes 1+1\cotimes x_1)\ldots(x_n\cotimes 1+1\cotimes x_n)]. \\
\end{split}
\end{displaymath}
Now $d\cotimes 1$ and $1\cotimes d$ are derivations on the algebra $\de[Ass]{\widehat{T}W}\cotimes\de[Ass]{\widehat{T}W}$, therefore
\begin{equation} \label{eqn_nrmshudummy}
Ns(x)=  \thiso[Ass]^{-1}\mu\left(\sum_{i=1}^n (-1)^{|x_1|+\ldots+|x_{i-1}|}(x_1\cotimes 1+1\cotimes x_1)\ldots (dx_i\cotimes 1+1\cotimes dx_i)\ldots (x_n\cotimes 1+1\cotimes x_n)\right).
\end{equation}

We will need the following auxiliary calculation in order to complete the proof: Let $u,v,w \in \de[Ass]{\widehat{T}W}$;
\begin{displaymath}
\begin{split}
\mu[(u\cotimes 1+1\cotimes u)\cdot(v\cotimes w)] & = \mu(uv\cotimes w + (-1)^{|u||v|}v\cotimes uw), \\
& = uvw + (-1)^{|u||v|}vuw. \\
\mu[(v\cotimes w)\cdot(u\cotimes 1+1\cotimes u)] & = \mu((-1)^{|u||w|}vu\cotimes w + v\cotimes wu), \\
& = (-1)^{|u||w|}vuw + vwu. \\
\end{split}
\end{displaymath}
We see that
\[ \mu\left([u\cotimes 1 + 1\cotimes u , v\cotimes w]\right) = [u,\mu(v \cotimes w)] = 0 \mod [\de[Ass]{\widehat{T}W},\de[Ass]{\widehat{T}W}]. \]

Let $\beta_i:=(-1)^{(|x_i|+1)(|x_1|+\ldots+|x_{i-1}|+|x_{i+1}|+\ldots+|x_n|)}$. Applying the preceding calculation to equation \eqref{eqn_nrmshudummy} yields;
\begin{displaymath}
\begin{split}
Ns(x) = & \thiso[Ass]^{-1}\mu\left[\sum_{i=1}^n \alpha_i(dx_i\cotimes 1)\Delta(x_{i+1})\ldots \Delta(x_n)\Delta(x_1)\ldots\Delta(x_{i-1})\right] \\
& + \thiso[Ass]^{-1}\mu\left[ \sum_{i=1}^n \alpha_i\beta_i\Delta(x_{i+1})\ldots \Delta(x_n)\Delta(x_1)\ldots\Delta(x_{i-1})(1\cotimes dx_i)\right], \\
= & \thiso[Ass]^{-1}\left[\sum_{i=1}^n \alpha_i dx_i\cdot\mu\Delta(x_{i+1}\cotimes\ldots\cotimes x_n\cotimes x_1\cotimes\ldots\cotimes x_{i-1})\right] \\
& + \thiso[Ass]^{-1}\left[ \sum_{i=1}^n \alpha_i\beta_i\mu\Delta(x_{i+1}\cotimes\ldots\cotimes x_n \cotimes x_1 \cotimes \ldots \cotimes x_{i-1})\cdot dx_i \right], \\
= & 2\sum_{i=1}^n \alpha_i x_i\cotimes\mu\Delta(x_{i+1}\cotimes\ldots\cotimes x_n\cotimes x_1\cotimes\ldots\cotimes x_{i-1}), \\
= & 2\tilde{s}N(x). \\
\end{split}
\end{displaymath}
\item[(ii)]
Let $x \in W$ and $y \in W^{\cotimes n}$ for $n \geq 0$: from equation \eqref{eqn_shucomdummya} we see that
\[ s(1-z)(xy) = s([x,y]) = [x,s(y)] = (1-z)\tilde{s}(xy). \]
\end{enumerate}
\end{proof}

\begin{cor} \label{cor_nrmshu}
Let $W$ be a profinite graded module and consider the maps
\[ e(j),\tilde{e}(j):\prod_{i=1}^\infty W^{\cotimes i} \to \prod_{i=1}^\infty W^{\cotimes i} \]
defined by Definition \ref{def_idemdf}:
\begin{enumerate}
\item[(i)]
For all $j \geq 0$,
\[ \tilde{e}(j) \circ N = N \circ e(j+1). \]
\item[(ii)]
For all $j \geq 1$,
\[ e(j) \circ (1-z) = (1-z) \circ \tilde{e}(j). \]
\end{enumerate}
\end{cor}

\begin{proof}
This is just a formal consequence of Lemma \ref{lem_nrmshu} and Lemma \ref{lem_shidem}. $e(j)$ and $\tilde{e}(j)$ are just the projections onto the eigenspaces
\[ \{ x \in \prod_{i=1}^\infty W^{\cotimes i} : s(x)=\lambda_j x \} \quad \text{and} \quad \{ x \in \prod_{i=1}^\infty W^{\cotimes i} : \tilde{s}(x)=\lambda_j x \} \]
respectively (where $\lambda_i:=2^i$).

Part (ii) follows immediately from part (ii) of Lemma \ref{lem_nrmshu} and part (i) follows from part (i) of Lemma \ref{lem_nrmshu} and the observation $\lambda_{j+1}=2\lambda_j$.
\end{proof}

We are now in a position to construct the Hodge decomposition of cyclic cohomology:

\begin{theorem}
Let $V$ be a \ci-algebra:
\begin{enumerate}
\item[(i)]
The cyclic Hochschild complex \cchoch{V} splits as a direct product of the subcomplexes $\Sigma\left(e(j)\left[\prod_{i=1}^\infty \left((\Sigma V^*)^{\cotimes i}\right)_{Z_i}\right]\right)$:
\[ \cchoch{V}=\prod_{j=1}^\infty\left(\Sigma\left(e(j)\left[\prod_{i=1}^\infty \left((\Sigma V^*)^{\cotimes i}\right)_{Z_i}\right]\right),b'\right). \]
\item[(ii)]
Let $Q_j$ and $\widetilde{Q}_j$ denote the modules
\[ \Sigma\left(e(j)\left[\col{V}\right]\right) \quad \text{and} \quad \Sigma\left(\tilde{e}(j)\left[\col{V}\right]\right) \]
respectively. The Tsygan bicomplex \ctsygan{V} splits as a direct product of subcomplexes,
\[ \ctsygan{V} = \prod_{j=0}^\infty \Gamma_j; \]
where $\Gamma_j$ is the subcomplex,
\[ \Gamma_j:= \xymatrix{(\widetilde{Q}_j,b) \ar^{1-z}[r] & (Q_j,b') \ar^{N}[r] & (\widetilde{Q}_{j-1},b) \ar^{1-z}[r] & (Q_{j-1},b') \ar^{N}[r] & \ldots \ar^{1-z}[r] & (Q_1,b') \ar^{N}[r] & (\widetilde{Q}_0,b)}. \]
\end{enumerate}
\end{theorem}

\begin{proof}
\
\begin{enumerate}
\item[(i)]
Part (ii) of Corollary \ref{cor_nrmshu} tells us that the maps $e(j):\col{V} \to \col{V}$ could be lifted uniquely to maps on the module of coinvariants;
\[ e(j):\prod_{i=1}^\infty \left((\Sigma V^*)^{\cotimes i}\right)_{Z_i} \to \prod_{i=1}^\infty \left((\Sigma V^*)^{\cotimes i}\right)_{Z_i}, \quad j \geq 1. \]
By Lemma \ref{lem_cyciso} we know that
\[ \cchoch{V} = \left(\Sigma\left(\prod_{i=1}^\infty \left((\Sigma V^*)^{\cotimes i}\right)_{Z_i}\right),b'\right). \]
By Lemma \ref{lem_shidem} and Lemma \ref{lem_shucom} we see that this complex splits as a direct product of subcomplexes:
\[ \left(\Sigma\left(\prod_{i=1}^\infty \left((\Sigma V^*)^{\cotimes i}\right)_{Z_i}\right),b'\right)=\prod_{j=1}^\infty\left(\Sigma\left(e(j)\left[\prod_{i=1}^\infty \left((\Sigma V^*)^{\cotimes i}\right)_{Z_i}\right]\right),b'\right). \]
\item[(ii)]
Lemma \ref{lem_shidem} tells us that the bigraded module \ctsygan{V} splits as a direct product of bigraded submodules
\[ \ctsygan{V}=\prod_{j=0}^\infty \Gamma_j. \]
Corollary \ref{cor_nrmshu} tells us that the operators $(1-z)$ and $N$ restrict to maps
\begin{align*}
& N: e(j+1)\left[\col{V}\right] \to \tilde{e}(j)\left[\col{V}\right], \\
& (1-z): \tilde{e}(j)\left[\col{V}\right] \to e(j)\left[\col{V}\right] .
\end{align*}
Combining this with Lemma \ref{lem_shucom} we conclude that each $\Gamma_j$ is actually a subcomplex of \ctsygan{V} whence the result.
\end{enumerate}
\end{proof}

\begin{rem}
Note that since the bicomplex \ctsygan{V} splits as the direct product of the $\Gamma_j$'s, where the latter are
bicomplexes located within vertical strips of \emph{finite width}, it follows that both spectral sequences associated with \ctsygan{V} converge to its cohomology.
\end{rem}

\begin{rem}
Recall that in Definiton \ref{def_connes} we defined the normalised Connes complex computing the cyclic cohomology of a unital \ai-algebra $V$. It is possible to construct a Hodge decomposition of this complex using Corollary \ref{cor_nrmshu} to obtain an extension of Theorem 4.6.7 of \cite{loday} for unital \ai-algebras.
\end{rem}

\begin{rem}
Let $q:\cchoch{V} \overset{\sim}{\to} \ctsygan{V}$ be the quasi-isomorphism defined by Lemma \ref{lem_cyciso} part (ii). Corollary \ref{cor_nrmshu} part (i) implies that $q$ respects the Hodge decomposition of these complexes, that is to say that for all $j \geq 0$ it restricts to a map,
\[ q: \left(\Sigma\left(e(j+1)\left[\prod_{i=1}^\infty \left((\Sigma V^*)^{\cotimes i}\right)_{Z_i}\right]\right),b'\right) \overset{\sim}{\to} \Gamma_j. \]
In particular, the Hodge decomposition of \cchoch{V} and the Hodge decomposition of \ctsygan{V} agree on the level of cohomology.
\end{rem}

\begin{defi}
Let $V$ be a \ci-algebra. We will define the complex \ecchoch{V}{j} as the subcomplex of \cchoch{V} given by the formula;
\[ \ecchoch{V}{j}:=\left(\Sigma\left(e(j+1)\left[\prod_{i=1}^\infty \left((\Sigma V^*)^{\cotimes i}\right)_{Z_i}\right]\right),b'\right), \quad j\geq 0 \]
and denote its cohomology by \ehchoch{V}{j}.
\end{defi}

We will now describe how the Hodge decomposition splits the long exact sequence of Proposition \ref{prop_sbiseq}:

\begin{prop} \label{prop_sbidec}
Let $V$ be a unital \ci-algebra. The Hodge decomposition of the Hochschild and cyclic Hochschild cohomologies respects the long exact sequence of Proposition \ref{prop_sbiseq}, that is to say that for all $j\geq 0$ we have the following long exact sequence in cohomology:
\[  \ldots\xymatrix{ \ehchoch[n-2]{V}{j} \ar^{S}[r] & \ehchoch[n]{V}{j+1} \ar^{I}[r] & \ehhoch[n]{V}{V^*}{j+1} \ar^{B}[r] & \ehchoch[n-1]{V}{j} \\ }\ldots. \]
\end{prop}

\begin{proof}
It is a simple check using the definitions of the maps $I$ and $S$ (cf. Remark \ref{rem_sbiseq}) and the Hodge decompositions to see that $I$ and $S$ respect the Hodge decomposition as claimed. It only remains to prove that the map $B$ restricts to a map
\[ B: \ehhoch[n]{V}{V^*}{j+1} \to \ehchoch[n-1]{V}{j}. \]

By Proposition \ref{prop_htpret}, every cocycle in \choch{V}{V^*} is cohomologous to a normalised cocycle. Let $x$ be such a normalised cocycle, i.e. $x \in \Sigma V^* \cotimes \widehat{T}(\Sigma V/\gf)^*$ and $b(x)=0$: it follows from figure \eqref{eqn_sbiseq} that
\[ B(x) = -h(1-z)[x] = -h(x). \]
Clearly $\tilde{e}(j+1)[x]$ is still normalised so;
\begin{displaymath}
\begin{split}
B(\tilde{e}(j+1)[x]) & = -h(\tilde{e}(j+1)[x]), \\
& = -e(j+1)h(x); \\
\end{split}
\end{displaymath}
hence the image of $B$ lies in \ehchoch[n-1]{V}{j} as claimed.
\end{proof}

Let $V$ be a \ci-algebra and consider the map $l:\drzf[Lie]{\clalg{V}} \to \drzf[Ass]{\ctalg{V}}$ defined by equation \eqref{eqn_liasmp}. By Lemma \ref{lem_plcomm} this is a map of complexes,
\[ l: \ccaq{V} \to \cchoch{V}. \]
We shall now show that this map is in fact split:

\begin{prop} \label{prop_caqspt}
Let $V$ be a \ci-algebra. There is an isomorphism
\[ \ccaq{V} \cong \ecchoch{V}{1}. \]
In particular the map $l: \ccaq{V} \to \cchoch{V}$ is split.
\end{prop}

\begin{proof}
We shall establish the following identity:
\begin{equation} \label{eqn_caqsptdummy}
l\left[\drzf[Lie]{\clalg{V}}\right] = e(2)\left[\prod_{i=1}^\infty \left((\Sigma V^*)^{\cotimes i}\right)_{Z_i}\right].
\end{equation}
The proposition will then follow as a result.

Consider the following commutative diagram:
\begin{displaymath}
\xymatrix{
& 0 & 0 & 0 & \\
0 \ar[r] & \frac{\prod_{i=1}^\infty \left((\Sigma V^*)^{\cotimes i}\right)_{Z_i}}{e(2)\left[\prod_{i=1}^\infty \left((\Sigma V^*)^{\cotimes i}\right)_{Z_i}\right]} \ar[u] \ar^{N}[r] & \frac{\col{V}}{\tilde{e}(1)\left[\col{V}\right]} \ar[u] \ar[r] & \frac{\col{V}}{N \cdot \col{V} + \tilde{e}(1)\left[\col{V}\right]} \ar[u] \ar[r] & 0 \\
0 \ar[r] & \prod_{i=1}^\infty \left((\Sigma V^*)^{\cotimes i}\right)_{Z_i} \ar^{d}[r] \ar[u] & \drof[Ass]{\ctalg{V}} \ar^{\overline{\thiso[Ass]^{-1}}}[u] \ar[r] & \frac{\drof[Ass]{\ctalg{V}}}{d\left(\drzf[Ass]{\ctalg{V}}\right)} \ar^{\overline{\thiso[Ass]^{-1}}}[u] \ar[r] & 0 \\
0 \ar[r] & \drzf[Lie]{\clalg{V}} \ar^{l}[u] \ar^{d}[r] & \drof[Lie]{\clalg{V}} \ar^{l}[u] \ar[r] & \frac{\drof[Lie]{\clalg{V}}}{d\left(\drzf[Lie]{\clalg{V}}\right)} \ar^{l}[u] \ar[r] & 0 \\
& 0 \ar[u] & 0 \ar[u] & 0 \ar[u] & \\
}
\end{displaymath}
where $\overline{\thiso[Ass]^{-1}}$ denotes the map induced by $\thiso[Ass]^{-1}$ by composing it with the relevant projection. It follows from Lemma \ref{lem_shidem} and Corollary \ref{cor_nrmshu} that the top row is exact and it follows from Lemma \ref{lem_poinca} that the middle and bottom rows are also exact. Lemma \ref{lem_nrmdif} and Lemma \ref{lem_clsdtf} imply that the right column is exact whilst equation \eqref{eqn_liespl} implies that the middle column is also exact. It follows from the $3\times 3$-Lemma that the left column is exact which in turn implies \eqref{eqn_caqsptdummy} and this proposition.
\end{proof}

Let $V$ be a \ci-algebra. By Remark \ref{rem_infdef} it has the structure of a complex given by a differential
\[ \check{m}_1:V \to V \]
and of course the dual $V^*$ has the structure of a complex given by the dual of the map $\check{m}_1$. We will denote the cohomology of this complex by $H^\bullet(V^*)$. We will now use the results of this section to prove the following important result:

\begin{lemma} \label{lem_aqlesc}
Let $V$ be a \ci-algebra. We have the following long exact sequence in cohomology:
\[ \xymatrix{ \ldots \ar^-{B}[r] & H^{n-2}(V^*) \ar^{S}[r] & \hcaq[n]{V} \ar^-{I}[r] & \haq[n]{V}{V^*} \ar^-{B}[r] & H^{n-1}(V^*) \ar^-{S}[r] & \ldots \\ }. \]
\end{lemma}

\begin{proof}
This follows from setting $i=0$ in the long exact sequence of Proposition \ref{prop_sbidec}. Proposition \ref{prop_aqsplt} allows us to identify \ehhoch{V}{V^*}{1} with \haq{V}{V^*} whilst Proposition \ref{prop_caqspt} allows us to identify \ehchoch{V}{1} with \hcaq{V}.

It follows from Lemma \ref{lem_nrmdif} part (i) and equation  \eqref{eqn_idelie} that
\[ \ecchoch{V}{0} = \left(\Sigma\left(e(1)\left[\prod_{i=1}^\infty \left((\Sigma V^*)^{\cotimes i}\right)_{Z_i}\right]\right),b'\right) = \left( V^* , m_1 \right); \]
where $m_1$ is the linear part of the \ci-structure on $V$. The result now follows.
\end{proof}

\begin{rem} \label{rem_ordgrd}
Suppose now that $V$ is a strictly graded commutative algebra (in which case it is a \ci-algebra). In this case there is a bigrading on \caq{V}{V^*} and \ccaq{V}. We say a 0-form $\alpha \in \ccaq{V}$ has bidegree $(i,j)$ if it is a 0-form of order $i$ and has degree $j$ as an element in the profinite graded module $\Sigma\drzf[Lie]{\clalg{V}}$. Similarly we say a 1-form $\alpha \in \caq{V}{V^*}$ has bidegree $(i,j)$ if it is a 1-form of order $i$ and has degree $j$ as an element in the profinite graded module $\Sigma\drof[Lie]{\clalg{V}}$. The differentials on \ccaq{V} and \caq{V}{V^*} both have bidegree $(1,1)$. In this situation we can formulate and prove the following corollary:
\end{rem}

\begin{cor} \label{cor_caqiso}
Let $V$ be a unital strictly graded commutative algebra. The map
\[ I:\hcaq[i+1,j]{V} \to \haq[ij]{V}{V^*} \]
of Lemma \ref{lem_aqlesc} is;
\begin{enumerate}
\item[(i)]
a monomorphism if $i=1$,
\item[(ii)]
an epimorphism if $i=2$,
\item[(iii)]
an isomorphism if $i \geq 3$.
\end{enumerate}
\end{cor}

\begin{proof}
A straightforward check utilising the definitions shows that the long exact sequence of Proposition \ref{prop_sbiseq} respects the bigrading on the cohomology referred to above, that is to say we have a long exact sequence in cohomology
\[ \xymatrix{ \ldots \ar^-{B}[r] & \hchoch[i-1,j-2]{V} \ar^{S}[r] & \hchoch[i+1,j]{V} \ar^{I}[r] & \hhoch[ij]{V}{V^*} \ar^{B}[r] & \hchoch[i,j-1]{V} \ar^-{S}[r] & \ldots \\ }. \]

Since the long exact sequence of Lemma \ref{lem_aqlesc} is derived from this long exact sequence we obtain a bigraded long exact sequence
\[ \xymatrix{ \ldots \ar^-{B}[r] & H^{i-1,j-2}(V^*) \ar^{S}[r] & \hcaq[i+1,j]{V} \ar^-{I}[r] & \haq[ij]{V}{V^*} \ar^-{B}[r] & H^{i,j-1}(V^*) \ar^-{S}[r] & \ldots \\ }. \]
Since $H^{\bullet\bullet}(V^*)$ is obviously concentrated in bidegree $(1,\bullet)$ the map
\[ B: \haq[ij]{V}{V^*} \to H^{i,j-1}(V^*) \]
is zero for all $j \in \mathbb{Z}$ and $i \neq 1$ and the map
\[ S: H^{i-1,j-2}(V^*) \to \hcaq[i+1,j]{V} \]
is zero for all $j \in \mathbb{Z}$ and $i \neq 2$, whence the result.
\end{proof}

\subsection{Hodge decomposition and normalised cyclic cohomology}
 Similar results also hold in the context of  normalised  cyclic  cohomology. Let $V$ be a unital $\ci$-algebra and consider the long exact sequence \eqref{eqn_nrmles} of Lemma \ref{lem_nrmles}. It is clear that the operator $s$ on $\cchoch{V}$ descends to an operator on the normalised cochains $\ccnhoch{V}$, hence we obtain a Hodge decomposition
\[\ccnhoch{V}=\prod_{j=1}^\infty\left(\Sigma\left(e(j)\left[\prod_{i=1}^\infty \left((\Sigma \bar{V}^*)^{\cotimes i}\right)_{Z_i}\right]\right),b'\right).\]
of $\ccnhoch{V}$ which is compatible with the Hodge decomposition of $\cchoch{V}$ under the inclusion. As before, we define the complex $\eccnhoch{V}{j}$ by the formula
\[\eccnhoch{V}{j}=\left(\Sigma\left(e(j+1)\left[\prod_{i=1}^\infty \left((\Sigma \bar{V}^*)^{\cotimes i}\right)_{Z_i}\right]\right),b'\right),\quad j\geq 0.\]
We have the following simple lemma.

\begin{lemma} \label{lem_uniresp}
Let $V$ be a non-negatively graded minimal unital $\ci$-algebra. The Hodge decomposition of cyclic Hochschild cohomology respects the long exact sequence \eqref{eqn_nrmles} of Lemma \ref{lem_nrmles}; in fact, for all $n\in\mathbb{Z}$ and $j\geq 0$ we have the following short exact sequence in cohomology:
\[\xymatrix{ 0 \ar[r] & \ehcnhoch[n]{V}{j} \ar[r]^{\iota} & \ehchoch[n]{V}{j} \ar[r]^{\pi} & \ehchoch[n]{\gf}{j} \ar[r] & 0 }\]
\end{lemma}

\begin{proof}
It is trivial to observe that the maps $i$ and $\pi$ preserve the respective summands of the Hodge decomposition. Since $V$ is a $\ci$-algebra, it follows from Lemma \ref{lem_nrmles} that the connecting map $\partial$ is equal to zero, therefore our long exact sequence degenerates to a series of short exact sequences.
\end{proof}

Now we introduce some of the corresponding terminology for cyclic Harrison cohomology. Let $V$ be a unital $\ci$-algebra. Just as in Lemma \ref{lem_nrmles}, the subspace of $\ccaq{V}$ given by
\[ \ccnaq{V} := \Sigma\left[\drzf[Lie]{\widehat{L}[\Sigma\bar{V}^*]}\right] \]
forms a subcomplex of $\ccaq{V}$, which we refer to as the subcomplex of normalised Lie 0-forms. Just as in equation \eqref{eqn_caqsptdummy}, the map $l$ from $\ccnaq{V}$ to $\ccnhoch{V}$ provides an isomorphism of complexes:
\[\ccnaq{V}\cong\eccnhoch{V}{1}.\]

Now suppose that $V$ is a non-negatively graded unital strictly graded commutative algebra and again introduce a bigrading on $\ccaq{V}$ and $\ccnaq{V}$ as in Remark \ref{rem_ordgrd}. The previous lemma yields the following proposition.

\begin{prop}
The map
\[\iota:\hcnaq[ij]{V}\to\hcaq[ij]{V}\]
induced by the canonical inclusion of $\ccnaq{V}$ into $\ccaq{V}$ is an isomorphism for all $(i,j)\neq (3,2)$.
\end{prop}

\begin{proof}
Clearly the maps $\iota$ and $\pi$ of Lemma \ref{lem_uniresp} preserve the bigrading of cyclic Hochschild cohomology, therefore setting $j=1$ in the corresponding short exact sequence gives us the short exact sequence:
\[\xymatrix{ 0 \ar[r] & \hcnaq[ij]{V} \ar[r]^{\iota} & \hcaq[ij]{V} \ar[r]^{\pi} & \hcaq[ij]{\gf} \ar[r] & 0 }.\]
To prove the proposition we must make an elementary calculation of the cyclic Harrison cohomology of the field. First of all, we note that the identity
\[\clalg{\gf} = \Sigma\gf^*\oplus(\Sigma\gf^*\otimes\Sigma\gf^*)\]
holds, since all brackets of order $\geq 3$ must vanish. From this it follows that
\[\drzf[Lie]{\clalg{\gf}} = (\Sigma\gf^*)^{\otimes 3};\]
i.e. there is one copy of $\gf$ in degree = 2, therefore $\hcaq[ij]{\gf}=\gf$ for $(i,j)=(3,2)$ and zero otherwise. The proposition now follows immediately from the above short exact sequence.
\end{proof}
\appendix
\section{Formal $\gf$-algebras} \label{app_todual}

We felt it necessary to include a section dealing with our formal passage to the dual language of formal $\gf$-algebras that we use in our paper since it is difficult to find any references for the material in the literature.

We will denote the category of free graded $\gf$-modules by $Mod_\gf$.  In what follows we shall omit the adjective `graded' when talking about graded $\gf$-algebras and graded $\gf$-modules since the ungraded ones are not considered.

\begin{defi} \label{def_promod}
A profinite $\gf$-module $V$ is a $\gf$-module which is an inverse limit of a diagram of free $\gf$-modules of finite rank; $V=\inlim{\alpha}{V_\alpha}$. A fundamental system of neighbourhoods of zero of $V$ is generated by the kernels of the projections $V\rightarrow V_\alpha$. The induced topology is known as the inverse limit topology. Profinite $\gf$-modules form a category $\pvect$ in which the morphisms are \emph{continuous} $\gf$-linear maps between $\gf$-modules.
\end{defi}

\begin{rem}
A profinite $\gf$-module can also be defined as a \emph{topologically free} $\gf$-module, i.e. a module $M$ for which there exists a collection  of elements $\{t_i\}_{i\in I}\subset M$ (topological basis) such that any element could be uniquely represented as a (possibly uncountably infinite) linear combination of the $t_i$'s.
\end{rem}

\begin{rem}
When we talk about a submodule of a profinite module generated by a set $X$ we mean the module generated by all convergent infinite linear combinations of elements in the set $X$. Note also that a submodule as well as a quotient of a profinite $\gf$-module need not be a profinite $\gf$-module in general. However in all cases that we encounter the $\gf$-modules and their submodules are obtained from profinite $\mathbb{Q}$-vector spaces by extending the scalars to $\gf$ and so this complication never arises.
\end{rem}

We would like to define the appropriate notion of tensor product in the category $\pvect$. This is given by the following definition:

\begin{defi}
Let $V=\inlim{\alpha}{V_\alpha}$ and $U=\inlim{\beta}{U_\beta}$ be profinite $\gf$-modules, then their completed tensor product $V \cotimes U$ is given by the formula,
\[ V \cotimes U:= \inlim{\alpha ,\beta}{V_\alpha \otimes V_\beta}. \]
\end{defi}

\begin{rem}
The completed tensor product could also be introduced as the universal object solving the problem of factorising continuous bilinear forms. Given two continuous linear maps between profinite modules $\phi:V \to W$ and $\psi:U \to X$ we could form the continuous linear map $\phi \cotimes \psi:V\cotimes U \to W \cotimes X$ in an obvious way. With this definition the category {\pvect} becomes a \emph{symmetric monoidal} category. If we endow {\pvect} with the direct product $\Pi$ then it also becomes an additive category.

Note that the construction of $V \cotimes U$ is not canonical as it depends on a preferred system $V_\alpha$, $U_\beta$ or rather a preferred choice of topological basis. Obviously different choices will produce continuously linearly isomorphic profinite modules.
\end{rem}

\begin{prop} \label{prop_antieq}
The functor $F:Mod_\gf \to \pvect$ given by sending $V$ to its linear dual $F(V):=\Hom_\gf(V,\gf)$ establishes an anti-equivalence of additive symmetric monoidal categories whose inverse functor $G$ is given by sending $U$ to the continuous linear dual $G(U):=\Hom_\mathrm{cont}(U,\gf)$.
\end{prop}
\noproof

\begin{rem}
Since a free $\gf$-module is a direct limit (union) of its finite rank free submodules, it is easy to see that the anti-equivalence of proposition \ref{prop_antieq} identifies the tensor product in $Mod_\gf$ with the completed tensor product in \pvect.
\end{rem}

We shall now discuss the notions of \emph{formal} associative, commutative and Lie algebras. To treat them uniformly we consider \emph{nonunital} associative and commutative $\gf$-algebras. When we say `$\gf$-algebra', that will simply mean that the corresponding statement could be applied to either commutative, associative or Lie algebras.

\begin{defi}
Let $A$ be a $\gf$-algebra whose underlying $\gf$-module is free of finite rank. For technical convenience we assume that $A$ is obtained from $\mathbb{Q}$ (or any subfield of $\gf$) by the extension of scalars. We say that $A$ is \emph{nilpotent} if there exists a positive integer $N$ for which the product (or bracket in the Lie case) of any $N$ elements in $A$ is zero.
An inverse limit of nilpotent $\gf$-algebras is called a formal $\gf$-algebra.
A morphism between two formal $\gf$-algebras is simply a continuous homomorphism of the corresponding structures. The categories of formal associative, commutative and Lie $\gf$-algebras will be denoted by $\mathcal{F}_{Ass}Alg, \mathcal{F}_{Com}Alg, \mathcal{F}_{Lie}Alg$ respectively. The notation $\mathcal{F}Alg$ will be used to denote either of the three categories.  A formal $\gf$-algebra supplied with a \emph{continuous} differential will be called a formal differential graded $\gf$-algebra.
\end{defi}

\begin{rem}
Every nonunital $\gf$-algebra (commutative or associative) gives rise to a unital one obtained by the well known procedure of adjoining a unit. We will call a formal associative or commutative algebra with an adjoined unit a \emph{formal augmented $\gf$-algebra}. In the main text by a `formal commutative or associative algebra' we will always mean a `formal commutative or associative augmented $\gf$-algebra'.
\end{rem}

An important example of a formal $\gf$-algebra is the so-called \emph{pro-free} $\gf$-algebra. To start, we define the \emph{free} associative $\gf$-algebra on a free $\gf$-module $V$ as $T_+V:=\bigoplus_{i=1}^\infty V^{\otimes i}$. Similarly the free commutative algebra on $V$ is defined as $S_+V:=\bigoplus_{i=1}^\infty (V^{\otimes i})_{S_i}$. Here we used the subscript $+$ to avoid confusion with the unital free associative algebra $TV:=\bigoplus_{i=0}^\infty V^{\otimes i}$ and similarly in the commutative case. The free Lie algebra on $V$ could be defined as the submodule in $T_+V$ spanned by all Lie monomials in $V$.

\begin{defi}
Let $V$ be a free $\gf$-module of finite rank.  Then the pro-free (associative, commutative or Lie) $\gf$-algebra on $V$ is the $\gf$-algebra formed by formal power series (associative, commutative or Lie) in elements of $V$. It will be denoted  by $\widehat{T}_+V, \widehat{S}_+V, \widehat{L}V$ in the associative, commutative and Lie cases respectively.

If $V=\inlim{\alpha}{V_\alpha}$ is a profinite $\gf$-module then we define $\widehat{T}_+V$ as $\inlim{\alpha} {\widehat{T}_+V_{\alpha}}$ and similarly in the commutative and Lie cases. We will denote by $\widehat{T}V$ and $ \widehat{S}V$ the unital versions of $\widehat{T}_+V$ and $\widehat{S}_+V$ respectively.
 \end{defi}

Clearly pro-free $\gf$-algebras are formal $\gf$-algebras. Furthermore we have the following result whose proof is a simple check.

\begin{prop}
 The functor $V\mapsto F(V)$ from  $\pvect$ to $\mathcal{F}Alg$ is left adjoint to the forgetful functor.
\end{prop}\noproof

\begin{rem}
The category of formal $\gf$-algebras is equivalent to the category of \emph{cocomplete} $\gf$-coalgebras,
i.e. $\gf$-coalgebras $C$ for which the kernels of iterated coproducts $\Delta^n:C\rightarrow C^{\otimes n}$ form an exhaustive filtration.  The functor from cocomplete $\gf$-coalgebras to formal $\gf$-algebras is simply taking the $\gf$-linear dual, and the
inverse functor is taking the \emph{continuous} dual. Because of this equivalence the theory of $\infty$-algebras is often formulated in terms of coalgebras and coderivations.
\end{rem}

We will often consider vector fields (=continuous derivations) and diffeomorphisms (=continuous invertible homomorphisms) of formal $\gf$-algebras. The following proposition is straightforward:

\begin{prop} \label{prop_dergen}
Let $\widehat{T}V$ ($\widehat{S}V$, $\widehat{L}V$) be a pro-free associative (commutative or Lie) algebra on a profinite $\gf$-module $V$.
Then any vector field or diffeomorphism of $\widehat{T}V$ ($\widehat{S}V$, $\widehat{L}V$) is uniquely determined by its restriction on $V$. In particular, any vector field $\xi$ has the form
\[\xi=\sum f_i(\bf t)\partial_{t_i},\]
where $\boldsymbol{t}:=\{t_i\}_{i \in I}$ is a topological basis of $V$ and $f_i(\boldsymbol{t})$ is a formal power series (associative, commutative or Lie) in the $t_i$'s.
\end{prop}\noproof

We will now define the formal universal enveloping algebra  of a formal Lie algebra $\mathcal{G}$.

\begin{defi}
\
\begin{enumerate}
\item
First let $\mathcal{G}$ be a nilpotent Lie $\gf$-algebra whose underlying $\gf$-module is free of finite rank. Denote by ${\mathcal{U}}\mathcal{G}$ its usual universal enveloping algebra and by $I(\mathcal{G})$ the augmentation ideal in ${\mathcal{U}}\mathcal{G}$. Then  $\overline{\mathcal{U}}\mathcal{G}:=\inlim{n}I(\mathcal{G})/I^n(\mathcal{G})$.
\item
Now let ${\mathcal{G}}=\inlim{} {\mathcal{G}}_i$ be a formal Lie $\gf$-algebra. Here $\{{\mathcal{G}}_i\}_{i\in i}$ is an inverse system of nilpotent Lie $\gf$-algebras. Then $\overline{\mathcal{U}}\mathcal{G} := \inlim{n}\overline{\mathcal{U}}\mathcal{G}_i$.
\end{enumerate}
\end{defi}

\begin{rem}
Note that  $\overline{\mathcal{U}}\mathcal{G}$ is a formal associative algebra, in particular it is nonunital. We will denote by $\widehat{\mathcal{U}}\mathcal{G}$ the unital $\gf$-algebra obtained from $\overline{\mathcal{U}}\mathcal{G}$ by adjoining a unit.
\end{rem}

Just like the usual universal enveloping algebra the formal one is characterised by a certain universal property. Namely, associated with any formal associative algebra $A$ is a formal Lie algebra $l(A)$ whose underlying $\gf$-module coincides with that of $A$ and the Lie bracket is defined as
\[ [a,b]:=ab-(-1)^{|a||b|}ba; \quad a,b\in A. \]
Thus, we have a functor $l$ from formal associative $\gf$-algebras to formal Lie $\gf$-algebras. Then we have the following result.

\begin{prop}\label{adjoint}
The functor $\overline{\mathcal{U}}:\mathcal{F}_{Lie}Alg \to \mathcal{F}_{Ass}Alg$ is left adjoint to the functor $l:\mathcal{F}_{Ass}Alg \to \mathcal{F}_{Lie}Alg$.
\end{prop}

\begin{proof}
We need to show that there is a canonical isomorphism
\[\mathcal{F}_{Ass}Alg(\overline{\mathcal{U}}\mathcal{G}, A)\cong \mathcal{F}_{Lie}Alg(\mathcal{G},l(A)).\]
First assume that $\mathcal{G}$ and $A$ are finite rank $\gf$-modules.  Choose a positive integer $N$ for which $A^N=0$.  Then we have the following isomorphisms:
\begin{equation} \label{frank}
\begin{split}
\mathcal{F}_{Ass}Alg(\overline{\mathcal{U}}\mathcal{G}, A) & \cong \mathcal{F}_{Ass}Alg(\overline{\mathcal{U}}\mathcal{G}/(\overline{\mathcal{U}}\mathcal{G})^N, A), \\
& \cong \Hom_{\gf-alg}(I(\mathcal{G}),A), \\
& \cong \mathcal{F}_{Lie}Alg(\mathcal{G},l(A)).
\end{split}
\end{equation}

Now let $\mathcal{G}=\inlim{\alpha}\mathcal{G}_\alpha$, $A=\inlim{\beta}{A}_\beta$ where $\mathcal{G}_\alpha, A_\beta$ are nilpotent Lie $\gf$-algebras and nilpotent associative $\gf$-algebras respectively.  We obtain using \eqref{frank}:
\begin{align*}
\mathcal{F}_{Ass}Alg(\overline{\mathcal{U}}\mathcal{G}, A) & \cong \inlim{\alpha}\dilim{\beta}\mathcal{F}_{Ass}Alg(\overline{\mathcal{U}}\mathcal{G}_\alpha, A_\beta), \\
& \cong \inlim{\alpha}\dilim{\beta}\mathcal{F}_{Lie}Alg(\mathcal{G}_\alpha,l(A_\beta)), \\
& \cong \mathcal{F}_{Lie}Alg(\mathcal{G},l(A)).
\end{align*}
\end{proof}

\begin{example}
Let $\mathcal{G}:=\widehat{L}V$, the pro-free Lie algebra on a profinite $\gf$-module $V$. Then $\widehat{\mathcal{U}}\mathcal{G}=\widehat{T}V$, the pro-free associative $\gf$-algebra on $V$.
\end{example}

\begin{rem} \label{rem_ccmult}
The above example is especially important for us. It is well known that for $V\in Mod_{\gf}$ the algebra $TV$ is in fact a Hopf algebra. The cocommutative coproduct $\Delta: TV\rightarrow TV\otimes TV$ (sometimes called the \emph{shuffle coproduct}) is defined uniquely by the requirement $\Delta$ be a $\gf$-algebra homomorphism and that $\Delta(v)=v\otimes 1+1\otimes v$ for all $v\in V$. Then $LV$ is naturally identified with the Lie subalgebra of primitive elements in $TV$.

These facts have formal analogues: let $V$ now be a profinite $\gf$-module. Then $\widehat{T}V$ has the structure of a (formal) Hopf algebra with the coproduct $\Delta: \widehat{T}V\rightarrow \widehat{T}V\cotimes \widehat{T}V$ which is uniquely specified by the requirement that $\Delta$ be a continuous homomorphism of formal $\gf$-algebras and that  $\Delta(v)=v\otimes 1+1\otimes v$ for all $v\in V$. Then the Lie algebra $\widehat{L}V$ is naturally identified with the Lie subalgebra of primitive elements in $TV$.
\end{rem}

Next, we need to define the notion of a \emph{formal module} over a formal algebra.

\begin{defi}
\
\begin{enumerate}
\item
Let $A$ be a formal $\gf$-algebra (commutative or associative) and $V$ be a profinite $\gf$-module. Then $V$ is said to be a left \emph{formal $A$-module} if there is a map of $\gf$-modules $A\cotimes V\rightarrow V; a\cotimes v\mapsto av$ such that the usual associativity axiom holds: $(ab)v=a(bv)$ for all $a,b\in A, v\in V$. Similarly, $V$ is a right formal $A$-module if there is a map of $\gf$-modules $V\cotimes A\rightarrow V; v\cotimes a\mapsto av$ such that $v(ab)=(va)b$ for all $a,b\in A, v\in V$. Finally, $V$ is a formal $A$-bimodule if it is both a right and left formal $A$-module and if $(av)b=a(vb)$ for   $a,b\in A, v\in V$.
\item Let $\mathcal{G}$ be a formal Lie $\gf$-algebra and $V$ be a profinite $\gf$-module. Then $V$ is said to be a  \emph{formal $\mathcal{G}$-module} if there is a map of $\gf$-modules $\mathcal{G}\cotimes V\rightarrow V; g\cotimes v\mapsto gv$ such that
$[g,h]v=g(hv)-h(gv)$.
\end{enumerate}
\end{defi}

Note that Proposition \ref{adjoint} implies that the structure of a formal $\mathcal{G}$-module on $V$ is equivalent to the structure of a formal $\widehat{\mathcal{U}}\mathcal{G}$-module on $V$.

\end{document}